\documentclass[11pt,a4paper]{amsart}

\usepackage{hyperref,upref} 
\hypersetup{citebordercolor=.1 .6 1}

\IfFileExists{mmymtpro2.sty}{\usepackage[mtpccal,mtpscr,subscriptcorrection]{mymtpro2}
	\let\ssave\mathcal
  \let\mathcal\mathscr
  \let\mathscr\ssave
  \let\mathocal\mathscr
  \newcommand{\hXcP}{\what{X\mkern2mu}\mkern-3mu_{\cP}}
  \newcommand{\hXcG}{\what{X\mkern2mu}\mkern-3mu_{\cG}}
  \def\cup{\cupprod}
  \def\cap{\capprod}
  \def\bigcup{\bigcupprod}
  \def\bigcap{\bigcapprod}
  \def\bigcupdisjoint{\mathop{\kern10pt\raisebox{4pt}{$\cdot$}\kern-12pt\bigcup}\limits}
}%
	{\usepackage{amssymb}%
}

\numberwithin{equation}{section}
\newtheoremstyle{ttheorem}%
       {1.8ex\@plus1ex}                % space above
       {2.1ex\@plus1ex\@minus.5ex}      % space below
       {\itshape}           % body font
       {0pt}                   % indent amount          
       {\bfseries}          % Theoremhead font 
       {.}                  % Punctuation after theorem head
       {.5em}               % Space after theorem head
       {}                % Theorem head spec (can be left empty: `normal')

\newtheoremstyle{ddefinition}%
       {1.8ex\@plus1ex}                % space above
       {2.1ex\@plus1ex\@minus.5ex}      % space belowi 
       {}           % body font
       {0pt}                   % indent amount
       {\bfseries}           % Theoremhead font 
       {.}                  % Punctuation after theorem head
       {.5em}               % Space after theorem head
       {}                % Theorem head spec (can be left empty: `normal')

\newtheoremstyle{rremark}%
       {1.8ex\@plus1ex}                % space above
       {2.1ex\@plus1ex\@minus.5ex}      % space belowi 
       {\normalfont}        % body font
       {0pt}                   % indent amount 
       {\bfseries}           % Theoremhead font 
       {.}                  % Punctuation after theorem head
       {.5em}               % Space after theorem head
       {}                   % Theorem head spec (can be left empty: `normal')

\theoremstyle{ttheorem}
\newtheorem{theorem}{Theorem}[section]
\newtheorem{lemma}[theorem]{Lemma}
\newtheorem{proposition}[theorem]{Proposition}
\newtheorem{cor}[theorem]{Corollary}
\newtheorem*{auxlemma}{Auxiliary Lemma}

\theoremstyle{ddefinition}
\newtheorem{definition}[theorem]{Definition}
\newtheorem{ass}[theorem]{Assumption}

\theoremstyle{rremark}
\newtheorem{remark}[theorem]{Remark}
\newtheorem{myremarks}[theorem]{Remarks}
\newtheorem{myexamples}[theorem]{Examples}
\newtheorem{example}[theorem]{Example}

\newenvironment{remarks}{\begin{myremarks}\begin{nummer}}%
    {\end{nummer}\end{myremarks}}
    {\end{nummer}\end{myexamples}}

\newcounter{numcount}
\newcommand{\labelnummer}{(\roman{numcount})}%

\makeatletter
\providecommand{\showkeyslabelformat}[1]{\relax}        %	 for compatibility with package showkeys
\let\mysaveformat\showkeyslabelformat                   %
\def\myformat#1{\raisebox{-1.5ex}{\mysaveformat{#1}}}   %

\newenvironment{nummer}%
  {\let\curlabelspeicher\@currentlabel%
    \begin{list}{\textup{\labelnummer}}%
      {\usecounter{numcount}\leftmargin0pt%
        \topsep0.5ex\partopsep2ex\parsep0pt\itemsep0ex\@plus1\p@%
        \labelwidth2.5em\itemindent3.5em\labelsep1em%
      }%
    \let\saveitem\item%
    \def\item{\saveitem%
      \def\@currentlabel{\curlabelspeicher\kern.1em\labelnummer}}%
    \let\savelabel\label%
    \def\label##1{{\ifnum\thenumcount=1\let\showkeyslabelformat\myformat\fi\savelabel{##1}}%
										{\def\@currentlabel{\labelnummer}% 
									 	\let\showkeyslabelformat\@gobble%	 for compatibility with package showkeys
									 	\savelabel{##1item}%
										}%
	   							}%
  }{\end{list}}%
  
\newenvironment{indentnummer}%
  {\let\curlabelspeicher\@currentlabel%
    \begin{list}{\textup{\labelnummer}}%
      {\usecounter{numcount}\leftmargin0pt%
        \topsep0.5ex\partopsep2ex\parsep0pt\itemsep0ex\@plus1\p@%
        \labelwidth2.5em\itemindent0em\labelsep1em%
        \leftmargin2.5em}%
    \let\saveitem\item%
    \def\item{\saveitem%
      \def\@currentlabel{\curlabelspeicher\kern.1em\labelnummer}}%
    \let\savelabel\label%
    \def\label##1{{\ifnum\thenumcount=1\let\showkeyslabelformat\myformat\fi\savelabel{##1}}%
										{\def\@currentlabel{\labelnummer}% 
									 	\let\showkeyslabelformat\@gobble%	 	 for compatibility with package showkeys
									 	\savelabel{##1item}%
										}%
    							}%
  }{\end{list}}%
  
\def\itemref#1{\ref{#1item}}

\def\section{\@startsection{section}{1}%
  \z@{1.3\linespacing\@plus\linespacing}{.5\linespacing}%
  {\normalfont\scshape\centering}}

\def\subsection{\@startsection{subsection}{2}%
  \z@{.8\linespacing\@plus.5\linespacing}{-1em}%
  {\normalfont\bfseries}}

\def\nlsubsection{\@startsection{subsection}{2}%
  \z@{.8\linespacing\@plus.5\linespacing}{.1ex}%
  {\normalfont\bfseries}}

\newcommand{\cA}{\mathcal{A}}

\newcommand{\cC}{\mathcal{C}}
\newcommand{\cE}{\mathcal{E}}
\newcommand{\cF}{\mathcal{F}}
\newcommand{\cG}{\mathcal{G}}
\newcommand{\cP}{\mathcal{P}}
\newcommand{\cQ}{\mathcal{Q}}
\newcommand{\cS}{\mathcal{S}}
\newcommand{\cV}{\mathcal{V}}

\renewcommand{\AA}{\mathbb{A}}

\newcommand{\NN}{\mathbb{N}}
\newcommand{\PP}{\mathbb{P}}

\newcommand{\RR}{\mathbb{R}}
\newcommand{\VV}{\mathbb{V}}

\renewcommand{\le}{\leqslant}
\renewcommand{\ge}{\geqslant}
\renewcommand{\d}{\mathrm{d}}

\newcommand{\hatmu}{\what{\mu}}

\DeclareMathOperator{\supp}{supp}

\DeclareMathOperator{\vol}{vol}
\DeclareMathOperator{\diam}{diam}
\DeclareMathOperator{\card}{card}
\DeclareMathOperator{\sydi}{\raisebox{.6pt}{\scriptsize$\bigtriangleup$}}
\def\bigtimes{\mathop{\raisebox{-1.5pt}{\LARGE $\times$}}\limits}
\def\mybigotimes{\mathop{\raisebox{-1.5pt}{\LARGE $\otimes$}}\limits}

\providecommand{\comp}{\circ}
\providecommand{\mathbold}[1]{\mathbf{#1}}
\providecommand{\wtilde}[1]{\widetilde{#1}}
\providecommand{\what}[1]{\widehat{#1}}
\providecommand{\hXcG}{\what{X}_{\cG}}
\providecommand{\hXcP}{\what{X}_{\cP}}
\providecommand{\mathocal}[1]{#1}
\providecommand{\bigcupdisjoint}{\mathop{\kern7pt\raisebox{6pt}{$\cdot$}\kern-9.5pt\bigcup}\limits}

\makeatletter

\def\per{.}
\def\HarvardComma{}

\newcounter{aucount}
\newif\ifedplural
\newif\ifper\pertrue

\def\au#1#2{{#1 #2}}
\def\lau#1#2{{#1 #2},}
\def\ed#1#2{\ifnum\theaucount=0\relax\fi{#1 #2}\addtocounter{aucount}{1}}
\def\led#1#2{\ifnum\theaucount=0\relax\edpluralfalse\else\edpluraltrue\fi{#1
    #2} (\editorname.),\setcounter{aucount}{0}}
\def\editorname{\ifedplural Eds\else Ed\fi}

\def\et{\ifnum\theaucount=1\else\HarvardComma\fi{} and\ }
\def\ti#1{\emph{#1},\ifper\fi\pertrue}

\def\bti{\@ifnextchar[\bbti\bbbti}
\def\bbti[#1]#2{\emph{#2}, #1,}
\def\bbbti#1{\emph{#1},}

\def\z{\@ifnextchar[\zz\zzz}
\def\zz[#1]#2#3#4#5{\perfalse{#2} \textbf{#3}, #4 \ifx
  @#5@\relax\else (#5)\fi [#1]\ifper\per\fi\pertrue} 
\def\zzz#1#2#3#4{{#1} \textbf{#2}, #3 \ifx @#4@\relax\else
  (#4)\fi\ifper\per\fi\pertrue}

\def\pub{\@ifstar\pubstar\pubnostar}
\def\pubnostar{\@ifnextchar[\@@pubnostar\@pubnostar}
\def\@@pubnostar[#1]#2#3#4{#2, #3, #4, #1\ifper\per\fi\pertrue}
\def\@pubnostar#1#2#3{#1, #2, #3\ifper\per\fi\pertrue}
\def\pubstar[#1]#2#3#4{\perfalse #2, #3, #4 [#1]\per\pertrue}

\makeatother

\begin{document}

%%%%%%%%%%%%%%%%%%%%%%%%%%%%%%%%%%%%%%%%%%%%%%%%%%%%%%%%%%%%%%
%%%%%%%%%%%%%%%%%%%%%%%%%%%%%%%%%%%%%%%%%%%%%%%%%%%%%%%%%%%%%%

\title[Randomly coloured point sets]{Ergodic properties of randomly coloured
  point sets\,\footnotemark{\raisebox{1mm}{*}}}

\author[P.\ M\"uller]{Peter M\"uller}
\address{Mathematisches Institut,
  Ludwig-Maximilians-Universit\"at M\"unchen,
  Theresienstra\ss{e} 39,
  80333 M\"unchen, Germany}
\email{mueller@lmu.de}

\author[C.\ Richard]{Christoph Richard}        
\address{Department f\"ur Mathematik,
  Friedrich-Alexander-Universit\"at Erlangen-N\"urnberg,
  Cauerstra\ss{}e 11,
  91058 Erlangen, Germany}
\email{richard@mi.uni-erlangen.de}

\thanks{* This version is almost identical to the one published electronically on May 10, 2012 in the \emph{Canadian Journal of Mathematics}.\quad \href{http://dx.doi.org/10.4153/CJM-2012-009-7}{doi:10.4153/CJM-2012-009-7}}

\begin{abstract}
  We provide a framework for studying randomly coloured point sets in a
  locally compact, second-countable space on which a
  metrisable unimodular group acts continuously and properly. 
  We first construct and describe an
  appropriate dynamical system for uniformly discrete uncoloured point sets. For
  point sets of finite local complexity, we
  characterise ergodicity geometrically in terms of pattern frequencies. 
  The general framework allows to incorporate a random
  colouring of the point sets. We derive an ergodic theorem for randomly
  coloured point sets with finite-range dependencies. 
  Special attention is paid to the exclusion of exceptional instances for uniquely ergodic
  systems. The setup allows for a straightforward application to randomly
  coloured graphs.
\end{abstract}

\dedicatory{We dedicate this work to Hajo Leschke on the occasion of his 67th birthday}

\maketitle

%%%%%%%%%%%%%%%%%%%%%%%%%%%%%%%%%%%%%%%%%%%%%%%%%%%%%%%%%%%%%%
%%%%%%%%%%%%%%%%%%%%%%%%%%%%%%%%%%%%%%%%%%%%%%%%%%%%%%%%%%%%%%

\section{Introduction}

Delone sets are subsets of Euclidean space which are uniformly
discrete and relatively dense. In the natural sciences, they are used
to model pieces of matter. Over the last years, geometric and spectral
properties of Delone sets have been studied by many authors, using
methods from topological dynamical systems, see e.g.~\cite{RaWo92,
  Hof93, Hof95, So97, BeHeZ, LeMo02, LaPl03, LeSt03, KlLe03, LeSt05, LeS08}.  
Here, the dynamical system arises from the closure of the translation
orbit of the given Delone set with respect to a suitable topology. 
This approach is particularly useful if the dynamical
system is uniquely ergodic, since the uniform ergodic theorem can then be used to 
infer properties of the given Delone set one has started with. 
For a Delone set of finite local complexity, a 
geometric characterisation of unique ergodicity in terms
of uniform pattern frequencies appears in \cite{LeMo02}. If the Delone set is not periodic, then 
such a characterisation cannot be achieved with a discrete periodic
subgroup of the Euclidean group as the group acting on the dynamical system. 
Therefore one has to rely upon an ergodic theorem for the action of a more general group 
than the multi-dimensional integers.
  
This approach has been generalised considerably in recent years.  
Euclidean space has been replaced by a $\sigma$-compact, 
locally compact Abelian group, which admits a suitable averaging sequence, and 
on which the same group acts by translations \cite{Sch00}.  Within that setup, the 
important subclass of repetitive regular model sets, see e.g.~\cite{Mo00}, which have a pure point diffraction spectrum such as periodic point sets, could be characterised by 
certain properties of the underlying dynamical system including strict ergodicity and pure point dynamical spectrum \cite{BLM07}.
More generally, dynamical systems of translation bounded Borel measures
\cite{BaLe04, LeRi07} on such spaces have been considered. Discrete
subsets of a general locally compact topological space have been
studied in \cite{Yok05} via group actions of a locally compact group,
focussing on finite local complexity and on repetitivity. 

As we are interested in discrete geometry, our setup will be
formulated in terms of uniformly discrete sets. We will use a locally compact,
second-countable space as our basic space of points, which we will simply refer to as the point space. By choosing a metric, this allows us
to define a notion of uniform discreteness. Local compactness of the point space ensures
sufficient structure for the space of uniformly discrete point sets.
The \emph{first main goal} of this paper is to establish 
geometric criteria for (unique) ergodicity of the dynamical system associated to a collection of 
uniformly discrete point sets in 
terms of pattern frequencies. To do so, we rely on properness of the group action. In addition, we will measure 
the ``size'' of a subset of the point space in terms of the Haar measure on the group, which is pushed forward by the group action and a reference point in the point space. 
We require that size to be
the same for group-equivalent reference points. This is ensured for unimodular groups. 
Thus our setup comprises 
non-Abelian groups such as the Euclidean group. In particular, this accommodates the 
pinwheel tilings of the plane \cite{Rad94, Rad95, Rad97} and their relatives \cite{Fre08}.
Apart from \cite{Yok05}, it seems that non-Abelian groups have not been treated in our 
general context so far.

In order to define an appropriate dynamical system, we require that uniformly discrete 
point sets, which are group-equivalent, have the same radius of discreteness. This will 
be guaranteed if the metric on the point space is group-invariant. Indeed, it has recently been 
shown \cite[Thm.~1.1]{AbMaNo11} that under our assumptions on the group action, such a metric always exists among the set of all metrics that are compatible with the topology of the point space.
We will supply the space of uniformly discrete point sets (of a given radius of discreteness) with the vague topology. This ensures compactness of the relevant dynamical
systems. In \cite{Yok05}, a stronger ``local matching topology'' is
favoured instead. For a proper and transitive group action, both topologies
coincide, if the point sets are of finite local complexity. This
follows with Lemma~\ref{FLClocal}, compare also \cite{BaLe04} for the
Abelian case. 

If the point space $M$ admits a uniform structure compatible with the given topology, but is not necessarily metrisable, it is still possible to define uniform discreteness via the uniformity. In fact, a continuous and proper action of a group $T$ on $M$ gives rise to a uniform structure on $M$, which induces the topology on $M$ and has the desirable invariance properties with respect to $T$, as follows from \cite[Thm.~1.2]{AbMaNo11}. A corresponding framework could include our setup for metrisable spaces $M$, as well as the approach of \cite{Sch00, BaLe04,Yok05} for a not necessarily metrisable space $M$ as special cases. We will not strive for such a generality here.

Our structural assumptions on the group and its action are minimal in a sense. 
By properness, the group inherits local compactness and $\sigma$-compactness from the point space. 
But local compactness is needed for the existence of a (well-behaved) Haar measure, and $\sigma$-compactness is required for amenability.
The r\^ole of unimodularity has been discussed above, and Lindenstrauss' pointwise ergodic theorem \cite{Lind01}, which we rely upon to a great extent, requires metrisability of the group.
We will however not assume transitivity of the group action. This is 
motivated by our desire to describe uniformly discrete sets, coloured versions 
thereof and also graphs built from such sets -- all within the same framework. 
Here, coloured point sets of possibly infinitely many colours and also graphs will appear as point sets
in some suitably enlarged point space on which one cannot expect to have a transitive group action. 
Due to the absence of transitivity we were also prompted to free the point 
space from being a group by itself. We mention that coloured Delone sets of finite local 
complexity -- and thus with at most finitely many colours 
-- have been studied by different methods in \cite{BeHeZ,LeMo02,LeSt03}.

As our choice of spaces is also canonical in stochastics, the
connection to stochastic geometry \cite{SKM95} may be broadened.  Indeed, the setup
allows to study random colourings of a point set on a rather general
level. Ergodic properties of random colourings of a repetitive
Delone set in the Euclidean plane have already been studied by Hof
\cite{Hof98}, motivated by the problem of site percolation on
the Penrose tiling. 
His approach has been used to infer diffraction properties of
random Euclidean point sets of finite local complexity \cite{BaaZin07} with finite-range dependencies  and beyond \cite{BaBiMo10}; see also \cite{Le08} for an alternative approach. A recent extension to certain Delone sets in $\sigma$-compact, locally compact Abelian groups is the subject of \cite{AI}.
Another recent generalisation to infinite-range dependencies, based on the theory of Gibbs measures and stochastic geometry, can be found in \cite{M10}. Diffraction properties of certain non-periodic stochastic point sets are also discussed in \cite{Ku1, Ku2}, where large-deviation estimates and concentration inequalities for the finite-volume scattering measure are derived. 

We are not concerned with diffraction in this paper, however. In fact, our \emph{second main 
goal} is to provide an optimal ergodic theorem for dynamical systems of randomly coloured point sets
with finite-range stochastic dependencies.
To do so, we also pursue an idea of Hof \cite{Hof98},
who used the law of large numbers for reducing the problem to that of the 
dynamical system for the underlying uncoloured point sets. Unfortunately, Hof's approach only works for point sets of finite local complexity, and is thus also restricted to finitely many colours. 
On the other hand,  Lenz \cite{Le08} proved an ergodic theorem for randomly coloured translation bounded measures on Euclidean space without the need for finite local complexity. In combining the two approaches, we obtain an ergodic theorem for randomly coloured point sets without requiring finite local complexity. And, in contrast to \cite{Le08}, it is optimal in the sense that exceptional instances are excluded as far as possible in the case of uniquely ergodic systems and continuous functions.  

In a subsequent article, we will apply the aforementioned ergodic results to 
describe spectral properties of subcritical
percolation graphs over such general point sets, compare \cite{KM} for the periodic case. 
As we are able to treat rather general colour spaces and group actions, our approach also 
opens the possibility to study finite-range operators on uniformly discrete point sets 
with quite general internal degrees of freedom and their randomised versions, such as (random) 
Schr\"odinger operators with magnetic finite-range interactions on non-periodic point sets. 

This paper is organised as follows. In Section~\ref{psdyn}, we
first recapitulate properties of dynamical systems of uncoloured point
sets, which carry over to our more general setup. As our general group
actions have apparently not been studied before in this context, we provide
proofs for the convenience of the reader. Based on the general pointwise ergodic 
theorem of Lindenstrauss~\cite{Lind01}, we state characterisations of ergodicity in 
Theorem~\ref{abstract-erg}, which are handy to use and which we have not found in the literature in sufficient generality. The same remark applies to 
Theorem~\ref{abstract-unierg}, which is the abstract analogue for uniquely ergodic systems,
extending the well-documented case of $\mathbb Z$-actions \cite{KB37, O52, Wa82, Fur81}. 
Whereas these two theorems are only of a propaedeutic nature, 
our main results of Section~\ref{psdyn} are Theorem~\ref{uniFLC} and Proposition~\ref{freq-cond12}. 
They provide geometric characterisations of ergodicity and of unique ergodicity in terms of pattern frequencies for uniformly discrete point sets of finite local complexity.

In Section~\ref{sec:ran-col}, we construct an ergodic measure for
randomly coloured point sets and present an optimal ergodic theorem as our second main result 
in Theorem~\ref{theo:perg}. Finite local complexity is not required in this section.

The formalism developed in Sections~\ref{psdyn} and~\ref{sec:ran-col} is applied 
to randomly coloured graphs in Section~\ref{secdyn}. Proofs are provided in the remaining sections.

%%%%%%%%%%%%%%%%%%%%%%%%%%%%%%%%%%%%%%%%%%%%%%%%%%%%%%%%%%%%%%%%%%%%%%%%%%
%
\section{Dynamical systems for point sets}
\label{psdyn}
%
%%%%%%%%%%%%%%%%%%%%%%%%%%%%%%%%%%%%%%%%%%%%%%%%%%%%%%%%%%%%%%%%%%%%%%%%%%

Here, we introduce our setup, discuss the basic ergodic theorem and give a geometric characterisation of ergodic point sets in terms of pattern frequencies.

\subsection{Topology on collections of point sets} \label{top-point-sets}

For the convenience of the reader, Section~\ref{proofs-top-point-sets} contains proofs of the material that we present here. 

A \emph{point space} $M$ is a non-empty, locally compact and second-countable topological
space. Throughout this paper we stick to the convention that every locally compact topological 
space enjoys the Hausdorff property.
We recall that in locally compact topological spaces, second countability is equivalent to $\sigma$-compactness and metrisability \cite[Chap.\ IX, \S 2.9, Cor.]{Bour}.

In addition to the point space $M$ we consider a metrisable topological group $T$ with left action $\alpha:=\alpha_{M}:T\times M\to M$, $(x,m)\mapsto xm$ on $M$. Throughout this paper we will rely on

\begin{ass}
	\label{basic-comp}
	The group $T$ is non-compact and its action on $M$ is continuous and proper. Moreover, we fix a 
 	$T$-invariant proper metric $d$ on $M$ that generates the topology on $M$.
\end{ass}

\begin{remarks}\label{basic-metric}
\item A metric $d$ on $M$ is $T$-invariant, iff $d(xm,xm') = d(m,m')$ for all $x \in T$ and 
	all $m,m' \in M$. A metric is proper, iff every metric ball has a compact closure. The existence of 
	a $T$-invariant proper metric $d$ on the metrisable space $M$ that generates the topology of $M$ 
	follows from $\sigma$-compactness of $M$ and properness of the group action \cite[Thm.~4.2]{AbMaNo11}.
\item We recall that the group action is continuous, iff the map $\alpha: T\times M\to M$ is 
	continuous with respect to the product topology on $T\times M$. Properness of the group action
	means that the (continuous) map 
	\begin{displaymath}
%  	\label{propermap}
  	\wtilde{\alpha}: T\times M \rightarrow M\times M, \;(x,m) \mapsto
  	(xm,m),
	\end{displaymath}
	is proper, that is, pre-images of compact sets are compact (see \cite[Chap.\ III.4]{Bour1}).
	We refer to Lemma~\ref{proper-char} below for different characterisations of properness.    
\item \label{prop-implies1}
	All results of the present section (and their proofs) remain valid for compact groups, but 
	reduce to trivial statements. For the following sections, however, non-compactness will be crucial. 		
	It will be used in the proof of the strong law of large numbers Theorem~\ref{lem:hoflemma}.
\item \label{prop-implies2} 
	Properness of the group action and local compactness of $M$ imply that $T$ is also locally compact, see Proposition~1.3 in \cite{AbMaNo11} and subsequent comments. Furthermore, since $M$ is also $\sigma$-compact, so is $T$.
\item We require metrisability of the group $T$ in order to satisfy the hypotheses of Lindenstrauss' ergodic theorem. The latter is essential for large parts of this paper. Thus, by the previous remark and by what was recalled at the beginning of this section, the group $T$ is also second countable. 
\item Below we want to define nice $T$-orbits of uniformly discrete subsets in $M$ with a given radius 
	of uniform discreteness. It is precisely for this purpose that we work with a 
	$T$-invariant metric $d$ on $M$ that is compatible with the topology. 
\end{remarks}

For convenience we recall the following alternative characterisations of properness. We use the notation $xU := \{xm: m\in U\}$ for $x\in T$ and $U \subseteq M$ and introduce the \emph{transporter} 
\begin{equation}
	\label{transporter}
 	S_{U_{1},U_{2}} := \{x \in T: xU_{1} \cap U_{2} \neq \varnothing \} \subseteq T
\end{equation}
of subsets $U_{1}, U_{2} \subseteq M$. The following lemma does not rely on Assumption~\ref{basic-comp}.

\begin{lemma}
	\label{proper-char}
	Assume that $M$ and $T$ are both locally compact and second-countable, and that $T$ acts continuously on $M$ from the left. 
	Then the following are equivalent.
	\begin{nummer}
 	\item $T$ acts properly on $M$.
	\item \label{transporterK} 
		For every choice of compact subsets $V_{1}, V_{2} \subseteq M$, 
		the transporter $S_{V_{1},V_{2}}$ is compact in $T$. 	
	\item \label{transporterU} 
		For every choice of relatively compact subsets $U_{1}, U_{2} \subseteq M$, 
		the transporter $S_{U_{1},U_{2}}$	is relatively compact in $T$. 	
	\item Given any two sequences $(x_{n})_{n\in\NN} \subseteq T$ and $(m_{n})_{n\in\NN} \subseteq M$ such that both $(m_{n})_{n\in\NN}$ and $(x_{n}m_{n})_{n\in\NN}$ converge in $M$, then $(x_{n})_{n\in\NN}$ has a convergent subsequence in $T$.
\end{nummer}
\end{lemma}

\begin{remarks}
\item We will only use the implications $\mathrm{(i)} \Rightarrow \mathrm{(ii)}  \Rightarrow \mathrm{(iii)}$ of the lemma in the sequel. But for the purpose of completeness, we have included the proof of all implications in Section~\ref{proofs-top-point-sets}. 
\item Similar characterisations can be found, e.g.,\ in \cite[Chap.\ III.4.5, Thm.~1]{Bour1} or \cite[Def.\ and Prop.~2.3]{AbMaNo11}.
\item Condition \itemref{transporterK} of the lemma implies that the map 
	$\alpha(\cdot,m): T \rightarrow M$ is proper for every $m\in M$, since 
	$(\alpha(\cdot,m))^{-1}(V)=S_{\{m\},V}$ for every compact $V\subseteq M$.
\end{remarks}

\begin{example}
An often studied special case of our setup arises when $M$ is a topological group by itself. 
Then one may choose $T:=M$ and $\alpha$ as the group multiplication from the left. 
An important example for this special case is the Abelian group $M=\RR^{\mathsf{d}}$ for $\mathsf{d}\in\NN$, equipped with the Euclidean metric, which acts on itself canonically by translation. 
This action is transitive, free and proper, and the Euclidean metric is $T$-invariant and also proper. (Recall that a group action is said to be \emph{transitive}, iff for every $m,m'\in M$ there exists $x\in T$ such that $xm=m'$. It is \emph{free}, iff for any $x\in T$ and
any $m\in M$ the property $xm=m$ implies $x=e$, the neutral element in $T$.)
Another prime example is $M=\mathbb R^{\mathsf{d}}$ 
and $T=E(\mathsf{d})$, the Euclidean group, with the Euclidean metric. Then the canonical group action from the left is transitive and proper, but not free. The Euclidean metric is also $T$-invariant and proper. 
We note, however, that in situations relevant to us, $M$ will not be a group.
\end{example}

The open ball (with respect to the metric $d$) of radius $s>0$ about $m\in M$ is denoted by $B_s(m)$.  
A subset $P\subseteq M$ is 
\emph{uniformly discrete} with radius
$r \in ]0,\infty[$, if any open ball in $M$ of radius $r$ contains at most one
element of $P$. A subset $P\subseteq M$ is called \emph{relatively dense} with radius
$R \in]0,\infty[$, if every closed ball of radius $R$ has non-empty
intersection with $P$. If $P$ is both uniformly discrete and relatively dense, then it is called a \emph{Delone set}. The collection of all subsets of
$M$, which are uniformly discrete of radius $r$, is denoted by $\mathcal
P_r(M)$. We call every element of $\mathcal P_r(M)$ a \emph{point set}.
 Throughout this paper, the radius of uniform discreteness $r$ will be
fixed.

We define a topology on $\cP_r(M)$ by requiring certain functions on $\cP_r(M)$, which are of the form \eqref{f-phi} below, to be continuous. These functions will serve as a ``scanning device'' on a point set. Let $C_{c}(M)$ denote the set of all real
valued, continuous functions on $M$ with compact support.

\begin{definition}\label{fphi}
  To $\varphi \in C_{c}(M)$ we associate
  \begin{equation} \label{f-phi}
    f_{\varphi}: \begin{array}{ccc} \cP_{r}(M) & \rightarrow  & \RR
      \quad\phantom{p}\\ 
      P & \mapsto & f_{\varphi}(P) := \displaystyle\smash{\sum_{p\in P}}
      \,\varphi(p) \end{array}. \phantom{\rule[-4ex]{1pt}{4ex}}
  \end{equation}
The \emph{vague topology} on $\cP_r(M)$ is the weakest topology such that 
$f_\varphi$ in \eqref{f-phi} is continuous for
  every $\varphi\in C_{c}(M)$.
\end{definition}

\begin{remarks} 
\item 
	Even though the set $\cP_{r}(M)$ itself depends on the metric $d$ on $M$, the nature of the vague topology on $\cP_{r}(M)$ is solely determined by the topology on $M$. 
\item
  Particular examples of open sets in $\cP_r(M)$ are given by pre-images of open balls in 
  $\mathbb R$. For $P\in\cP_r(M)$, $\varphi\in C_{c}(M)$ and $\varepsilon >0$ we define
  the open set
  \begin{displaymath} %\label{base-set} 
    U_{\varphi,\varepsilon}(P):=
    \left\{\wtilde{P}\in\mathcal P_r(M): \big|
        f_{\varphi}(\wtilde P)-f_{\varphi}(P)\big|<\varepsilon \right\}.
  \end{displaymath}
  It is readily checked that the family obtained from finite intersections of open sets 
  $U_{\varphi,  \varepsilon}(P)$ as above forms a neighbourhood base of the vague topology.
\item
  The above neighbourhood base arises naturally
  when identifying a point set $P$ with a point measure on $M$ that
  has an atom of unit mass at each point of $P$, see e.g.\
  \cite{BeHeZ, Sch00, BaLe04,Le08}. It is from this perspective that the
  topology of Definition~\ref{fphi} appears as the vague topology on
  this space of measures. For the case where $M$ is also a group,
   \cite{BaLe04} coined
  the name \emph{local rubber topology} for the vague topology (and they defined it using transitivity
  of the canonical group action of $M$ on itself). For the particular example $M=\RR^{\mathsf{d}}$ 
  the vague 
  topology was studied in \cite{LeSt03a} under the name \emph{natural topology} and
  earlier on in \cite{BeHeZ,LaPl03}.    
\item Instead of uniformly discrete subsets of $M$, one may consider more
general \emph{locally finite} sets. These are sets $P\subseteq M$ for which
$P\cap V$ is finite for every compact set $V\subseteq M$. But the
space of locally finite sets equipped with the vague topology
is not closed. For example, a sequence of locally finite point sets may give
rise to accumulation points in $M$.
\item Local compactness and second countability of $M$ imply metrisability of
the topology on  $\mathcal P_r(M)$, see \cite[Thm.~31.5]{Bau01}. 
\end{remarks}

Convergence in the topological space $\mathcal P_r(M)$
is characterised in the following lemma.

\begin{lemma}\label{konchar}
  Fix a sequence $(P_k)_{k\in\mathbb N}\subseteq \mathcal P_r(M)$.  Then the
  following statements are equivalent.
  \begin{indentnummer}
  \item The sequence $(P_k)_{k\in \mathbb N}$ converges in $\mathcal P_r(M)$.
  \item There exists $P\in\mathcal P_r(M)$ such that for all $\varphi\in C_{c}(M)$ we have
   \begin{displaymath}
      \lim_{k\to\infty} f_{\varphi}(P_{k}) = f_{\varphi}(P).
    \end{displaymath}
  \item \label{con-ball}
  	For every $m\in M$ exactly one of the following two cases
    occurs.
    \begin{enumerate}
    \saveitem[(a)]
      For every $\varepsilon>0$ we have $P_k\cap  
      B_\varepsilon(m)\ne\varnothing$ for finally all $k\in\mathbb N$. 
    \saveitem[(b)] There exists $\varepsilon>0$ such that $P_k\cap
      B_{\varepsilon}(m)=\varnothing$  
      for finally all $k\in\mathbb N$.
    \end{enumerate}
	\item \label{con-thick} 
    There exists $P\in\mathcal P_r(M)$ such that for every compact
    set $V\subseteq M$ we have, for every $\varepsilon>0$ and finally
    all $k\in\mathbb N$, the inclusions
    \begin{displaymath}
      P_k\cap V\subseteq (P)_\varepsilon \qquad \text{and}\qquad 
      P\cap V\subseteq (P_k)_\varepsilon.
    \end{displaymath}
    Here, the ``thickened'' point set
  	$(P)_\varepsilon := \bigcup_{p\in P} B_{\varepsilon}(p)$ is the set of points in
  	$M$ lying within distance less than $\varepsilon$ to $P$.
	\end{indentnummer}
  In either case, the limit $P$ is the set of all points $m\in M$ satisfying
  \textup{\itemref{con-ball}-(a)}.  
\end{lemma}

Below we will be concerned with ergodic properties of
$\cP_{r}(M)$ as a topological dynamical system. This relies on

\begin{proposition} \label{PrMcompact} The space $\mathcal{P}_r(M)$ is compact with
  respect to the vague topology.
\end{proposition}

\begin{remarks}
\item 
	In order to give a self-contained presentation, we prove (sequential) compactness of 
	the metrisable space $\cP_r(M)$ in Section~\ref{proofs-top-point-sets}. 
	Thus, Proposition~\ref{PrMcompact} yields complete metrisability of $\cP_{r}(M)$, in other words,	
	$\cP_r(M)$ is even a Polish space.	For the more 
	general case of $M$ being only $\sigma$-compact and locally compact, compactness of $\cP_{r}(M)$
	has already been shown in \cite[Thm.~3]{BaLe04}, see also \cite[Thm.~31.2]{Bau01} and \cite{BeHeZ}.
\item
	In tiling dynamical systems, the topology on $\cP_{r}(M)$ is often characterised in terms
	of a particular metric resembling a connection to symbolic dynamics, see 
	e.g.\ \cite{RaWo92, Hof98, LeMo02}. The corresponding notion of distance 
	means that two point sets are close, if they almost agree on a large ball in the point space. 
	This can be formalised as follows.  The map 
	$\mathrm{dist}: \mathcal P_r(M)\times\mathcal P_r(M)\to\RR_{\ge0}$,
	given by  
  \begin{displaymath}
    \mathrm{dist}(P,\wtilde P) :=\min\Big\{ \tfrac{1}{\sqrt{2}}\,, 
    \inf\big\{ \varepsilon>0: 
    P \cap B_{\frac{1}{\varepsilon}}\subseteq  (\wtilde P)_\varepsilon 
    \text{~~and~~} \wtilde P \cap B_{\frac{1}{\varepsilon}} \subseteq (P)_{\varepsilon} \big\}
    \Big\},
  \end{displaymath}
  where $B_{\frac{1}{\varepsilon}} := B_{\frac{1}{\varepsilon}}(m_{o})$ for some fixed 
  reference point $m_{o} \in M$, defines a metric on $\cP_{r}(M)$.  The topology induced by the above 
  metric does not depend on the choice of reference point $m_o$ and 
  coincides with the vague topology.  (Since $\cP_r(M)$ is compact, it is complete w.r.t.~the above metric, and the metric is proper.) In this paper we prefer to work with the vague topology instead of the metric.
\end{remarks}

%%%%%%%%%%%%%%%%%%%%%%%%%%%%%%%%%%%%%%%%%%%%%%%%%%%%%%%%%%%%%%%%%%%%%%%%%%
\subsection{Ergodic theorems for group actions}
\label{subsec:action-ergodic}

Our basic workhorse will be the general ergodic theorem of Lindenstrauss \cite{Lind01}. In order to apply it we need to recall some further notions. We fix a left Haar measure on the locally compact, second-countable group $T$ and write 
$\vol (S) = \int_{S}\d x$ for this Haar measure of a Borel set $S\subseteq T$. 
Below we will also impose that the group $T$ is 
\emph{unimodular}. This is equivalent to the requirement that the Haar measure is 
inversion invariant, i.e., $\int_T f(x^{-1})\d x=\int_T f(x)\d x$ for every measurable function 
$f:T\to [0,\infty]$ (and hence also for every integrable function $f:T\to\mathbb R$). In particular, this implies $\vol(S^{-1})=\vol(S)$ for every Borel
set $S\subseteq T$, where $S^{-1}:=\{x\in T: \exists s\in S \text{~~such that~~} x=s^{-1}\}$.

Since we want to compute certain group means below, we require that
$T$ admits suitable averaging sequences. As usual, for $K\subseteq T$ we denote
by $\mathring{K}=\mathrm{int}(K)$ the interior of $K$ and by $\overline{K}$ the closure of $K$, and
for $A,B\subseteq T$ we write $AB:=\{x\in T:\exists (a,b)\in A\times B \text{~~such that~~}x=ab\}$ for the Minkowski product of $A$ and $B$.

\begin{definition}
  Let $(D_n)_{n\in\mathbb N}$ be an increasing sequence of non-empty, compact
  subsets of $T$ such that $\bigcup_{n\in\mathbb N}\mathring{D_n}=T$. 
  \begin{nummer}
  \item 
    The sequence $(D_n)_{n\in\mathbb N}$ is called a
    \emph{F{\o}lner sequence}, if for every compact $K\subseteq T$ we
    have
    \begin{equation}\label{def-folner}
      \lim_{n\to\infty}\frac{\vol (\delta^{K}
        D_n)}{\vol (D_n)} = 0, 
    \end{equation}
    where $\delta^{K} D_n$ is the symmetric difference of $D_n$ and $KD_n$,
    \begin{displaymath}
    \delta^{K} D_n :=(KD_n)\setminus D_n\cup (KD_n)^c\setminus D_n^c.
    \end{displaymath}
  \item 
    The sequence $(D_n)_{n\in\mathbb N}$ is called a \emph{van Hove
      sequence}, if for every compact $K\subseteq T$ we have
    \begin{equation}\label{def-van-Hove}
      \lim_{n\to\infty}\frac{\vol (\partial^{K}
        D_n)}{\vol (D_n)} = 0, 
    \end{equation}
    where $\partial^{K} D_n:=(KD_n)\setminus \mathring{D_n}\cup
    (K\overline{D_n^c})\setminus D_n^c$.
  \item 
    The sequence $(D_n)_{n\in\mathbb N}$ is \emph{tempered}
    (or obeys \emph{Shulman's condition}),
    if there exists $C\ge1$ such that for all $n\in\mathbb N$ we have
    the estimate
    \begin{displaymath}
      \vol \left(\bigcup_{k=1}^{n-1} D_{k}^{-1}D_{n}\right)
      \le C \vol (D_n). 
    \end{displaymath}
  \end{nummer}
\end{definition}

\begin{remarks}
\item  For every $n\in\mathbb N$, the set $\partial^K D_n$ in \eqref{def-van-Hove} is compact. 
If $T$ is Abelian, then our definition of van Hove sequence is equivalent to that in \cite{Sch00}.
\item \label{hove-foelner}
  We have $\delta^{K} D\subseteq \partial^{K} D$, which follows
  from the inclusion $(AB)^c\subseteq AB^c$ for arbitrary
  $A,B\subseteq T$. Consequently, every van Hove sequence is
  a F{\o}lner sequence.
\item The existence of a F{\o}lner sequence in $T$ is equivalent to amenability of the group \cite[Thm. 4.16]{Pa00}. Every F{\o}lner sequence has a tempered subsequence \cite[Prop.~1.4]{Lind01}.
\item 
According to \cite{Str74}, every second countable, locally compact group has a
left-invariant proper metric that generates the topology. Suppose the sequence $(B_n)_{n\in\mathbb N}$ of
closed balls (with respect to this metric) about the neutral element $e\in T$ of radius
$n \in\mathbb N$ constitutes a (tempered) F{\o}lner sequence in $T$ satisfying $B_mB_n=B_{m+n}$ for all $m,n\in\mathbb N$.
Assume in addition $\mathrm{vol}(\partial B_n)/\mathrm{vol}(B_n)\to 0$ as $n\to\infty$, with
$\partial B_n$ the topological boundary of $B_n$. It can be shown that $(B_n)_{n\in\mathbb N}$ is also a
(tempered) van Hove sequence under these assumptions. 
\item If $T$ is Abelian, the existence of a tempered van Hove sequence in $T$ is
  guaranteed under our hypotheses. Indeed, \cite[p.~145]{Sch00} ensures the existence of 
  a van Hove sequence in $T$, which is also a F\o{}lner sequence by \itemref{hove-foelner}. 
  As every F\o{}lner sequence has a tempered subsequence, the 
  argument is complete.
\item Consider the semidirect product \cite{HeRo63}, denoted by $T=N\rtimes H$, of a 
unimodular group $N$ and a compact group $H$. Then $T$ is unimodular.
It can be shown that if $(D_n)_{n\in\mathbb N}$ is an $H$-invariant tempered van Hove sequence in $N$, 
then $(D_n\times H)_{n\in \mathbb N}$ is a tempered van Hove sequence in $T$. This can 
be used to provide examples of a 
non-Abelian non-compact unimodular group $T$ with a tempered van Hove sequence. 
(Take $H$ non-Abelian and $N$ Abelian but not compact.) A prominent example is the Euclidean group
$E(\mathsf{d})=\mathbb R^{\mathsf{d}}\rtimes O(\mathsf{d})$, with centered closed balls of radius $n\in\mathbb N$ as tempered van Hove sequence 
in $\mathbb R^{\mathsf{d}}$. The existence of F{\o}lner sequences in semidirect products is discussed in 
\cite{Ja08, Wil09}.
\end{remarks} 

The following lemma states that a F{\o}lner sequence $(D_n)_{n\in\NN}$ and its ``thickened''
version $(LD_n)_{n\in\NN}$, where $L\subseteq T$ is a compact set, have asymptotically 
the same volume. It also states that 
thickened versions of van Hove boundaries $\partial^{K} D_n$, with $K\subseteq T$ compact,
are of small volume, asymptotically as $n\to\infty$. These properties will be used
repeatedly below.

\begin{lemma}\label{volest}
  Let $L\subseteq T$ be a compact set. Then the following statements hold.
  \begin{indentnummer}
  \item \label{volestfoelner}
    If $(D_n)_{n\in\mathbb N}$ is a F{\o}lner sequence in $T$, we have the
    asymptotic estimate
    \begin{displaymath}
      \vol (LD_n)=\vol (D_n)+ \mathocal{o}\big(\vol (D_n)\big) \qquad
      (n\to\infty). 
    \end{displaymath} 
  \item \label{volesthove}
  	If $(D_n)_{n\in\mathbb N}$ is a van Hove sequence in $T$, we
    have for every compact $K\subseteq T$ the asymptotic estimate
    \begin{displaymath}
      \vol (L\partial^{K} D_n)=\mathocal{o}\big(\vol (D_n)\big) \qquad (n\to\infty).
    \end{displaymath}
  \end{indentnummer}
\end{lemma}

Next, we state the basic pointwise ergodic theorem that will be applied
several times in the sequel. Let $\cQ$ be a compact metrisable space 
(hence, with a countable base of the topology and complete with respect 
to every metric generating the topology), and assume that the group $T$ 
acts measurably from the left on $\cQ$, 
i.e., there exists a measurable map $\alpha_{\cQ}: T\times \cQ\to \cQ$,
$(x,q)\mapsto \alpha_{\cQ}(x,q) =:xq$.  Here, $T\times \cQ$ is endowed with
the product topology. 

A $T$-invariant probability measure on the Borel $\sigma$-algebra of $\cQ$ 
is called \mbox{($T$-)} \emph{ergodic}, if every $T$-invariant Borel set has either 
measure 0 or 1.  The existence of an ergodic probability measure on $\cQ$ 
follows from the compactness of $\cQ$ by standard arguments (compare \cite[\S6.2]{Wa82} for the discrete
case). In other words, $\cQ$ is ergodic w.r.t.\ the group $T$.  A 
dynamical system is called \emph{uniquely} \mbox{($T$-)} \emph{ergodic}, if it carries 
exactly one $T$-invariant probability measure, which is then ergodic, see below.

We rely on the general Birkhoff ergodic theorem of Lindenstrauss \cite[Thm.~1.2]{Lind01}. 
For related abstract ergodic theorems, see also \cite{Cha,Kre85, Pa00, Nev06}. The shorthand 
$\mu(f) := \int_{\cQ} \d\mu(q)\, f(q)$ in the next theorem denotes the $\mu$-integral of a function $f$ on $\cQ$. We remark that the assumptions on the group $T$ in the next theorem are more general than those required by Assumption~\ref{basic-comp}.

\begin{theorem}[Pointwise Ergodic Theorem]\label{abstract-erg}
  Let $\cQ$ be a compact metrisable space, on which a locally
  compact second-countable group $T$ acts measurably from the left. Assume that $T$
  admits a tempered F{\o}lner sequence $(D_n)_{n\in\mathbb N}$. Fix a $T$-invariant Borel probability
  measure $\mu$ on
  $\cQ$ and let $f\in L^{1}(\cQ,\mu)$ arbitrary be given. Then
  \begin{equation} \label{cesaro-av}
    I_n(q,f):= \frac{1}{\vol(D_{n})}  \int_{D_n}\!\d x\; f(xq)
  \end{equation}
  is finite for $\mu$-a.a.\ $q\in\cQ$ and all
  $n\in\NN$. Furthermore, there exists a $T$-invariant function $f^\star\in
  L^1(\cQ,\mu)$ such that $\mu(f^\star)=\mu(f)$ and 
  \begin{equation}\label{chaterg}
    \lim_{n\to\infty} I_n(q,f)= f^\star(q) \qquad\qquad\text{for $\mu$-a.a.\
      $q\in \cQ$.} 
  \end{equation}
  Moreover, the following statements are equivalent.
  \begin{indentnummer}
  \item The measure $\mu$ is ergodic.
  \item For every $f\in L^1(\cQ,\mu)$,
    Eq.~\eqref{chaterg} holds with $f^\star=\mu(f)$.
  \item There exists a $(\|\cdot\|_{\infty}$-$)$ dense subset
    $\mathcal{D}\subseteq C(\cQ)$ such that for every $f\in
    \mathcal{D}$, Eq.~\eqref{chaterg} holds with $f^\star=\mu(f)$.
  \end{indentnummer}  
\end{theorem}

\begin{remarks}
\item For $\mu$ ergodic, the limit \eqref{chaterg} is obviously independent of the tempered
  F{\o}lner sequence.
  \item \label{general-proof}
  As the proof shows, the statement of the theorem remains true for locally compact 
  Polish spaces (with $\mathcal{D}\subseteq C_c(\cQ)$ in (iii)). Moreover, our proof uses 
  local compactness only for the implication (iii)$\Rightarrow$(i). Non-compact dynamical systems have been studied in \cite{O52}.
  The reason why we assume even compactness of $\cQ$ in the hypotheses of the theorem is to guarantee the existence of an ergodic probability measure on $\cQ$.
  It is not obvious to us how
  to dispense with metrisability of $\cQ$.
\end{remarks}

In the case of a uniquely ergodic system, one may adapt arguments
from \cite{Fur81,Wa82} to exclude the exceptional set in \eqref{chaterg}, provided
that $f$ is continuous. Note that the next theorem is stated under more general assumptions 
on the group $T$ than those required by Assumption~\ref{basic-comp}.

\begin{theorem}[Unique ergodicity]\label{abstract-unierg}
 Let $\cQ$ be a compact metrisable space, on which a locally
  compact group $T$ acts measurably from the left. 
  Assume that $T$ admits a  F{\o}lner sequence $(D_n)_{n\in\mathbb N}$
  and define $I_{n}(\cdot,\cdot)$ as in \eqref{cesaro-av}.  
  Then the following statements are equivalent.
  \begin{indentnummer}
  \item For every $f\in C(\cQ)$ the sequence $\big(I_n(q,f)\big)_{n\in\NN}$
    converges uniformly in $q\in\cQ$ and there is a constant $I(f)\in\RR$
    such that 
    \begin{displaymath}
    \lim_{n\to\infty} I_n(q,f)=I(f)
    \end{displaymath}
    for all $q\in \cQ$. 
  \item There exists a dense subset $\mathcal{D}\subseteq C(\cQ)$
    and for every $f\in \mathcal{D}$ there exists a constant $I(f)\in\RR$
    such that pointwise for every $q\in \cQ$ we have
    \begin{displaymath}
	    \lim_{n\to\infty} I_n(q,f)=I(f).
    \end{displaymath}
  \item There exists exactly one $T$-invariant Borel probability
    measure $\mu$ on $\cQ$.
  \end{indentnummer}  
  In either case, the measure $\mu$ is ergodic and the above statements hold
  with $I(f)=\mu(f)$.
\end{theorem}

\begin{remark}
  In particular, the limits in the above theorem are again independent of the
  choice of the F{\o}lner sequence. In contrast to Theorem~\ref{abstract-erg},
  the F{\o}lner sequence here does not need to be tempered. Neither does
  one need second countability of the group $T$.
\end{remark}

The r\^{o}le of the compact space $\cQ$ in the above ergodic theorems will be
played by the closure of $T$-orbits of point sets.

\begin{definition} 
  Given a collection of point sets $\mathcal{P}\subseteq \mathcal P_r(M)$, we
  introduce its closed $T$-orbit
  \begin{equation}
  	\label{hull-def}
    X_\mathcal{P}:= \overline{\{xP:x\in T, P\in\mathcal{P}\}} \subseteq \cP_{r}(M),
  \end{equation}
  where $xP := \{ xp: p \in P\}$. Being closed, $X_\mathcal{P}$ is a compact subset of the compact space $\cP_{r}(M)$. The induced group action $\alpha_{X_{\cP}}: T \times X_{\cP} \rightarrow X_{\cP}, (x,P) 
  \mapsto xP$ is continuous.

\end{definition}

\begin{remarks}
\item
	The validity of the set inclusion in \eqref{hull-def} depends crucially on $\cP_{r}(M)$ being defined in terms of balls with respect to a $T$-invariant metric on $M$. The compatibility between the point space $M$ and the group $T$ is necessary in order to have a fruitful concept of orbits.   
\item
	The compact metrisable space $X_{\cP}$ is
	particularly useful if the closure is not too large in comparison to the
	(unclosed) $T$-orbit. This has been analysed mainly for $M=T=\mathbb R^{\mathsf{d}}$ with
	the canonical group action and the Euclidean metric.
	In that case, there are two simple examples where the closure does not add anything new to
	the (unclosed) $T$-orbit of $\cP$.  This is when $\cP$ consists of a single
	periodic point set, or when $\cP$ is a suitable collection of random tilings
	\cite{RHHB98,GrSh89}. 
	The definition of the vague topology suggests that elements in
	$X_{\mathcal P}$ added by the closure share local properties of point sets
	from $\mathcal P$. If $\cP=\{P\}$ consists of a single point set $P$, there is
	a geometric characterisation of $X_\cP$ as the so-called local 
	isomorphism class of the point set, iff  $P$ is repetitive.
	The latter property is in fact equivalent to minimality of $X_\cP$, if $\cP$ is of finite local complexity as in Definition~\ref{def:FLC}; cf.\ \cite{LaPl03} for
	$M=T=\RR^{\mathsf{d}}$ and \cite{Yok05} for the general case. Another criterion for a 
	``nice'' closure is unique ergodicity of $X_\cP$. 
	We give a geometric characterisation of unique ergodicity in Theorem~\ref{uniFLC}.
\end{remarks}

The triple $(X_\mathcal{P},T,\alpha_{X_{\cP}})$ constitutes a compact topological
dynamical system. Thus, we have

\begin{cor}\label{erg-cor}
  Let $\mathcal P \subseteq \mathcal P_r(M)$ be a collection of uniformly discrete
  point sets of radius $r$. Then the Ergodic Theorems~\ref{abstract-erg}
  and~\ref{abstract-unierg} hold for $\cQ =X_{\cP}$.
\end{cor}

\begin{remark}
	Ergodic theorems for systems of point sets in $\RR^{\mathsf{d}}$ or in a locally compact Abelian group have been 
  given and applied before, see e.g. \cite{Sch00, LeMo02}. In addition we mention
  \cite[Thm.~1]{LeSt05} for Banach space-valued functions in the case of minimal ergodic 
  systems of Delone sets of
  finite local complexity (see below for a definition) in $\RR^{\mathsf{d}}$.
\end{remark}

\subsection{Geometric characterisation of ergodicity} \label{geo-char}

In this section we relate ergodicity of a dynamical system of point sets to 
the spatial frequencies with which patterns occur therein. This will require 
to count the number of equivalent patterns within a given region of the point space, where
equivalence is defined by the group action.

First we introduce the relevant notation. Given a point set $P\in\mathcal P_r(M)$, we
call a finite subset $Q\subseteq P$ a \emph{pattern of $P$ (in $M$)}. Given a collection $\cP\subseteq 
\cP_r(M)$ of point sets, we say that $Q$ is a pattern of $\cP$, if there exists $P\in\cP$ such that
$Q$ is a pattern of $P$. We write $\cQ_{\cP}$ for the set of all patterns of $\cP$, 
see also Definition~\ref{pat-col}.
For a pattern $Q$ of $P$, every 
compact set $V\subseteq M$ such that $Q=P\cap \mathring V$ is called a \emph{support} of $Q$, and we say 
that $Q$ is a $V$-pattern of $P$. Two subsets $V,V'\subseteq M$ 
are called \mbox{($T$-)} \emph{equivalent}, if $xV=V'$ for some $x\in T$. 

For $ P \in \cP_{r}(M)$, $Q \subseteq P$ a pattern of $P$ and $D\subseteq T$, we analyse the
number of equivalent patterns of $Q$ in $P$. Fixing $m\in M$, one may consider the two different sets
\begin{align*}
M_{D}(Q) & :=\{\wtilde Q \subseteq P: \exists x\in D^{-1}: x Q=\wtilde Q\},\\
M'_{D}(Q) & :=\{\wtilde Q \subseteq P\cap D^{-1}m: \exists x\in T: x Q =\wtilde Q\}
\end{align*}
for this purpose. Note that the set $M'_{D}(Q)$ depends on the choice of  $m\in M$, in contrast to the set $M_{D}(Q)$. For this reason we will use $M_{D}(Q)$ for pattern counting in Definition~\ref{pattcount}. The set $M_{D}(Q)$ is a subset
of the equivalence class of $Q$. The next lemma describes how these two sets grow with the volume of $D$.

\begin{lemma} \label{count-compare2}
Assume that $T$ is even unimodular. 
Fix a F{\o}lner sequence $(D_n)_{n\in\mathbb N}$. 
Let $Q$ be a pattern of $\cP_r(M)$ and fix $P\in\mathcal P_r(M)$. 
Then we have the asymptotic estimates
\begin{displaymath} 
	\displaystyle\begin{array}{l}
		\mathrm{card}\big(M_{D_{n}}(Q)\big)=\mathocal{O}\left(\vol (D_n)\right), \\[1ex]
		\mathrm{card}\big(M'_{D_{n}}(Q)\big)=\mathocal{O}\left(\vol (D_n)\right),
 	\end{array}
 \qquad\quad (n\to\infty).
\end{displaymath}
If $(D_n)_{n\in\mathbb N}$ is even a van Hove sequence, and if $Q\subseteq Tm$, then
\begin{displaymath}
\mathrm{card}\big(M'_{D_{n}}(Q)\big) 
=\mathrm{card}\big(M_{D_{n}}(Q)\big) + \mathocal{o}\big(\vol (D_n)\big) \qquad (n\to\infty).
\end{displaymath}
The $\mathocal{O}$-terms and the $\mathocal{o}$-term may be chosen uniformly in $P\in\mathcal P_r(M)$.
\end{lemma}

\begin{remarks}
\item 
The condition $Q\subseteq Tm$ is satisfied for a transitive group action, since $Tm=M$ in that case.
\item
The number of equivalent copies of $Q$ in $P$ may also be analysed by counting corresponding 
group elements of the group $T$. One may consider the two different sets
\begin{align*}
T_{D}(Q) & :=\{x\in D^{-1}: xQ \subseteq P\},\\
T'_{D}(Q) & :=\{x\in T: xQ \subseteq P\cap D^{-1}m\}.
\end{align*}
The set $T'_{D}(Q)$ is commonly used for pattern counting, see 
\cite{Sch00, LeMo02}, but depends on the choice of $m\in M$. In order 
to relate this to the above approaches of pattern counting, consider 
the map $f:T_{D}^{(\prime)}(Q)\to M^{(\prime)}_{D}(Q)$,
given by $x\mapsto f(x):=xQ$. This map is onto. It is readily checked that $f$ is one-to-one,
if $Q\ne\varnothing$, the group $T$ is free on $Q$ and does not contain a nontrivial element of finite order. 
Hence, in that case, both approaches coincide.
\end{remarks}

Our central notion of pattern counting is

\begin{definition}
 	\label{pattcount}
 	Let $(D_n)_{n\in\mathbb N}$ be a F{\o}lner sequence in $T$ and let $P,Q \in \cP_{r}(M)$ be point sets with $|Q|<\infty$. (In particular, $Q$ may be a pattern of $P$). 
	If the limit
  \begin{displaymath}
    \nu(Q)  \equiv \nu^{P}\big(Q; (D_{n})_{n\in\NN}\big) := 
    \lim_{n\to\infty}\frac{\mathrm{card}\big(M_{D_n}(Q)\big)}{\vol(D_{n})}
  \end{displaymath} 
	exists, we call it the \emph{pattern frequency} of $Q$. In most cases we suppress its dependence on $P$ and the F{\o}lner sequence in our notation.
\end{definition}

\begin{lemma}\label{freq-gen}
	Assume that $T$ is even unimodular. Let $(D_n)_{n\in\mathbb N}$ be a 
	F{\o}lner sequence in $T$ and let $P,Q \in \cP_{r}(M)$ be point sets with $|Q|<\infty$. Then
	\begin{indentnummer}
	\item  
		the quotient which arises in the definition of the pattern frequency is bounded,
    \begin{displaymath}
    	\sup_{n\in\NN}\frac{\mathrm{card}\big(M_{D_n}(Q)\big)}{\vol(D_{n})} < \infty.
    \end{displaymath}
    In other words, since $\limsup$ and $\liminf$ of ${\mathrm{card}(M_{D_n}(Q))}/{\vol(D_{n})}$ 
    are always both finite, existence of the pattern frequency $\nu(Q)$ is only a matter of whether 
    they coincide.
	\item \label{freq-gen-shift} 
  	If $(D_n)_{n\in\mathbb N}$ is even a van Hove sequence and if the pattern frequency 
  	$\nu^{P}(Q; (D_{n})_{n\in\NN})$ exists, then $\nu^{P}(xQ; (D_{n})_{n\in\NN})$, 
		$\nu^{P}(Q; (xD_{n})_{n\in\NN})$ 
  	and $\nu^{xP}(Q; (D_{n})_{n\in\NN})$ exist and are all equal, i.e.,
    \begin{align*}
    	\nu^{P}\big(xQ; (D_{n})_{n\in\NN}\big) &= \nu^{P}\big(Q; (xD_{n})_{n\in\NN}\big) \\
	 			&= \nu^{xP}\big(Q; (D_{n})_{n\in\NN}\big) = \nu^{P}\big(Q; (D_{n})_{n\in\NN}\big)
    \end{align*}
   	for every $x\in T$.
	\end{indentnummer}
\end{lemma}

In order to relate ergodicity to pattern counting, we require a certain type
of rigidity for point sets.

\begin{definition}
Let $\cP\subseteq \cP_r(M)$ be a collection of point sets and $\cQ_{\cP}\subseteq \cP_r(M)$ the collection of its patterns.
\begin{nummer}
\item \label{def:FLC}
	$\cP$ is of \textit{finite local complexity} (FLC), if for every compact set 
	$V\subseteq M$ there is a finite collection $\cF_{\cP}(V)\subseteq \cQ_{\cP}$ of 
	(w.l.o.g.\ mutually non-equivalent) patterns, such that every pattern 
	of $\cP$, which admits a support equivalent to $V$, is equivalent to some 
	pattern in $\cF_{\cP}(V)$.
\item $\cP$ is \textit{locally rigid}, if for every $Q\in\cQ_{\cP}$ there exists $\varepsilon>0$ 
	such that for all $\wtilde Q\in\cQ_{\cP}$ and for all $x\in T$ the properties 
	$x\wtilde Q\subseteq (Q)_\varepsilon$ and $Q\subseteq (x\wtilde Q)_\varepsilon$ imply that 
	$Q$ and $\wtilde Q$ are equivalent.
\end{nummer}	
\end{definition}

The following lemma discusses and relates the above notions. For $\cP\subseteq \cP_r(M)$
and $V\subseteq M$, define $\cP\wedge V:=\{P\cap V:P\in\cP\}\subseteq \cP_r(M)$.

\begin{lemma}\label{FLClocal}
Let $\mathcal P\subseteq \cP_r(M)$  be a collection of point sets. Then
\begin{indentnummer}
\item \label{FLC-transfer}
	$\cP$ is FLC if and only if $X_\cP$ is FLC.
\item If $\cP$ is FLC, then $\cP$ is locally rigid.
\item If $\cP$ is locally rigid and if $\cQ_{\cP}\wedge V$ is closed in $\cP_r(M)$
for all compact $V\subseteq M$, then $\cP$ is FLC.
\end{indentnummer}

\end{lemma}

\begin{remarks}\label{remFLC}
\item If $\cP$ is finite, then FLC is equivalent to local rigidity. This holds since for finite $\cP$
the set $\cQ_{\cP}\wedge V$ is finite for all compact $V\subseteq M$,
due to uniform discreteness. In particular, it is closed in $\cP_r(M)$.
\item \label{QP=QXP}
	The proof of Lemma~\ref{FLC-transfer} shows that in the FLC case every pattern of $X_{\cP}$
is equivalent to some pattern of $\cP$.
\end{remarks}

Restricting to collections of point sets of finite local complexity, we can now state a 
geometric characterisation of ergodicity and of unique ergodicity.

\begin{theorem}[Ergodicity for FLC sets]\label{uniFLC}
Assume that $T$ is even unimodular and that $T$ has a tempered van Hove sequence 
$(D_n)_{n\in\mathbb N}$.
Let $\mathcal{P}\subseteq \mathcal P_r(M)$ be a collection of 
point sets of finite local complexity. 
Let $\mu$ be a $T$-invariant Borel probability measure on $X_\cP$. Then the following statements are 
equivalent.
\begin{indentnummer}
\item The measure $\mu$ is ergodic.
\item \label{erg-pattern-char}
	For every pattern $Q$ in  the set $\cQ_{\cP}$ of all patterns of $\cP$, 
	there is a subset $X\subseteq X_\cP$ 
  of full $\mu$-measure such that the pattern frequency $\nu(Q)=\nu^{P}(Q; (D_{n})_{n\in\NN})$
  exists for all $P\in X$ and is independent of $P\in X$.
\end{indentnummer} 
    If any of the above statements applies, then every pattern frequency $\nu(Q)$, $Q\in\cQ_{\cP}$, is
    independent of the choice of the tempered van Hove sequence. 
    
     The system $X_\cP$ is \emph{uniquely} ergodic iff \itemref{erg-pattern-char} holds
    for all patterns $Q \in\cQ_{\cP}$ with $X=X_\cP$, that is, for \emph{every} $P\in X_\cP$.
    In that case, the van Hove sequence needs not to be tempered, and every
    pattern frequency $\nu(Q)$, $Q\in\cQ_{\cP}$, is independent of the choice of the van Hove sequence.
    Furthermore, the convergence to the limit underlying the definition of each $\nu(Q)$ 
    is even uniform in $P\in X_\cP$.         
\end{theorem}

In the following proposition, we give a characterisation of unique ergodicity in terms of properties
of $\cP$ instead of $X_\cP$. This characterisation is often referred to as
\emph{uniform pattern frequencies}, compare \cite[Thm.~3.2]{Sch00}, 
\cite[Thm.~2.7]{LeMo02}, and \cite[Def.~6.1]{LaPl03}.

\begin{definition}
Fix a van Hove sequence $(D_n)_{n\in\mathbb N}$ and let $\cP\subseteq \cP_r(M)$ be given. 
We say that $\cP$ has \textit{uniform pattern frequencies}, 
if for every pattern $Q$ of $\cP$  the sequence 
$\big(\nu_n^{y,P}(Q)\big)_{n\in\mathbb N}$, defined by
\begin{equation}\label{freq-cond-uni}
\nu_n^{y,P}(Q):=\frac{\card\big(\{\wtilde Q\subseteq P:\exists x\in D_n y: 
x\wtilde Q=Q\}\big)}{\vol (D_n)},
\end{equation}
converges uniformly in $(y,P)\in T\times \cP$, and if its limit is independent of $(y,P)\in T\times\cP$.
\end{definition}

\begin{remarks}
\item
If $M=T=\mathbb R^{\mathsf{d}}$ with the canonical group action, and if $\cP=\{P\}$ is linearly repetitive, 
then $\cP$ has uniform pattern frequencies, see \cite{LaPl03, DL}. 
\item 
If $\cP$ has uniform pattern frequencies, then the limit of \eqref{freq-cond-uni} is also independent 
of the choice of the van Hove sequence according to Theorem~\ref{uniFLC} and
\end{remarks}

\begin{proposition}[Unique ergodicity for FLC sets]\label{freq-cond12}
Assume that $T$ is even unimodular and has a van Hove sequence $(D_n)_{n\in\mathbb N}$. 
Let $\mathcal{P}\subseteq \mathcal P_r(M)$ be a collection of 
point sets of finite local complexity. Then the following statements are equivalent.
\begin{indentnummer}
\item $X_\cP$ is uniquely ergodic.
\item $\cP$ has uniform pattern frequencies.
\end{indentnummer}
\end{proposition}

At the end of this section we investigate which values an ergodic measure on $X_{\cP}$ can assign to cylinder sets. Cylinder sets play a prominent r\^ole in the constructions of \cite{LeMo02}, compare also \cite{Le08}, and are defined as follows: it is well-known \cite[Lemma 4.5]{Kel75} that the compact metrisable 
topological space $\mathcal P_r(M)$ can be embedded into the 
compact product space $\prod_{\varphi\in  C_c(M)}f_\varphi(\mathcal P_r(M))$, 
with injection map $i$ given by $i(P)(f_\varphi)=f_\varphi(P)$.  This motivates to call 
$f^{-1}_\varphi(O)\subseteq \cP_r(M)$ an \emph{open cylinder}
if $O\subseteq\mathbb R$ open and $\varphi\in C_c(M)$, and finite intersections 
thereof are called \emph{cylinder sets}.

We give an example of open cylinders in $\cP_r(M)$ with a simple geometric interpretation.
Let $U\subseteq M$ be an open ball such that $\mbox{diam}(U)<r$.
Take $\varphi\in C_c(M)$ such that $\varphi^{-1}(\{0\})=M\setminus
U$. (A possible choice is $\varphi=d(\cdot, U^c)$, where $d(\cdot,\cdot)$ denotes the metric on $M$.) 
Consider the open cylinder
\begin{displaymath}
C_{U} := \left\{P\in\mathcal P_r(M): f_{\varphi}(P)\ne 0\right\}
=f_\varphi^{-1}(\mathbb R\setminus \{0\}).
\end{displaymath}
It consists of all those point sets of $\mathcal P_r(M)$ which have exactly 
one point in $U$. Note that $C_U$ is independent of the particular choice 
of $\varphi\in C_c(M)$ with $\mathrm{supp}(\varphi)=\overline{U}$. For open
cylinders $C_{U_1},\ldots, C_{U_k}$ as above, we denote the cylinder set of their 
intersection by 
\begin{equation}
\label{cyl-set}
 	C_\mathbold{U} := \bigcap_{i=1}^{k} C_{U_{i}}.
\end{equation}

In the case of finite local complexity, the pattern frequencies determine the values 
which an ergodic measure assigns to such cylinder sets. The following proposition extends 
\cite[Cor.~2.8, Lemma~4.3]{LeMo02}.

\begin{proposition}\label{cylinder1}
Assume that $T$ is even unimodular and has a van Hove sequence. 
Let $\mathcal P\subseteq\mathcal P_r(M)$ be a collection 
of point sets of finite local complexity.
Let $\mu$ be an ergodic Borel probability measure on $X_\cP$.
Furthermore, let $Q=\{q_1,\ldots, q_k\}$, $k\in\mathbb N$, be a nonempty pattern of $X_\cP$. Choose $\varepsilon\in ]0,r/2[$ such that all patterns of $X_\cP$ in $(Q)_\varepsilon$ of
cardinality $k$ are equivalent to $Q$.
For $i\in\{1,\ldots,k\}$, consider the pairwise disjoint sets $U_i:=B_\varepsilon(q_i)$ and define $U:=\bigcup_{i=1}^k U_i$. Then the corresponding cylinder set \eqref{cyl-set} has $\mu$-measure
\begin{displaymath}
\mu(C_\mathbold{U})=\nu(Q) \, \vol (D_\varepsilon),
\end{displaymath}
with $D_\varepsilon:=\{x\in T: xQ\subseteq U\}\subseteq T$ being open
and relatively compact.

If $T$ is even Abelian and acts transitively on $M$, then we have the equality
\begin{displaymath}
\mathrm{vol}(D_\varepsilon)=\mathrm{card}(\mathcal S_k(Q))\, \zeta_{\varepsilon},
\end{displaymath}
where $\zeta_{\varepsilon} := \mathrm{vol}(\{x\in T: x m\in B_\varepsilon(m)\})$ does not depend on 
$m\in M$ and where $\mathcal{S}_k(Q)$ is the group of ``$T$-realisable'' permutations of $Q$, i.e.,
\begin{displaymath}
\mathcal S_k(Q):=\{\pi\in \mathcal S_k: \exists x\in T \text{~~such that~~} 
xq_{\pi(i)}=q_i \text{~~for all~~} i\in\{1,\ldots, k\}\},
\end{displaymath}
with $\mathcal S_k$ denoting the permutation group on $\{1,\ldots,k\}$.
\end{proposition}

%%%%%%%%%%%%%%%%%%%%%%%%%%%%%%%%%%%%%%%%%%%%%%%%%%%%%%%%%%%%%%%%%%%%%%%%%%%%
%%%%%%%%%%%%%%%%%%%%%%%%%%%%%%%%%%%%%%%%%%%%%%%%%%%%%%%%%%%%%%%%%%%%%%%%%%%%
%
\section{Dynamical systems for randomly coloured point sets}
\label{sec:ran-col}
%
%%%%%%%%%%%%%%%%%%%%%%%%%%%%%%%%%%%%%%%%%%%%%%%%%%%%%%%%%%%%%%%%%%%%%%%%%%%%
%%%%%%%%%%%%%%%%%%%%%%%%%%%%%%%%%%%%%%%%%%%%%%%%%%%%%%%%%%%%%%%%%%%%%%%%%%%%

In this section, we supply the point sets of the previous section with a
random colouring. The results obtained here will be applied to randomly
coloured graphs in the next section. All proofs are deferred to Section~\ref{sec:proofs-ran-col}.

In addition to a point space $M$ and a metrisable group $T$ satisfying Assumption~\ref{basic-comp} in the previous section, we consider a non-empty, locally compact, second-countable topological space $\AA$, 
which we call \emph{colour space}. 

\begin{lemma}
	\label{mhat-ass}
 	The product space $\what{M} := M \times\AA$, equipped with the product topology, constitutes 
	a point space in the sense of Section~\ref{psdyn}. The continuous and proper action $\alpha$ 
	of $T$ on $M$ induces a continuous and proper action $\what{\alpha}: T \times \what{M} 
	\rightarrow \what{M}$ of $T$ on $\what{M}$ by setting $\what{\alpha}\big(x, (m,a)\big) := (xm,a)$. 
	Thus, $\what{M}$ and $T$ satisfy Assumption~\ref{basic-comp}. We fix a $T$-invariant proper 
	metric $\what{d}$ on $\what{M}$ that is compatible with the topology on $\what{M}$.
\end{lemma}

\begin{remarks}
\item When acting on $\what{M}$, the group $T$ simply transports the colour $a$ of $m$ along with $m$. 
\item The maximum metric of the $T$-invariant proper metric $d$ on $M$ and some metric generating the topology on $\AA$ is an admissible choice for the metric $\what{d}$ in Lemma~\ref{mhat-ass}, because it is $T$-invariant and because every metric on $\what{M}$ which generates the product topology is equivalent to the proper metric $\what{d}$ of Lemma~\ref{mhat-ass}, and hence proper itself.
\end{remarks}

Lemma~\ref{mhat-ass} and the arguments in the previous section imply that the space 
$\cP_{r}(\what{M})$ of uniformly discrete point sets with radius $r$ in $\what{M}$ 
is a compact metrisable space with respect to the vague topology. The vague topology is defined as in Definition~\ref{fphi}, but with $M$ replaced by $\what{M}$.
The continuous action $\what{\alpha}$ induces a continuous
group action on $\cP_{r}(\what{M})$ by setting 
\begin{equation}
  \label{hat-trans}
  x\what{P} := \big\{(xm,a) \in \what{M}: (m,a) \in \what{P}\big\} \in \cP_{r}(\what{M})
\end{equation}
for $\what{P} \in \cP_{r}(\what{M})$. Again, the group action does not lead out of $\cP_{r}(\what{M})$ because of $T$-invariance of the metric $\what{d}$. To summarise, all results established
for $\cP_{r}(M)$ in Section~\ref{psdyn} remain true for $\cP_{r}(\what{M})$.

Rather than working with general subsets of $\what{M}$, we are interested in
those subsets for which each point of $M$ comes with exactly one colour.

\begin{definition}
  \begin{nummer}
  \item \label{col-ps-def}
    For a given point set $P\subseteq M$ we set $\Omega_{P} :=
    \bigtimes\nolimits_{p\in P}\AA$ and call
    $P^{(\omega)}=\big\{\big(p,\omega(p)\big):p\in P\big\} \subseteq \what{M}$ a
    \emph{coloured point set} with colour realisation $\omega\in\Omega_{P}$.
  \item Given a collection $\cP \subseteq \cP_{r}(M)$ of point sets, we
    introduce the collection of all associated coloured point sets $\cC_{\cP}
    := \{P^{(\omega)} : P\in\cP, \omega\in\Omega_{P}\}$. In particular, we
    write $\cC_{r}(M) := \cC_{\cP_{r}(M)}$ for the space of coloured,
    uniformly discrete points sets of radius $r$ and $\hXcP :=
    \cC_{X_{\cP}}$ for the space of coloured closed $T$-orbits.
  \end{nummer}
\end{definition}

\begin{remarks}
\item Let $\pi: \what{M} \rightarrow M$, $(m,a) \mapsto m$, be the canonical
  projection onto the space $M$. Then, $\what{P} \subseteq \what{M}$ is a
  coloured point set if and only if the restriction $\pi\big|_{\what{P}}$ is
  injective.
\item If $P\in\cP_{r}(M)$, then $P^{(\omega)} \in
  \cP_{r}(\what{M})$ for every $\omega\in\Omega_{P}$. Thus, we have
  $\cC_{\cP}\subseteq \cP_{r}(\what{M})$ in the above definition, and
  $\cC_{\cP}$ inherits the vague topology from $\cP_{r}(\what{M})$.
\item We equip $\Omega_P$ with the product topology (which is metrisable,
  since $\AA$ is a metric space and the product is countable). The product
  topology on $\Omega_{P}$ and the vague topology on $\cC_{P}$ coincide, when
  the two spaces are canonically identified by $\omega \leftrightarrow
  P^{(\omega)}$. This is seen by noting that for both topologies convergence
  means pointwise convergence.
\end{remarks}

Compactness of spaces of coloured point sets is established in

\begin{proposition}\label{cCr}
  Let $\mathcal P\subseteq\mathcal P_r(M)$ be given and assume that
  $\mathcal P$ is closed in $\mathcal P_r(M)$. Then the metrisable space $\mathcal
  C_\mathcal{P}$ is closed in $\mathcal P_r(\what{M})$ and hence compact. In
  particular, $\cC_r(M)$ and $\hXcP$ are compact.
\end{proposition}

It follows from \eqref{hat-trans} that the action of $x\in T$ on a coloured
point set $P^{(\omega)} \in \cC_{r}(M)$ can be described as
\begin{equation}
  \label{col-ps-trans}
  xP^{(\omega)} = (xP)^{(\tau_{x}\omega)}.
\end{equation}
Here, we introduced the measurable \emph{shift} $\tau_{x}: \Omega_{P} \rightarrow \Omega_{xP}$,
$\omega \mapsto \tau_{x}\omega$, between probability spaces, given by $(\tau_{x}\omega)(xp) := \omega(p)$ for
all $p\in P$. 

We are particularly interested in $T$-invariant,
compact spaces of coloured point sets. Therefore the following is useful.

\begin{lemma}
  \label{col-T-inv}
  For $\mathcal P \subseteq \mathcal P_r(M)$ we have
  \begin{displaymath}
    \hXcP = \overline{\big\{ xP^{(\omega)}: x\in T, \, P^{(\omega)} \in
      \cC_{\cP}\big\}} 
  \end{displaymath}
  where the closure is taken with respect to the vague topology. The group
  action $\alpha_{\hXcP} : T\times\hXcP \rightarrow\hXcP$, $(x,\what{P})
  \mapsto x\what{P}$ is continuous.  
\end{lemma}

The preceding proposition and lemma imply

\begin{cor}
  Let $\mathcal P \subseteq \mathcal P_r(M)$ be a collection of
  point sets. Then $(\hXcP, T, \alpha_{\hXcP})$ is a compact
  topological dynamical system and the Ergodic Theorems~\ref{abstract-erg}
  and~\ref{abstract-unierg} hold for $\cQ =\hXcP$. \qed
\end{cor}

Since we want to describe randomly coloured point sets, we will now
introduce suitable probability measures on the Borel spaces $(\Omega_{P},
\cA_{P})$ for different $P$. Here, $\cA_{P} := \bigotimes_{p\in P}\,
\mathcal \cA$ is the product over all points in $P$ of the Borel
$\sigma$-algebra $\cA$ on $\AA$. It coincides with the Borel $\sigma$-algebra on
$\Omega_{P}$ \cite[Lemma~1.2]{Kal01}.  For $V\subseteq M$ and
$P\in \cP_{r}(M)$ we define the local $\sigma$-algebra $\cA_{P}^{(V)}$ as
the smallest $\sigma$-algebra on $\Omega_{P}$ such that the canonical
projection $(\Omega_{P}, \cA_{P}^{(V)}) \rightarrow (\Omega_{P\cap V},
\cA_{P \cap V})$ is measurable. It is the $\sigma$-algebra of events
concerning only colours attached to points in $P\cap V$.

For $\cP\subseteq \cP_{r}(M)$ consider a family of Borel probability measures
$\PP_{P}$ on $(\Omega_{P},\cA_{P})$, which is indexed by $P\in X_{\cP}$.

\begin{ass}
  \label{p-ass}
  This is a list of properties which the family $(\PP_{P})_{P\in X_{\cP}}$ may
  satisfy or not.
  \begin{indentnummer}
  \item \label{p-ass-cov}
    \emph{$T$-covariance.} \quad $\PP_{xP} = \PP_{P} \comp \tau_{x}^{-1}$ for
    all $x\in T$ and all $P\in X_{\cP}$.
  \item \label{p-ass-iad}
  	\emph{Independence at a distance.}
		There exists a length $\varrho>0$ such that for every $P \in X_{\cP}$ 
		and every $V_{1}, V_{2} \subset M$ with $d(V_{1}, V_{2}) > \varrho$ 
		the local $\sigma$-algebras
    $\cA_{P}^{(V_{1})}, \cA_{P}^{(V_{2})}$ are $\PP_{P}$-independent. 
	\item \label{p-ass-wcomp} 
    \emph{M-compatibility.} \quad For every $f\in C(\hXcP)$
    the colour average $E_{f} :X_{\cP}\to\RR$, defined by
    \begin{displaymath}
      \label{col-av}
      E_f(P):=\int_{\Omega_P}\d\PP_P(\omega)\, f(P^{(\omega)}), \qquad\qquad P
      \in X_{\cP},
    \end{displaymath}
    is a measurable function.
  \item \label{p-ass-scomp} 
    \emph{$C$-compatibility.} \quad For every $f\in C(\hXcP)$
    we have $E_{f} \in C(X_{\cP})$.
  \end{indentnummer}
\end{ass}

Independence at a distance can be interpreted as interactions of finite range between points in the point set. Such type of interaction has relevance in statistical mechanics, see for example \cite{Sla80}.

In the next lemma we give two examples for a family $(\PP_{P})_{P\in X_{\cP}}$ of measures which
satisfy all of the above assumptions. The second example involves a random field  $\xi: \Sigma \times M \rightarrow \AA$, $(\sigma,m) \mapsto \xi^{(\sigma)}(m)$
over $M$ with values in $\AA$, where $(\Sigma,\mathcal{A}',\PP')$ is some given underlying probability space \cite{A}. We will also need the \emph{strong mixing coefficient}  \cite{Dou94} of $\xi$, defined by 
\begin{displaymath}
  \RR_{\ge 0} \ni L \mapsto \kappa(L) := \sup\bigl\{ \kappa(V_{1},V_{2}):
  V_{1},V_{2}\subseteq M, \;
   d(V_{1},V_{2})> L \bigr\},
\end{displaymath}
  where
\begin{align*}
	\kappa(V_{1},V_{2})&:= \sup\bigl\{ |\PP'(A'_{1}  \cap A'_{2}) - \PP'(A'_{1})\PP'(A'_{2})| : A'_{j}\in\mathcal{A}'(V_{j}) \text{~~for~} j\in\{1,2\}\bigr\}  \\
	& \phantom{:}\le 1/4
\end{align*}
measures the correlation of the local sub-$\sigma$-algebras
$\mathcal{A}'(V_{j})$ of events generated by the family of
random variables $\{\xi^{(\cdot)}(m) :\, m\in V_{j} \}$.

\begin{lemma} \label{IAD}
  \begin{nummer}
  \item \label{iid-ex}
    Let $\PP$ be a Borel probability measure on $(\AA, \cA)$, and define
    $\PP_{P} := \bigotimes_{p\in P}\PP$ for every $P \in X_{\cP}$. Then, the
    family of measures $(\PP_{P})_{P\in X_{\cP}}$ satisfies the
    Assumptions~\ref{p-ass-cov} --~\itemref{p-ass-scomp}.
  \item 
  	Let $(\Sigma,\mathcal{A}',\PP')$ be a probability space and let $\xi:
    \Sigma \times M \rightarrow \AA$, $(\sigma,m) \mapsto \xi^{(\sigma)}(m)$
    be an $\AA$-valued random field over $M$ which is jointly measurable,
    $T$-stationary, has a compactly supported strong
    mixing coefficient, and has continuous realisations $\xi^{(\sigma)}: M 
		\rightarrow
    \AA$ for $\PP'$-a.a.\ $\sigma\in\Sigma$. For a given point set $P \in X_{\cP}$ we
    define the map $\Xi_{P} : \Sigma \rightarrow \Omega_{P}$, $\sigma \mapsto
    \Xi_{P}(\sigma) := \xi^{(\sigma)}\big|_{P}$.  Then, $\PP_{P} :=
    \PP'\comp\Xi_{P}^{-1}$ is a Borel probability measure on $(\Omega_{P},
    \cA_{P})$, and the family of measures $(\PP_{P})_{P\in X_{\cP}}$ satisfies
    the Assumptions~\ref{p-ass-cov} --~\itemref{p-ass-scomp}.
  \end{nummer}
\end{lemma}

The main goal of this section is to characterise an ergodic Borel probability
measure $\hatmu$ on $\hXcP$ in terms of an ergodic Borel probability measure
$\mu$ on uncoloured point sets $X_{\cP}$ and the colour probability measures
$(\PP_{P})_{P\in X_\cP}$. This will be achieved in Theorem~\ref{theo:perg} below. A
crucial ingredient of the proof is the following statement, which is inspired by
and generalises \cite[Lemma~3.1]{Hof98}.

\begin{theorem}[Strong law of large numbers]\label{lem:hoflemma}
  Assume that $T$ is even unimodular and admits a F{\o}lner 
  sequence $(D_n)_{n \in\mathbb N}$.
  Fix $f\in C(\hXcP)$, $P\in X_{\cP}$ and suppose that $\PP_{P}$ has
  the property ``independence at a distance'' as in
  Assumption~\ref{p-ass-iad}. For $n\in\mathbb N$ define the random variable
  $Y_n:\Omega_P\to\mathbb R$ by
  \begin{equation}\label{Yn}
    Y_n(\omega) :=  \frac{1}{\vol(D_{n})}\int_{D_n}{\rm
      d}x\,f(xP^{(\omega)}), \qquad \omega\in\Omega_P.  
  \end{equation}
  Then we have for $\mathbb P_P$-almost all $\omega\in\Omega_{P}$ the relation
  \begin{equation}
    \label{hofstatement}
    \lim_{n\to\infty}\left(Y_n(\omega) - \int_{\Omega_P}{\rm d}\mathbb
      P_P(\eta) \, Y_n(\eta)\right) =0. 
  \end{equation}
\end{theorem}

We are now ready for the main
result of this section.  

\begin{theorem}
  \label{theo:perg}
  Assume that $T$ is even unimodular. 
  Let $\cP \subseteq \mathcal P_r(M)$ be a collection of point sets.  Fix an
  ergodic Borel probability measure $\mu$ on $X_{\cP}$ and a family of Borel
  probability measures $(\PP_{P})_{P\in X_{\cP}}$ satisfying
  Assumptions~\ref{p-ass-cov} -- \itemref{p-ass-wcomp}. Then there exists a
  unique ergodic probability measure $\hatmu$ on $\hXcP$ such that the
  following statements hold.
  \begin{indentnummer}
  \item
    \label{fubini}
    For every $f\in L^1(\hXcP,\hatmu)$ we have
    \begin{equation} \label{fubini-eq}
      \int_{\hXcP} \!\d\hatmu(P^{(\omega)}) \; f(P^{(\omega)}) =
      \int_{X_{\cP}} \!\d\mu(P) \int_{\Omega_{P}}\!\d\PP_{P}(\omega) \;
      f(P^{(\omega)}).
    \end{equation}
  \item 
    \label{item-erg}
    For every $f\in L^1(\hXcP,\hatmu)$ and every tempered 
    F\o{}lner sequence $(D_{n})_{n\in\NN}$ in $T$ the limit
    \begin{equation}
      \label{erg-limit}
      \lim_{n\to\infty} \frac{1}{\vol(D_{n})} \int_{D_n}\!\d x\;
      f(xP^{(\omega)}) = \int_{\hXcP} \!\d\hatmu(Q^{(\sigma)})\;
      f(Q^{(\sigma)})   
    \end{equation}
    exists for $\hatmu$-a.a.\ $P^{(\omega)} \in \hXcP$. In fact, the limit
    exists for $\mu$-a.a. $P \in X_{\cP}$ and for $\mathbb P_P$-a.a.\
    $\omega\in\Omega_P$.  
  \end{indentnummer}
  If $X_{\cP}$ is even uniquely ergodic, if $(\PP_{P})_{P\in X_{\cP}}$
  satisfies also Assumption~\ref{p-ass-scomp} and if $f$ is
  continuous, then the limit \eqref{erg-limit} exists \emph{for all} $P\in
  X_{\cP}$ and for $\mathbb P_P$-a.a.\ $\omega\in\Omega_P$. In this 
  case the F\o{}lner sequence does not need to be tempered.
\end{theorem}

\begin{remarks} 
\item 
  The asserted uniqueness of the ergodic measure $\hatmu$ in the theorem
  does not mean that the dynamical system $\hXcP$ is uniquely ergodic. It only
  means that $\hatmu$ is uniquely determined by the given ergodic measure
  $\mu$ on $X_{\cP}$ and the measures $(\PP_{P})_{P\in X_\cP}$ on $\Omega_{P}$.
\item 
  The corresponding Theorem~3.1 of \cite{Hof98} is a statement about
  Bernoulli site percolation on the Penrose tiling. Our result is an
  extension, which covers both the aperiodic and the periodic situation, under
  weaker assumptions on the underlying point set, and for more
  general types of randomness.
\item In contrast to a corresponding result \cite[Lemma 10]{Le08}, our
  Theorem~\ref{theo:perg} does not require a group structure of the point
  space $M$. Theorem~\ref{theo:perg} also makes a stronger conclusion in that
  exceptional instances are characterised beyond being $\hatmu$-null
  sets. This is particularly useful in the uniquely ergodic case.
\end{remarks}

At the end of this section we discuss which values the measure $\what\mu$ assigns 
to cylinder sets of coloured point sets. In view of the decomposition \eqref{fubini-eq} 
we are interested in the relation to the $\mu$-measure of the corresponding uncoloured cylinder set, see 
Section~\ref{geo-char}.
Consider open sets $U_1,\ldots, U_k\subseteq M$ with
$\mathrm{diam}(U_i)<r$ for $i\in\{1,\ldots,k\}$. Choose $\varphi_i\in C_c(M)$ such 
that $\varphi_i^{-1}(\{0\})=M\setminus U_i$ for $i\in\{1,\ldots,k\}$. Similarly, consider open,
relatively compact sets $A_1,\ldots, A_k\subseteq \AA$ and choose 
$\psi_1,\ldots,\psi_k\in C_c(\AA)$ such that $\psi_i^{-1}(\{0\})=\AA\setminus A_i$.
With $f_{\varphi,\psi}$ as in \eqref{fphipsi}, we define the \emph{coloured cylinder set}
\begin{equation}\label{colcylinder}
C_\mathbold{U}^\mathbold{A} := \Big\{P^{(\omega)}\in\cC_r(M):
  f_{\varphi_1,\psi_1}(P^{(\omega)}) \cdot\ldots\cdot
  f_{\varphi_k,\psi_k}(P^{(\omega)})\ne 0 \Big\},
\end{equation}
where $\mathbold{U}:=U_1\times\ldots\times U_k$ and $\mathbold{A} := A_1\times\ldots\times A_k$.
The set $C_\mathbold{U}^\mathbold{A}$ is independent of the particular choice of 
the functions $\varphi_i$ and $\psi_i$ with $\mathrm{supp}(\varphi_i)=\overline{U_i}$ and 
$\mathrm{supp}(\psi_i)=\overline{A_i}$, and it consists of all coloured point sets 
such that $U_i$ contains 
exactly one point of the underlying point set, with corresponding colour value in $A_i$, for $i\in\{1,\ldots,k\}$. In the case of independent and identically distributed (i.i.d.) colours, we have a nice product formula  for the measure of such a cylinder set, 
which is stated in

\begin{proposition} \label{col-cyl-prop}
Assume that $T$ is even unimodular and admits a F{\o}lner sequence.
Fix $f\in C(\hXcP)$, $P\in X_{\cP}$ and an ergodic Borel probability measure $\mu$ on $X_\cP$. 
Let $\PP$ be a Borel probability measure on $(\AA, \cA)$ and for every $P \in X_{\cP}$ consider the 
product measure $\PP_{P} := \bigotimes_{p\in P}\PP$ on $\Omega_P$, see Lemma~\ref{iid-ex}. Assume in addition that the sets $U_{1}, \ldots, U_{k}$ of the coloured cylinder set $C_\mathbold{U}^\mathbold{A}$ in \eqref{colcylinder} are pairwise disjoint. Then
\begin{displaymath}
\hatmu(C_\mathbold{U}^\mathbold{A})=\mu(C_\mathbold{U})\,\,
\PP(A_1)\cdot\ldots\cdot\PP(A_k),
\end{displaymath}
where $C_\mathbold{U}\subseteq X_\cP$ is the corresponding uncoloured cylinder set \eqref{cyl-set}.
\end{proposition}

%%%%%%%%%%%%%%%%%%%%%%%%%%%%%%%%%%%%%%%%%%%%%%%%%%%%%%%%%%%%%%%%%%%%%%%%%%
%%%%%%%%%%%%%%%%%%%%%%%%%%%%%%%%%%%%%%%%%%%%%%%%%%%%%%%%%%%%%%%%%%%%%%%%%%
%
\section{Application to graphs}
\label{secdyn}
%
%%%%%%%%%%%%%%%%%%%%%%%%%%%%%%%%%%%%%%%%%%%%%%%%%%%%%%%%%%%%%%%%%%%%%%%%%%
%%%%%%%%%%%%%%%%%%%%%%%%%%%%%%%%%%%%%%%%%%%%%%%%%%%%%%%%%%%%%%%%%%%%%%%%%%

One of our reasons for dealing with point spaces without a group structure in the previous 
sections is that this allows for a description of simple graphs \cite{Die05}. Most statements 
in this section do not require extra proofs, because they follow from applying the general 
results of Sections~\ref{psdyn} and~\ref{sec:ran-col}. More generally, one could treat simple directed graphs or even hypergraphs \cite{Die05} by the same methods. 

%%%%%%%%%%%%%%%%%%%%%%%%%%%%%%%%%%%%%%%%%%%%%%%%%%%%%%%%%%%%%%%%%%%%%%%%%%
\subsection{Graphs as point sets}

Let $\VV$ be a point space, i.e.\ a non-empty, locally compact and second-countable Hausdorff
space and $T$ a metrisable group. As in Section~\ref{top-point-sets}, we require

\begin{ass}
	\label{basic-comp-graph}
	The group $T$ is non-compact and its left action $\alpha_{\VV}: T \times \VV 
	\rightarrow \VV$, $(x,v) \mapsto xv$, on $\VV$ is continuous and proper. Moreover, we fix a 
 	$T$-invariant proper metric $d_{\VV}$ on $\VV$ that generates the topology on $\VV$.
\end{ass}

\noindent
Hence, $T$ is also locally compact and second countable, compare Remarks~\ref{prop-implies1} and \itemref{prop-implies2}.

Next we consider the space $M := ( \VV \times \VV)/\!\sim$ 
with the quotient topology, arising from $\VV\times\VV$ with the product topology of $\VV$. Here, the equivalence relation $\sim$ identifies $(v,w)\in
\VV\times \VV$ with $(w,v)\in \VV\times \VV$, and we write $m=m_{v,w} =
m_{w,v}\in M$ for the corresponding equivalence class. Clearly, the definition
\begin{align*} %\label{metric-M}
  d(m_{v_{1},w_{1}}, m_{v_{2},w_{2}}) := \min \Big\{
  &\max\big\{d_{\VV}(v_{1},v_{2}), d_{\VV}(w_{1},w_{2}) \big\}, \nonumber\\
  &\max\big\{d_{\VV}(v_{1},w_{2}), d_{\VV}(v_{2},w_{1}) \big\} \Big\}  
\end{align*}
for all $v_{1},v_{2},w_{1},w_{2} \in\VV$ provides a $T$-invariant proper metric $d$ on $M$ that is compatible with the topology on $M$.

\begin{lemma}
 The space $M$ is a point space and, together with the induced action $\alpha:
T \times M \rightarrow M$, $\alpha (x, m_{v,w}) := x m_{v,w} := m_{xv,xw}$ of $T$ on $M$, satisfies Assumption~\ref{basic-comp}.
\end{lemma}

As to the proof the lemma, we note that $M$ is clearly non-empty, locally compact and second-countable, and thus a point space, equipped with the $T$-invariant proper metric above. Also, continuity of the induced action $\alpha$ is evident. The proof is then completed by the first part of 

\begin{lemma} \label{alpha-trans}
 	\begin{nummer}
 	\item \label{proper-trans}
  	The induced action $\alpha$ on $M$ is proper.
	\item \label{free-trans}
		If $T$ does not contain an element of order two and $\alpha_\VV$ acts freely, then
		so does $\alpha$.
	\end{nummer}
\end{lemma}

\begin{remarks}
	\item A proof of the lemma can be found in Section~\ref{sec:proofs-dyn}.
	\item Transitivity of $\alpha_{\VV}$ does not imply transitivity of $\alpha$.
\end{remarks}

\begin{definition}\label{def:graph}
  A point set $G\subseteq M$ is called a (simple) \emph{graph} (in $\VV$), if $m_{v,w}\in
  G$ for $v,w\in\VV$ implies $m_{v,v}\in G$ and $m_{w,w}\in G$.  A graph $G
  \subseteq M$ has the vertex set $\cV_{G}:=\{v\in \VV: m_{v,v} \in G\}$, which
  is a point set in $\VV$, and its edge set is given by $\cE_{G} :=\big\{
  \{v,w\}: v,w \in\VV, v\neq w, m_{v,w} \in G \big\}$.
\end{definition}

\begin{remark}
  Every graph $G \subseteq M$ is a simple graph \cite{Die05}, that is, without
  self-loops or multiple edges between the same pair of vertices.  It is easy
  to see that $G$ is a uniformly discrete subset of $M$ with radius $r$, if
  and only if $\cV_{G}$ is uniformly discrete in $\VV$ with the same
  radius. Also, relative denseness of $G$ with radius $R$ implies relative
  denseness of $\cV_{G}$ with radius $R$.  The converse statement does not
  hold. This is seen from a graph with relatively dense vertex set, but
  without edges. Relative denseness of the point set $G$ implies the existence
  of vertices with infinitely many incident edges.
\end{remark}

%%%%%%%%%%%%%%%%%%%%%%%%%%%%%%%%%%%%%%%%%%%%%%%%%%%%%%%%%%%%%%%%%%%%%%%%%%

\subsection{Ergodicity for dynamical systems of graphs}

In this section we will apply the ergodic results from Section~\ref{psdyn} to graphs.
The space
$\cP_{r}(M)$ of uniformly discrete point sets in $M = (\VV \times \VV) /\!\sim$ with radius $r$ is a
compact metrisable space with respect to the vague topology of Definition~\ref{fphi}. 
The action $\alpha$ on $M$ induces in turn a continuous group action
on $\cP_{r}(M)$ by pointwise shifts as in Section~\ref{psdyn}.
Consequently, all results established for $\cP_{r}(M)$ in
Section~\ref{psdyn} are available in the present context.

\begin{definition}
  For $r\in]0,\infty[$ we introduce the space
  \begin{displaymath}
    \cG_{r}(M) := \big\{ G \in\cP_{r}(M): G \text{~is a graph}\big\}
  \end{displaymath}
  of graphs in $\VV$ with uniformly discrete vertex sets of radius $r$. It
  inherits the vague topology from $\cP_{r}(M)$.
\end{definition}

We omit the obvious proof of

\begin{proposition}
  $\cG_{r}(M)$ is closed, hence compact in
  $\cP_{r}(M)$. Moreover, $\cG_{r}(M)$ is $T$-invariant. \qed
\end{proposition}

Again, we are interested in closed (hence compact) and $T$-invariant
subsets of $\cG_{r}(M)$. An example is given by
\begin{displaymath}
  X_{\cG} := \overline{ \{xG: x \in T, G \in\cG \} } \subseteq \cG_{r}(M)
\end{displaymath}
for some given subset $\cG \subseteq \cG_{r}(M)$. We denote the continuous
group action of $T$ on $X_{\cG}$ by $\alpha_{X_{\cG}}$ and are now ready
to apply the general ergodic results of Section~\ref{psdyn}.

\begin{cor}
  Let $\cG \subseteq \cG_r(M)$ be given. Then $(X_{\cG}, T, \alpha_{X_{\cG}})$ is a compact
  topological dynamical system, and the Ergodic Theorems~\ref{abstract-erg}
  and~\ref{abstract-unierg} hold for $\cQ =X_{\cG}$. \qed
\end{cor}

$G'$ is called a \emph{subgraph} of $G$, if $G'\subseteq G$ and if $G'$ is a graph.  
A subgraph $G'$ of $G$ with a finite vertex set is called a \emph{patch} of $G$. 
Every patch of $G$ is a pattern of $G$. A pattern $Q$ of $G$ is a patch of $G$
if and only if $Q$ is a subgraph of $G$. For a collection $\cG\subseteq \cG_r(M)$ of graphs, 
we call $G'$ a patch of $\cG$, if there exists $G\in\cG$ such that $G'$ is a patch of $G$.
Every compact set $V\subseteq\VV$ such that $\cV_{G'}=\cV_{G}\cap \mathring V$ 
is called a \emph{support} of the patch.

It is easy to see that a collection $\cG\subseteq \cG_r(M)$ of graphs has
finite local complexity in the sense of Definition~\ref{def:FLC} if and only 
if for every compact set $V\subseteq \VV$ there is a finite collection $\mathcal F_\cG(V)$ 
of patches, such that every patch of $\cG$, which admits a support equivalent to 
$V$, is equivalent to some patch in $\mathcal F_\cG(V)$.

In analogy to the case of point sets, we have the following characterisation of ergodicity
and of unique ergodicity.

\begin{theorem}[Ergodicity for FLC graphs]\label{uniFLCgraph}
Assume that $T$ is even unimodular and that $(D_n)_{n\in\mathbb N}$ is a 
tempered van Hove sequence in $T$. 
Let 
$\mathcal G\subseteq\mathcal G_r(M)$ be a collection of graphs of
finite local complexity.   Let 
$\mu$ be a $T$-invariant Borel probability measure on $X_\cG$. Then
the following statements are equivalent.
\begin{indentnummer}
\item The measure $\mu$ is ergodic.
\item \label{graph-FLC-pattern}
	For every patch $H$ of $\mathcal G$, there is a
	set $X\subseteq X_\cG$ of full $\mu$-measure, such that the limit
    \begin{equation}\label{freq-cond2}
      \nu(H) := \lim_{n\to\infty}\frac{\mathrm{card}\big(\{\wtilde H\subseteq G:\exists x\in 
      D_n: x\wtilde H=H\}\big)}
      {\vol(D_{n})}
    \end{equation} 
    exists for all $G\in X$ and is independent of $G\in X$.
\end{indentnummer} 
If any of the above statements applies, then the limit \eqref{freq-cond2} is independent of 
the choice of the tempered van Hove sequence.\\  
Furthermore, the dynamical system  $X_{\mathcal G}$ is \emph{uniquely} ergodic 
iff \itemref{graph-FLC-pattern} holds for 
all patches $H$ of $\cG$ with $X=X_\cG$, that is, for \emph{all} $G\in X_\cG$. 
In that case, the van Hove sequence needs not to be tempered, and the limit \eqref{freq-cond2} 
is independent of the choice of the van Hove sequence. Moreover, the convergence to the
limit in \eqref{freq-cond2} is even uniform in $G\in X_\cG$.
\end{theorem}

\begin{remarks}
\item A proof of the theorem can be found in Section~\ref{sec:proofs-dyn}.
\item Assume that condition \itemref{graph-FLC-pattern} in the above theorem is satisfied, and let $H_1$ and $H_2$ be
equivalent patches of $\cG$. We then have $\nu(H_1)=\nu(H_2)$, compare Lemma~\ref{freq-gen-shift}.
\item Analogously to Proposition~\ref{freq-cond12}, there is also a characterisation of unique ergodicity
in terms on uniform patch frequencies.
\end{remarks}

%%%%%%%%%%%%%%%%%%%%%%%%%%%%%%%%%%%%%%%%%%%%%%%%%%%%%%%%%%%%%%%%%%%%%%%%%%

\subsection{Randomly coloured graphs}

Dynamical systems for coloured graphs are constructed as in
Section~\ref{sec:ran-col}. The only difference is that $\cP_{r}(M)$ will be
replaced by $\cG_{r}(M)$. A \emph{coloured graph} $G^{(\omega{})}$, where
$G\subseteq M$ is a graph and $\omega\in\Omega_{G}$, is given as in
Definition~\ref{col-ps-def}. Copying the proofs of Proposition~\ref{cCr} and
Lemma~\ref{col-T-inv}, we get

\begin{proposition}
  If $\cG\subseteq\cG_r(M)$ is closed in $\cG_r(M)$, then the metrisable space 
  \begin{displaymath}
    \cC_{\cG} := \{G^{(\omega)} : G\in\cG, \omega\in\Omega_{G}\}
  \end{displaymath}
  is closed in $\cG_r(\what{M})$ and hence compact. In particular,
  $\cC_{\cG_{r}(M)}$ and 
  \begin{displaymath}
    \hXcG :=\cC_{X_{\cG}} = \overline{\big\{ xG^{(\omega)}: x\in T, \, G^{(\omega)} \in
      \cC_{\cG}\big\}} 
  \end{displaymath}
  are compact. Moreover, the group action $\alpha_{\hXcG} : T\times\hXcG
  \rightarrow\hXcG$, $(x,G^{(\omega)}) \mapsto x G^{(\omega)}$, which obeys
  \eqref{col-ps-trans}, is continuous. \qed
\end{proposition}

\begin{cor}
  Let $\cG \subseteq \cG_r(M)$ be given. Then 
  $(\hXcG, T, \alpha_{\hXcG})$ is a compact
  topological dynamical system and the Ergodic Theorems~\ref{abstract-erg}
  and~\ref{abstract-unierg} hold for $\cQ =\hXcG$. \qed
\end{cor}

Finally, we turn to randomly coloured graphs and, for given $\cG\subseteq
\cG_{r}(M)$, consider a family of Borel probability measures $\PP_{G}$ on
$(\Omega_{G},\cA_{G})$, which is indexed by $G\in X_{\cG}$.
Assumptions~\ref{p-ass} read exactly the same when formulated for the family
$(\PP_{G})_{G\in X_{\cG}}$. In fact, we refer to this, when we cite
Assumptions~\ref{p-ass} below. Noting Lemma~\ref{proper-trans}, the
Ergodic Theorem~\ref{theo:perg} takes the following form for randomly 
coloured graphs.

\begin{cor}
  \label{ran-G-erg}
  Let $\cG \subseteq \mathcal G_r(M)$ be given. Fix an ergodic Borel probability measure $\mu$ on $X_{\cG}$
  and a family of Borel probability measures $(\PP_{G})_{G\in X_{\cG}}$
  satisfying Assumptions~\ref{p-ass-cov} -- \itemref{p-ass-wcomp}. Then there
  exists a unique ergodic probability measure $\hatmu$ on $\hXcG$ such that
  the following statements hold.
  \begin{indentnummer}
  \item
    For every $f\in L^1(\hXcG,\hatmu)$ we have
    \begin{displaymath} 
      \int_{\hXcG} \!\d\hatmu(G^{(\omega)}) \; f(G^{(\omega)}) =
      \int_{X_{\cG}} \!\d\mu(G) \int_{\Omega_{G}}\!\d\PP_{G}(\omega) \;
      f(G^{(\omega)}).
    \end{displaymath}
  \item 
    For every $f\in L^1(\hXcG,\hatmu)$ and every tempered 
     F\o{}lner sequence $(D_{n})_{n\in\NN}$ in $T$ the limit
    \begin{equation}
      \label{erg-limit-G}
      \lim_{n\to\infty} \frac{1}{\vol(D_{n})} \int_{D_n}\!\d x\;
      f(xG^{(\omega)}) = \int_{\hXcG} \!\d\hatmu(H^{(\sigma)})\;
      f(H^{(\sigma)})   
    \end{equation}
    exists for $\hatmu$-a.a.\ $G^{(\omega)} \in \hXcG$. In fact, the limit
    exists for $\mu$-a.a. $G \in X_{\cG}$ and for $\mathbb P_G$-a.a.\
    $\omega\in\Omega_G$.  
  \end{indentnummer}
	If $X_{\cG}$ is even uniquely ergodic, if $(\PP_{G})_{G\in X_{\cG}}$
  satisfies also Assumption~\ref{p-ass-scomp} and if $f$ is continuous, then
  the limit \eqref{erg-limit-G} exists \emph{for all} $G\in X_{\cG}$ and for
  $\mathbb P_G$-a.a.\ $\omega\in\Omega_G$. In this case the F\o{}lner
  sequence does not need to be tempered.
\end{cor}

When colouring graphs randomly, one may wish to distribute colours on vertices differently from colours on 
edges. This is possible within our framework, as is shown by 

\begin{example}
  Let $\PP_{v}$ and $\PP_{e}$ be two Borel probability measures on $(\AA,
  \cA)$. Given a graph $G \in \cG \subseteq\cG_{r}(M)$, we define
  \begin{displaymath}
    \PP_{G} := \mybigotimes_{m_{v,v} \in G : \; v\in \cV_{G}}\PP_{v} \;\,
    \mybigotimes_{m_{v,w} \in G : \;e=\{v,w\}\in \cE_{G}}\PP_{e}
  \end{displaymath}
  on $(\Omega_{G}, \cA_{G})$, which corresponds to an i.i.d.~distribution of
  colours on vertices and an independent i.i.d.~distribution of colours on
  edges.  Then, the family of measures $(\PP_{G})_{G\in X_{\cG}}$ satisfies
  Assumptions~\ref{p-ass-cov} --~\itemref{p-ass-scomp}. Indeed, $T$-covariance
  and independence at a distance are obvious, and $C$-compatibility can be
  verified as in the proof of Lemma~\ref{iid-ex}. In doing so, we use that the identity
  \eqref{PQ-close} in the proof of that lemma has an analogue in the present
  context of graphs because, if $G,G' \in \cG_{r}(M)$, $m\in G$, $m'\in G'$
  and $d(m,m') < r$, then $m$ and $m'$ are either both vertices or both edges due
  to uniform discreteness.
\end{example}

%%%%%%%%%%%%%%%%%%%%%%%%%%%%%%%%%%%%%%%%%%%%%%%%%%%%%%%%%%%%%%%%%%%%%%%%%%%%%
%%%%%%%%%%%%%%%%%%%%%%%%%%%%%%%%%%%%%%%%%%%%%%%%%%%%%%%%%%%%%%%%%%%%%%%%%%%%%
%
\section{Proofs of results in Section~\ref{psdyn}}
\label{sec:proofs-setup}
%
%%%%%%%%%%%%%%%%%%%%%%%%%%%%%%%%%%%%%%%%%%%%%%%%%%%%%%%%%%%%%%%%%%%%%%%%%%%%%
%%%%%%%%%%%%%%%%%%%%%%%%%%%%%%%%%%%%%%%%%%%%%%%%%%%%%%%%%%%%%%%%%%%%%%%%%%%%%

For the convenience of the reader we have also included proofs of the more
elementary results in Section~\ref{top-point-sets}.

\nlsubsection{Proofs of results in Section~\ref{top-point-sets}}
\label{proofs-top-point-sets}

\begin{proof}[Proof of Lemma~\ref{proper-char}]
  (i) $\Rightarrow$ (ii).~Let $V_{1}, V_{2} \subseteq M$ be compact and observe the identity 
	$S_{V_{1},V_{2}} = \pi_{T} 
  \big(\wtilde{\alpha}^{\,-1}(V_{2} \times V_{1}) \big)$,
	where $\pi_{T}$ stands for the canonical projection $T\times M \rightarrow T$.
	Hence, properness of the map $\wtilde{\alpha}$ and continuity of $\pi_{T}$ yield the claim.
  
  (ii) $\Rightarrow$ (iii).~Let $U_{1}, U_{2} \subseteq M$ be relatively compact. Then, $S_{U_{1}, U_{2}}$ is relatively compact because it is contained in $S_{\overline{U_{1}}, \overline{U_{2}}}$, which is compact by hypothesis.
	
	(iii) $\Rightarrow$ (iv).~ The sets $U_{1} := \{m_{n} \in M : n\in\NN\}$ and $U_{2} := \{x_{n}m_{n} \in M : n\in\NN\}$ are relatively compact in $M$ because the sequences converge. Since 
	$\{x_{n} \in T : n\in\NN\} \subseteq S_{U_{1}, U_{2}}$, and the latter is relatively compact by hypothesis, we infer the claim from the Bolzano-Weierstra{\ss} theorem.

	(iv) $\Rightarrow$ (i).~ Let $V^{(2)}\subseteq M\times M$ be compact. Then $Z:= \wtilde{\alpha}^{\,-1}(V^{(2)})$ is closed in $T\times M$ by continuity of $\wtilde{\alpha}$. Let $\big((x_{n},m_{n})\big)_{n\in\NN}$ be any sequence in $Z$. Then $\big((x_{n}m_{n},m_{n})\big)_{n\in\NN} \subseteq V^{(2)}$, and compactness allows to choose a convergent subsequence $(x_{n_{l}}m_{n_{l}},m_{n_{l}}) \rightarrow (m',m) \in V^{(2)}$ as $l\to\infty$. In particular, every component converges and (iv) yields the existence of a further subsequence $(x_{n_{l_{k}}})_{k\in\NN} \subseteq T$, which converges in $T$. Let us denote its limit by $x$. Since $Z$ is closed, we have $(x_{n_{l_{k}}},m_{n_{l_{k}}}) \rightarrow (x,m) \in Z$ as $k\to\infty$. This proves compactness of $Z$.  
\end{proof}

\begin{proof}[Proof of Lemma~\ref{konchar}]
  (i) $\Rightarrow$ (ii).~This holds by continuity.
  
  (ii) $\Rightarrow$ (iii).~ 
  Fix $m\in M$. \emph{Case 1: assume $m\in P$.} ~
  Let $\varepsilon >0$ and  define $\varphi\in C_c(M)$
  by $\varphi(m')=1-d(m,m')/{\varepsilon}$ for $m'\in B_{\varepsilon}(m)$ and $\varphi(m')=0$ otherwise. 
  Then we have $f_\varphi(P) \ge 1$
  and hence $f_\varphi(P_k)>0$ for finally all $k\in\mathbb N$ by (ii). But this
  means that $ P_k\cap B_\varepsilon(m) \neq \varnothing $ for finally 
all $k\in\mathbb N$.\\
  \emph{Case 2: assume $m\notin P$.}~ Choose $\varepsilon>0$ such that $P\cap B_{2\varepsilon}(m)=\varnothing$. Define $\varphi\in C_c(M)$
  by $\varphi(m')=1-d(m,m')/{2\varepsilon}$ for $m'\in B_{2\varepsilon}(m)$ and $\varphi(m')=0$ otherwise. 
  Then we have $f_\varphi(P)=0$
  and hence $f_\varphi(P_k)<1/2$ for finally all $k\in\mathbb N$ by (ii). But this means
  that $P_k\cap B_\varepsilon(m)=\varnothing$ for finally all $k\in\mathbb N$.
  
  (iii) $\Rightarrow$ (iv).~ Let $P$ be the set of points satisfying condition (iii)-(a) and define
  $N:=M\setminus P$. We first show that  $P\in\mathcal P_r(M)$. To see this, take arbitrary 
  $m\in M$ and assume that $p,q\in P\cap B_r(m)$. Then there exist for finally all $k\in\mathbb N$ 
  points $p_k,q_k\in P_k$ such that $p_k\to p$ and $q_k\to q$ as $k\to\infty$.  But then
  $p_k,q_k\in P_k\cap B_r(m)$ for finally all $k\in\mathbb N$, which implies $p_k=q_k$ for finally all
  $k\in\mathbb N$ due to uniform discreteness of $P_k$. Hence $p=q$ and $\mathrm{card}(P\cap B_r(m))=1$.\\
  Since $V$ is compact, $P_f:=P\cap V$ is finite. For $m\in N$ denote by
  $\varepsilon(m)>0$ a number satisfying $P_k\cap B_{\varepsilon(m)}(m)=\varnothing$ for 
  finally all $k\in\mathbb N$. Now fix $\varepsilon>0$. Compactness of $V$ yields the existence 
  of a finite set $N_f\subseteq N$ such that $ V\subseteq (P_f)_\varepsilon\cup\bigcup_{m\in 
  N_f}B_{\varepsilon(m)}(m)$. We therefore have for finally all $k\in\mathbb N$ the inclusions
  \begin{displaymath}
 P_k\cap V\subseteq P_k\cap\Bigg( (P)_\varepsilon\cup\bigcup_{m\in N_f}B_{\varepsilon(m)}(m)\Bigg)
 =P_k\cap (P)_\varepsilon\subseteq (P)_{\varepsilon},
  \end{displaymath}
  where we used the assumption (iii)-(b) for the equality sign.
  The remaining inclusion follows from $P\cap V=P_f\subseteq (P_k)_\varepsilon$ for finally all $k\in
\mathbb N$ due to assumption (iii)-(a).
  
      (iv) $\Rightarrow$ (i).~
    Let $\varphi_1,\ldots,\varphi_n\in C_c(M)$ and $\varepsilon_1,\ldots,\varepsilon_n>0$ be given.
    For an arbitrary $i\in\{1,\ldots,n\}$, define the compact set $V_i=\mathrm{supp}(\varphi_i)$ and 
denote by 
    $n_i\in\mathbb N$ an upper bound for the
    number of points which a uniformly discrete point set of radius $r$ may have in $V_i$.
    By continuity of $\varphi_i$, we may choose $\delta_i\in]0,r[$ such that we have $|\varphi_i(m)-\varphi_i(m')|<\varepsilon_i/n_i$ for all $m,m'\in V_i$ satisfying $d(m,m')<\delta_i$. This means in particular that $|\varphi_i(m)|<\varepsilon_i/n_i$ for all $m\in V_i$ such that $d(m,V_i^c)< \delta_i$. By assumption (iv), with $\varepsilon=\delta_i$ and 
$V=V_i$, this implies for finally all $k\in\mathbb N$ the estimate
    \begin{displaymath}
   |f_{\varphi_i}(P_k)-f_{\varphi_i}(P)|<\varepsilon_i.
    \end{displaymath}
    Since $i\in\{1,\ldots,n\}$ was arbitrary, this means that for finally all $k\in\mathbb N$ we have $P_k
\in U_{\varphi_1,\varepsilon_1}(P)\cap\ldots\cap
  U_{\varphi_n,\varepsilon_n}(P)$.
  \end{proof}

\begin{proof}[Proof of Proposition~\ref{PrMcompact}]
  Let $(P_n)_{n\in\mathbb N}\subseteq \mathcal P_r(M)$ be given. It suffices to show
  that $(P_n)_{n\in\mathbb N}$ contains a convergent subsequence, since $\mathcal P_r(M)$ 
  is metrisable. 
  
  Since $M$ is $\sigma$-compact, we can find a
  countable open cover $(B_{r}(m_j))_{j\in\mathbb N}$ of $M$ with $m_{j}\in M$ for $j\in\NN$. For
  $j\in\mathbb N$ fixed, consider the sequence $(P_n\cap
  B_{r}(m_j))_{n\in\mathbb N}$. Exactly one of the following two
  cases occurs. Either there is $N_{j}\in\NN$ such that $P_n\cap
  B_{r}(m_j)=\varnothing$ for all $n > N_{j}$, or there is a
  subsequence $(n_k)_{k\in\mathbb N} \subseteq\NN$ such that $
  \varnothing \neq P_{n_k}\cap B_{r}(m_j) =: \{p_{n_{k}}^{(j)}\}$. Due
  to relative compactness of $B_{r}(m_j)$ we assume w.l.o.g.\ that the
  induced point sequence $(p_{n_{k}}^{(j)})_{k\in\NN}$ converges in
  $M$.

  Now, consider the sequence $(P_n\cap B_{r}(m_1))_{n\in\mathbb N}$.
  In the second case of the above scenario, choose a subsequence
  $\big(P_{n^{(1)}_k}\big)_{k\in\mathbb N}$ of $(P_n)_{n\in\mathbb
    N}$ such that the induced point sequence
  $(p_{n_{k}^{(1)}}^{(1)})_{k\in\NN}$ converges to some $p^{(1)} \in
  M$.  Otherwise, set $n_k^{(1)}:=N_{1}+k$ for all $k\in\NN$.  We
  repeat this procedure with the sequence
  $\big(P_{n_{k}^{(1)}}\cap B_{r}(m_{2})\big)_{k\in\mathbb N}$,
  yielding a subsequence $(n_{k}^{(2)})_{k}$ of
  $(n_{k}^{(1)})_{k}$, and then successively for all $j\ge 3$.  In
  this way we obtain nested subsequences $(n_k^{(j+1)})_{k}
  \subseteq (n_k^{(j)})_k$ for all $j\in\NN$. We claim that by
  Cantor's diagonal sequence trick,
  $\big(P_{n^{(k)}_k}\big)_{k\in\mathbb N}$ fulfills the convergence
  criterion of Lemma~\ref{konchar}(iii) and thus converges to
  $P:=\{p^{(j)} \in M: j\in\NN\}$ in the vague topology. Indeed,
  for every $j$ the set $B_{r}(m_{j})\cap \big(\bigcup_{k\in\NN}
  P_{n_{k}^{(k)}}\big)$ is either empty or consists of the points of the convergent sequence
  $(p_{n_{k}^{(k)}}^{(j)})_{k\in\NN}$ with limit $p^{(j)} \in M$.
  Thus, alternative (b) must hold for every $m\in B_{r}(m_{j})$, $m\neq
  p^{(j)}$, otherwise (a) applies.
\end{proof}

\nlsubsection{Proofs of results in Section~\ref {subsec:action-ergodic}}

\begin{proof}[Proof of Lemma~\ref{volest}]
(i). This follows readily from the F{\o}lner property, as 
$(LD_n)\setminus D_n\subseteq \delta^{L} D_n$.

(ii). For $D,E\subseteq T$, it is straightforward to verify
$L((KD)\setminus E) \subseteq LKD\cap LE^c$. This results in the relations 
\begin{equation}\label{van-Hove-imp}
\begin{split}
L\big((KD)\setminus E\big) \subseteq 
\left\{
\begin{array}{cc}
(LKD)\setminus D \cup (LE^c)\setminus D^c\\
(LKD)\setminus E \cup (LE^c)\setminus E^c
\end{array}
\right..
\end{split}
\end{equation}
Consider the first relation in \eqref{van-Hove-imp} for $D=D_n$ and for $E=\mathring{D_n}$. 
This yields
\begin{displaymath}
L\big((KD_n)\setminus \mathring{D_n}\big) \subseteq (LKD_n)\setminus \mathring{D_n} \cup 
(L\overline{D_n^c})\setminus D_n^c,
\end{displaymath}
where we used $(\mathring{D})^c=\overline{D^c}$.
Consider the second relation in \eqref{van-Hove-imp} for $D=\overline{D_n^c}$ and for $E=D_n^c$. 
This yields
\begin{displaymath}
L\big((K\overline{D_n^c})\setminus D_n^c\big) \subseteq (LK\overline{D_n^c})\setminus D_n^c \cup 
(L D_n)\setminus \mathring{D_n}.
\end{displaymath}
When combining these two implications, we obtain
\begin{displaymath}
L(\partial^{K} D_n) \subseteq \partial^{LK}D_n\cup\partial^{L} D_n.
\end{displaymath}
Now the van Hove property yields the claim of the lemma.
\end{proof}

\begin{proof}[Proof of Theorem~\ref{abstract-erg}]
  Until further notice in this proof we only assume that $\cQ$ is a Polish space (i.e., completely metrisable
  with a countable base of the topology), see Remark \ref{general-proof}. The $\mu$-almost-sure existence of 
  the integral \eqref{cesaro-av} follows from Fubini's theorem. To prove \eqref{chaterg} 
  we apply the general pointwise ergodic theorem of Lindenstrauss \cite[Thm.~1.2]{Lind01} for
  tempered F\o{}lner sequences. Since this theorem requires to work with a
  standard probability space (also called Lebesgue space, see \cite[Thm.~2.4.1]{Ito84}), 
  we consider the completed probability space $(\cQ,\bar{\mu})$ with the completed measure 
  $\bar{\mu}$ living on the completion of the Borel $\sigma$-algebra. Then \cite{Lind01} yields the
  existence of a $T$-invariant function $f^{\star} \in L^{1}(\cQ,
  \bar{\mu})$ which obeys $\bar{\mu}(f^{\star}) = \bar{\mu}(f)$ and
  \begin{equation} \label{lind-erg}
    \lim_{n\to\infty} I_n(q,f)= f^\star(q) \qquad\qquad\text{for
      $\bar{\mu}$-a.e.\ $q\in \cQ$.}  
  \end{equation}
  But the limit on the left-hand side of \eqref{lind-erg} is clearly Borel
  measurable in $q$, since $f$ is. Thus, we conclude $f^{\star} \in
  L^{1}(\cQ,\mu)$, $\mu(f^{\star})= \mu(f)$ and that the
  exceptional set in \eqref{lind-erg} can be chosen to be a $\mu$-null set.

  It remains to establish the chain of equivalences.

  (i) $\Rightarrow$ (ii).\quad 
  Since $\mu$ is ergodic, $f^\star$ is $\mu$-a.e.\ constant
  \cite[Thm.~1.3]{BeMa00}. Hence $f^\star(q)=\mu(f^\star)=\mu(f)$
  for $\mu$-a.a.\ $q\in \cQ$.

  (ii) $\Rightarrow$ (iii). \quad
  This is obvious.
  
    (iii) $\Rightarrow$ (i). We are inspired by ideas in \cite[Thm.\ 3.1]{Hof98} and 
    establish first an
  \begin{auxlemma}
    For $k\in\mathbb N$ let $f_k,f\in L^1(\mathcal{Q},\mu)$ be
    given. Assume that  $\|f-f_k\|_1\to 0$ as
    $k\to\infty$ and that $f_k^\star=\mu(f_k)$ holds  $\mu$-a.e.\ for all
    $k\in\mathbb N$. Then, $f^\star =\mu(f)$ holds $\mu$-a.e.
  \end{auxlemma}

  \begin{proof}[Proof of the Auxiliary Lemma]
    Eq.\ \eqref{chaterg}, the triangle inequality, Fatou's Lemma,
    Fubini's theorem and $T$-invariance of $\mu$ provide the
    inequality $\|(f-f_k)^\star\|_1 \le \|f-f_k\|_1$. Thus
    $\|(f-f_k)^\star\|_1 \to 0$ as $k\to\infty$, which in turn implies
    the existence of a subsequence $(k_l)_{l\in\mathbb N}$ such that
    $(f-f_{k_l})^\star\to 0$ pointwise $\mu$-a.e.\ as $l\to\infty$. Now,
    the assertion of the Auxiliary Lemma can be seen from
    \begin{displaymath}
      0 = \lim_{l\to\infty} (f-f_{k_l})^{\star} 
      =  f^\star- \lim_{l\to\infty} \mu(f_{k_l})
      = f^\star-\mu(f),
    \end{displaymath}
    which holds $\mu$-a.e.\ and where the rightmost equality follows
    from $L^{1}$-convergence. 
  \end{proof}

  From now on we assume in addition that $\cQ$ is 
    locally compact. We use the Auxiliary Lemma to establish $\mu$-a.e. 
  \begin{equation}
    \label{open-conv}
    (1_{K})^{\star} = \mu(1_{K})
  \end{equation}
  for indicator functions $1_K$ of compact sets $K\subseteq \cQ$. To see
  this, fix a metric on the metrisable space $\cQ$ that is compatible with the topology,
  and note that by local 
  compactness of $\cQ$, there exists $\varepsilon>0$ such that $(K)_\varepsilon$ 
  is relatively compact. For $n\in\mathbb N$ such that $n\ge 1/\varepsilon$ consider 
  the relatively compact thickened sets $K_n:=(K)_{1/n}$. Using the metric, we define 
  continuous functions $g_n:\cQ\to[0,1]$, such that $g_n=1$ on $K$, and
  $g_n=0$ on $K_n^c$ and $\mathrm{supp}(g_n)\subseteq \overline{K_n}$.
  We thus have $g_n\in C_c(\cQ)$ for all $n\ge N$. We also have 
  $L^1$-convergence of $g_n$ to $1_K$, since
  \begin{displaymath}
  \|g_n-1_K\|_1=\int_{K_n\setminus K} \mathrm{d}\mu(q) \, g_n(q)\le \mu(K_n\setminus K).
  \end{displaymath}
  The latter expression vanishes as $n\to\infty$ by dominated convergence,
  since $\mu(\cQ)=1$, and since closedness of $K$ implies $\lim_{n\to\infty}
  1_{K_n\setminus K}\equiv0$. 
  On the other hand, denseness of $\mathcal{D}$ in $C_c(\cQ)$ with respect to $\|\cdot\|_{\infty}$
  implies denseness with respect to $\|\cdot\|_{1}$ so that we infer
  the existence of a sequence $(f_{n})_{n\in\NN} \subseteq \mathcal{D}$
  with $\| f_n-1_{K} \|_1 \to 0$ as $n\to\infty$. By hypothesis
  of (iii) we also have $\mu$-a.e.\ the equality $f_n^\star=\mu(f_n)$
  for all $n\in\NN$. The Auxiliary Lemma then yields
  \eqref{open-conv}.

  Now local compactness and second countability of $\cQ$ ensure 
  \cite[Thm.~29.12]{Bau01} inner regularity  of the Borel measure  $\mu$, 
  and another application of the Auxiliary Lemma yields $(1_{B})^{\star} =
  \mu(1_{B})$ almost surely for arbitrary Borel sets $B\subseteq \cQ$.
  In particular, for every $T$-invariant Borel set $B\subseteq \cQ$ we
  conclude from this $\mu(B)=1_B(q)$ for $\mu$-a.a.\ $q \in \cQ$.
  Hence, either $\mu(B)=0$ or $\mu(B)=1$, proving (i).

\end{proof}

\begin{proof}[Proof of Theorem~\ref{abstract-unierg}] 
  We adapt the line of reasoning in \cite[Thm.~6.19]{Wa82}. An alternative
  proof can be given using \cite[Thm.~3.5]{Fur81}, compare also \cite[Thm.~3.2]{Sch00}.
  
  The implication (i) $\Rightarrow$ (ii) is obvious.

  (ii) $\Rightarrow$ (iii).\quad Let $\mu_{j}$, $j\in\{1,2\}$, be two $T$-invariant
  Borel probability measures on $\cQ$. The estimate $|I_n(q,f)|\le
  \|f\|_{\infty}$ holds for all $n\in\NN$, $q\in\cQ$ and $f\in
  \mathcal{D}$. This and dominated convergence imply $\lim_{n\to\infty}
  \mu_{j}\big(I_n(\cdot,f)\big) = \mu_{j}\big(I(f)\big) = I(f)$. On the other
  hand, Fubini's theorem yields $\mu_{j}\big(I_n(\cdot,f)\big)=\mu_{j}(f)$ for
  all $n\in\NN$ and $j\in\{1,2\}$. Hence, we get $\mu_{1}(f)=\mu_{2}(f)$ for all
  $f\in \mathcal{D}$. Now, denseness of $\mathcal{D}$ and boundedness of
  $\mu_{j}$ give $\mu_{1}(f)=\mu_{2}(f)$ for all $f\in C(\cQ)$. Thus,
  $\mu_{1}=\mu_{2}$, as both belong to the dual space of $C(\cQ)$.

  (iii) $\Rightarrow$ (i).\quad 
  We prove that (i) holds with $I(f)= \mu(f)$. Suppose, this were
  false. Then there exists $g\in C(\cQ)$, $\varepsilon>0$, a subsequence
  $(n_k)_{k \in\mathbb N} \subseteq\NN$ and a sequence $(q_k)_{k\in\mathbb
    N}\subseteq \cQ$ such that for all $k\in\mathbb N$ we
  have
  \begin{equation}
    \label{contradiction}
    \left|I_{n_k}(q_k,g) -\mu(g)\right|\ge \varepsilon.
  \end{equation}
  On the other hand, for every $k\in\mathbb N$ the linear functional
  $I_{n_k}(q_k,\cdot)$ belongs to the closed unit ball in the dual of the
  Banach space $C(\cQ)$, which is separable, since $\cQ$ is metrisable 
  \cite[Sec.~X.3.3, Thm.~1]{Bour}. The sequential Banach-Alaoglu theorem
  \cite[Thm.~3.17]{Rud91} asserts that this closed unit ball is
  weak*-sequentially compact so that for every $f\in C(\cQ)$ the sequence
  $\big(I_{n_k}(q_k,f)\big)_{k\in\mathbb N}$ contains a convergent
  subsequence. Pick a countable dense subset $\mathcal{C} \subseteq 
  C(\cQ)$ with $g,1\in\mathcal C$. Cantor's diagonal trick 
  gives the existence of a \emph{common} subsequence $(n_{k_{l}})_{l\in\NN}
  \subseteq\NN$ such that $\lim_{l\to\infty} I_{n_{k_{l}}}(q_{k_{l}},f) =: J(f)$ exists for all
  $f\in\mathcal{C}$. Furthermore we have $|J(f)| \le \|f\|_{\infty}$ for
  all $f\in\mathcal{C}$. Thus, $J: \mathcal{C} \rightarrow \RR$,
  $f\mapsto J(f)$, is a bounded linear functional. It is $T$-invariant, due to the
  F{\o}lner property of $(D_n)_{n}$. It admits a unique bounded
  linear extension to $C(\cQ)$, which we denote again by $J$.  This extension
  is positivity preserving, i.e.\ $J(f) \ge 0$, if $f\ge 0$. Therefore the
  Riesz-Markov representation theorem \cite[Thm.~IV.14]{RS80} yields the existence of a positive Borel
  measure $\nu$ on $C(\cQ)$ such that $J(f)=\nu(f)$ for all $f\in C(\cQ)$. But $J(1)=1$,
  which can be seen from the definition of $J$. Moreover, $\nu$ inherits $T$-invariance
  from $J$. Hence, $\nu$ is a $T$-invariant probability measure
  and thus $\nu=\mu$ by uniqueness. But $\nu(g)\ne\mu(g)$ on account of
  \eqref{contradiction}, which is a contradiction.

  So far we have obtained the equivalences (i) to (iii) and that in either
  case the limit $I(f)$ equals $\mu(f)$. In particular, this limit does not
  depend on the chosen F\o{}lner sequence. If $\mu$ is not ergodic,
  there exists a $T$-invariant Borel set $E$ such that $0<\mu(E)<1$. Then a $T$-invariant 
  probability measure $\nu\ne\mu$ is given by $\nu(B):=\mu(B\cap E)/\mu(E)$
  for all Borel sets $B$. Since $\mu$ is the only $T$-invariant probability
  measure, we conclude that $\mu$ is ergodic.
\end{proof}

%%%%%%%%%%%%%%%%%%%%%%%%%%%%%%%%%%%%%%%%%%%%%%%%%%%%%%%%%%%%%%%%%%

\subsection{Proofs of results in Section~\ref{geo-char}}

The Haar measure on $T$ allows to estimate how many
points of a given point set $P\in\mathcal P_r(M)$ fall into some
compact region in $M$. The following statement is a preparation for the proof of
Lemma~\ref{count-compare2}. 

\begin{proposition}\label{standardest}
 Assume that $T$ is even unimodular. Let $(D_n)_{n\in\mathbb N}$ be a 
 F{\o}lner sequence in $T$. Then:
  \begin{indentnummer}
  \item \label{standardest-vol}
  	For every relatively compact set $U\subseteq M$ we have the asymptotic estimate 
    \begin{displaymath}
      \mathrm{card} (P\cap D_n^{-1} U )=\mathocal{O}\big(\vol (D_n)\big) \qquad \qquad
      \text{as} \quad    n\to\infty,
    \end{displaymath}
    uniformly in $P\in\mathcal P_r(M)$.
  \item \label{standardest-bdry}
  	If $(D_n)_{n\in\mathbb N}$ is even a van Hove sequence, we have for 
  every compact set $K\subseteq T$ and for every relatively compact set $U\subseteq M$ the asymptotic estimate
    \begin{displaymath}
      \begin{split}
        &\mathrm{card}\big(P\cap (\partial^{K} D_n)^{-1} U\big)=
        \mathocal{o}(\vol (D_n)) \qquad \qquad
        \text{as} \quad  n\to\infty,
      \end{split}
    \end{displaymath}
    uniformly in $P\in\mathcal P_r(M)$.
  \end{indentnummer}
\end{proposition}

Before we give the proof of the proposition we turn to a useful transformation property of transporters. 

\begin{remark}
 	\label{S-transform}
	In slight abuse of the notation for transporters introduced in \eqref{transporter}, we write 
	$S_{m,U} := S_{\{m\},U} = \{x \in T: xm \in U \}$ for $m\in M$ and a subset $U\subseteq M$. 
	Then, given any group element $x\in T$,
	we observe the identity $S_{xm,U} = S_{m,U} \,x^{-1}$.
\end{remark}

\begin{proof}[Proof of Proposition~\ref{standardest}]
We fix $\varepsilon >0$ and a relatively compact subset $U\subseteq M$. W.l.o.g.\ we assume that $U \neq \varnothing$. 
The open thickened subset $U_{\varepsilon} := (U)_\varepsilon$ is still relatively compact due properness of the metric $d(\cdot, \cdot)$ on $M$.   
We define $\varphi_{\varepsilon}\in C_c(M)$ for $p\in M$ by
$\varphi_{\varepsilon}(p):=d(p,(U_{\varepsilon})^c)$.
For given $p\in M$, the function $x\mapsto \varphi_{\varepsilon}(xp)$ lies in $C_{c}(T)$, since $\varphi_{\varepsilon} \in C_{c}(M)$  and the group action is continuous and proper. In particular, $x\mapsto \varphi_{\varepsilon}(xp)$ is integrable. For $D\subseteq T$ compact we evaluate
\begin{equation}\label{firsteval}
    \int_{D}\!{\d}x \; f_{\varphi_{\varepsilon}}(xP)=
    \int_{D}\!\d x   \sum_{p\in P} \varphi_{\varepsilon}(xp) 
    =   \sum_{p\in P\cap D^{-1}U_{\varepsilon}}
    \int_{D}\!\d x \;\varphi_{\varepsilon}(xp).
\end{equation}
Note that $P\cap D^{-1}U_{\varepsilon}$ is finite, because $P$ is uniformly discrete and $D^{-1}\overline{U_{\varepsilon}}$ is compact. (This argument uses continuity of the group action.)  
Below we wish to approximate the last integral of \eqref{firsteval} by 
\begin{equation}
\label{I-deff}
 I(p):=\int_T\!{\d}x \;\varphi_{\varepsilon}(xp) = \int_{S_{p,U_{\varepsilon}}}\!{\d}x \;\varphi_{\varepsilon}(xp).
\end{equation}
Here, we made use of the transporter $S_{p,U_{\varepsilon}} \subseteq T$, which was introduced in Remark~\ref{S-transform} and is relatively compact by Lemma~\ref{transporterU}. 
It serves to restrict the integration in $I(p)$ to all those arguments where the integrand is strictly positive.
	But first, we will rewrite $I(p)$ for $p\in P \cap D^{-1}U_{\varepsilon}$. For such $p$ there exists $y\in D$ and $m\in U_{\varepsilon}$ such that $p=y^{-1}m$. Hence, Remark~\ref{S-transform} implies
\begin{displaymath}
 S_{p,U_{\varepsilon}} = S_{m,U_{\varepsilon}} \, y \subseteq S_{m,U_{\varepsilon}} D \subseteq L_{U_{\varepsilon}} D,
\end{displaymath}
  where $L_{U_{\varepsilon}} := S_{\overline{U_{\varepsilon}},\, \overline{U_{\varepsilon}}}$ is compact in $T$ by Lemma~\ref{transporterK}. Therefore we have the identity
  \begin{equation}
  	\label{I-rewrite}
    I(p) = \int_{L_{U_{\varepsilon}}D}\!\d x \;\varphi_\varepsilon(xp)
  \end{equation}
  for every $p \in P \cap D^{-1}U_{\varepsilon}$. 
  
  Next we derive a positive lower bound on the integral $I(p)$, which is uniform in $p\in P \cap TU$ 
  (not $U_{\varepsilon}$!). 
  For such $p$ there is $y\in T$ and $q\in U$
  such that $yp=q$. This implies
\begin{displaymath}
I(p) = \int_T\mathrm{d}x\, \varphi_{\varepsilon}(xy^{-1}q)
=\int_T\mathrm{d}x\, \varphi_{\varepsilon} \big((yx)^{-1}q\big) 
=\int_T \mathrm{d}x \,\varphi_{\varepsilon}(xq) 
=I(q),
\end{displaymath}
where we used unimodularity of the group in the second and third equality
and left invariance of the Haar measure in the third equality. We conclude
that for every $p\in P \cap TU$ we have 
\begin{equation}
\label{I-min}
I(p)\ge \inf\{I(q):q\in \overline{U}\} =: I_U >0.
\end{equation}
The strict positivity follows from a compactness argument, using continuity of the map
$q\mapsto I(q)$, and from $I(q)>0$ for all $q\in \overline{U}$. To see the latter we observe $S_{q,U_{\varepsilon}} \supseteq S_{q, B_{\varepsilon/2}(q)}$ for every $q\in \overline{U}$. The transporter  
$S_{q, B_{\varepsilon/2}(q)}$ contains the identity $e\in T$ and is open by continuity of the group action and openness of the ball $B_{\varepsilon/2}(q) \subseteq M$. Therefore there exists an open ball $B^{T} \subseteq T$ about $e$  such that $S_{q, B_{\varepsilon/2}(q)} \supseteq B^{T}$  and $\varphi_{\varepsilon}\big|_{B^{T}} >0$. Since $\vol (B^{T}) >0$ (one can cover the $\sigma$-compact group $T$ by countably many copies of $B^{T}$, all of which have the same Haar measure), it follows that $I(q) >0$.

Next, we establish an auxiliary estimate, which is a consequence of uniform discreteness (of given radius $r$): for every relatively compact subset $U\subseteq M$ there exists a constant $N(U) < \infty$ such that
for every $P \in \cP_{r}(M)$ and every $x\in T$ the bound
\begin{equation}
	\label{Nbound}
 	\mathrm{card}(P\cap xU) \le N(U) < \infty
\end{equation}
holds. To prove this, we first set $x=e$, the identity in $T$, and note that a covering argument then implies \eqref{Nbound} uniformly in $P \in \cP_{r}(M)$. This and the equality $\mathrm{card}(P\cap xU) = \mathrm{card} (x^{-1}P\cap U)$ yield the desired uniformity of $N(U)$ in $x\in T$.

Now, consider the difference of the right-hand side of \eqref{firsteval} 
and the corresponding expression where the integral is replaced by $I(p)$. 
This difference can be estimated 
as
  \begin{align}
  	\label{foldiff}
    \Bigg| \sum_{p \in P \cap D^{-1}U_{\varepsilon}}
    \int_{(L_{U_{\varepsilon}}D)\setminus D}  \!{\d}x \,
    \varphi_{\varepsilon}(xp) \Bigg|
     &\le \int_{(L_{U_{\varepsilon}}D)\setminus D}\!{\d}x  \sum_{p \in
      P \cap D^{-1} U_{\varepsilon}} \;| \varphi_{\varepsilon}(xp)| \\ 
     &\le \int_{L_{U_{\varepsilon}}D}\!{\d}x  \sum_{p \in
      P\cap x^{-1} U_{\varepsilon}} \;| \varphi_{\varepsilon}(xp)| \nonumber \\ 
    & \le \vol (L_{U_{\varepsilon}}D) \, F_{U_{\varepsilon}} ,  \nonumber 
    \end{align}
  where, using \eqref{Nbound}, the constant $F_{U_{\varepsilon}} := N(U_{\varepsilon})\,\|\varphi_{\varepsilon}\|_\infty <\infty $ does not depend on $P$, nor on $D$. Therefore \eqref{I-min} and \eqref{foldiff} imply
  \begin{align*}
  \mathrm{card}(P \cap D^{-1}U)\, I_U  &\le \sum_{p \in P \cap D^{-1}U}  I(p) \le \sum_{p \in
      P \cap D^{-1}U_{\varepsilon}}  I(p)   \\      
   &\le \int_{D}{\d}x \; f_{\varphi_{\varepsilon}}(xP) + 
      	\vol(L_{U_{\varepsilon}}D) \, F_{U_{\varepsilon}}  \\
	 &\le F_{U_{\varepsilon}} \big[ \vol(D) + \vol(L_{U_{\varepsilon}}D) \big].
  \end{align*}
  Thus, the first claim of the proposition follows with $D=D_{n}$ and Lemma~\ref{volestfoelner}, while the second claim follows with $D=\partial^KD_n$, the van Hove property and Lemma~\ref{volesthove}.

\end{proof}

\begin{proof}[Proof of Lemma~\ref{count-compare2}]
Fix any $\varnothing \neq D \subseteq T$ compact,  assume w.l.o.g.\ that $Q\ne\varnothing$ and fix $q\in Q$. Then the number of
points $\wtilde q\in P\cap D^{-1}m$ such that $\wtilde q=xq$ for some $x\in T$
is at most $\mathrm{card}(P\cap D^{-1}m)$. For $\wtilde q\in P\cap D^{-1}m$
we introduce the set
\begin{equation}
\label{Aqq}
A_{q,\wtilde q} := \left\{\wtilde Q\subseteq P : \exists x\in T:xQ=\wtilde Q \;\text{~and~}\;
xq=\wtilde q\right\}.
\end{equation}
Clearly, the estimate
\begin{displaymath}
\card\big(M_{D}^{\prime}(Q)\big)\le \card(P\cap D^{-1}m)
\cdot \max\big\{\card(A_{q,\wtilde q}):\wtilde q\in P\cap D^{-1}m\big\}
\end{displaymath}
holds with the fixed $q\in Q$. In order to estimate the cardinality of $A_{q,\wtilde q}$ for a given $\wtilde q$ 
(assuming $A_{q,\wtilde q}\ne\varnothing$), we fix a reference pattern $\wtilde Q_r\in 
A_{q,\wtilde q}$ and consider an arbitrary $\wtilde Q\in A_{q,\wtilde q}$. Thus there
exist $x_r,x\in T$ such that $x_rQ=\wtilde Q_r$, $xQ=\wtilde Q$, and 
$x_rq=\wtilde q$, $xq=\wtilde q$. The latter imply $x_r^{-1}x\in S_{q, \{q\}}=:S_{q}$, the compact stabiliser group of $q$ by Lemma~\ref{transporterK}. 
In addition, $\wtilde Q=x_rx_r^{-1}xQ\subseteq x_r S_q Q$.
By definition, we have $\wtilde Q\subseteq P$, hence $\wtilde Q\subseteq P\cap x_r S_q Q$
for every $\wtilde Q\in A_{q,\wtilde q}$. 

We conclude from \eqref{Nbound} that $\mathrm{card}(P\cap x_r S_q Q)\le N(\bigcup_{q\in Q}S_qQ)$, uniformly in $P\in\cP_r(M)$ and uniformly 
in $q \in Q$ and in $\wtilde q\in P\cap D^{-1}m$ (which enters through $x_r$).
Therefore there are at most
\begin{displaymath}
F(Q):=\binom{N\big(\bigcup_{q\in Q}S_qQ\big)}{\mathrm{card}(Q)}
\end{displaymath}
possibilities to choose a subset $\wtilde Q$ with $\mathrm{card}(\wtilde Q)=\mathrm{card}(Q)$
points out of the pattern $P\cap x_rS_qQ$. We conclude
\begin{equation}
\label{Aqq-bound}
\mathrm{card}(A_{q,\wtilde q}) \le F(Q)
\end{equation}
uniformly in $q \in Q$ and $\wtilde q\in P\cap D^{-1}m$ and $P\in\cP_r(M)$, and thus
\begin{displaymath}
\mathrm{card}\big(M_{D}^{\prime}(Q)\big) \le F(Q)\,\mathrm{card}(P\cap D^{-1}m).
\end{displaymath}
Hence the second estimate follows with $D=D_{n}$ from Proposition~\ref{standardest-vol}, since $T$ is unimodular 
and the group action is proper.  For the first estimate, we note
\begin{equation}
\label{pat-aux}
M_{D}(Q)\subseteq \big\{\wtilde Q\subseteq P\cap D^{-1}Q:\exists x\in T: xQ=\wtilde Q\big\}.
\end{equation}
Therefore we can argue as above and obtain
\begin{equation} 
\label{pat-card}
\mathrm{card}\big(M_{D}(Q)\big)\le F(Q) \, \mathrm{card}(P\cap D^{-1}Q).
\end{equation}
Since $Q$ is compact, we may now set $D=D_{n}$ and apply Proposition~\ref{standardest-vol}, which uses unimodularity and properness.
 
So it remains to prove that $\mathrm{card}(M_{D}(Q))$ and $\mathrm{card}(M'_{D}(Q))$ differ 
by a $\mathocal{o}(\vol(D))$-term under the specified stronger hypotheses. 
Assume w.l.o.g.~that $Q\ne\varnothing$. By assumption, we may write $Q=Km$ for some non-empty 
finite set $K\subseteq T$. Let $S:=S_m := S_{m,\{m\}}$ denote the stabiliser group 
of $m\in M$, which is compact by Lemma~\ref{transporterK}. We have $S=S^{-1}$ and $Am\cap Bm\subseteq (A\cap BS)m$ for $A,B\subseteq T$ arbitrary.  

Note that $\wtilde Q\in M_{D_n}^{\prime}(Q)\setminus M_{D_n}(Q)$ implies that there
exists $x\in(D_n^{-1})^c$ such that $xQ=\wtilde Q\subseteq P\cap D_n^{-1}m$. Hence we have
\begin{displaymath}
\begin{split}
\wtilde Q&\subseteq D_n^{-1}m\cap (D_n^{-1})^cKm
\subseteq \big(D_n^{-1}\cap (D_n^c)^{-1}KS \big)m\subseteq (\partial^{SK^{-1}}D_n)^{-1}m.
\end{split}
\end{displaymath}
Now the same argument as in (i) yields
\begin{displaymath}
\mathrm{card}\big(M_{D_n}^{\prime}(Q)\setminus  M_{D_n}(Q)\big)\le F(Q)\,
\mathrm{card}\big(P\cap (\partial^{SK^{-1}}D_n)^{-1}m\big),
\end{displaymath}
and the latter term is recognised as $\mathocal{o}(\vol (D_n))$ by Proposition~\ref{standardest-bdry}, 
since $T$ is unimodular and the group action is proper.

Similarly, $\wtilde Q\in 
M_{D_n}(Q)\setminus M_{D_n}^{\prime}(Q)$ implies $\wtilde Q\subseteq D_n^{-1}Q$ 
and $\wtilde Q\not\subseteq D_n^{-1}m$. Thus, there exist $x\in T$ and $q\in Q$ such 
that $xq\in\wtilde Q$ and $xq\in (D_n^{-1}m)^c$, which implies 
\begin{displaymath}
xq \in D_n^{-1}Km \cap(D_n^c)^{-1}m\subseteq \big(D_n^{-1}KS\cap (D_n^c)^{-1}\big)m 
\subseteq (\partial^{SK^{-1}} D_n )^{-1}m.
\end{displaymath}
We have 
\begin{multline*}
M_{D_n}(Q)\setminus  M_{D_n}^{\prime}(Q)\\
\subseteq \left\{\wtilde Q\subseteq P:\exists (x,q)\in T\times Q:xQ=\wtilde Q\wedge xq\in 
(\partial^{SK^{-1}}D_n)^{-1}m\right\} =: A.
\end{multline*}
This set can be represented as 
\begin{displaymath}
A = \bigcup_{q\in Q} \; \; \bigcup_{\wtilde q\in P \cap (\partial^{SK^{-1}}D_n)^{-1}m}  A_{q, {\wtilde{q}}}
\end{displaymath}
with $A_{q, {\wtilde{q}}}$ given by \eqref{Aqq}. Therefore we use \eqref{Aqq-bound} to conclude 
\begin{displaymath}
\mathrm{card}\big(M_{D_n}(Q)\setminus  M_{D_n}^{\prime}(Q)\big)\le F(Q)\,
\mathrm{card}(Q) \, \mathrm{card}\big(P\cap (\partial^{SK^{-1}}D_n)^{-1}m\big),
\end{displaymath}
and the latter term is recognised as $\mathocal{o}(\vol (D_n))$ by Proposition~\ref{standardest-bdry}, since $T$ is unimodular and the group action is proper.
\end{proof}

\begin{proof}[Proof of Lemma~\ref{freq-gen}]
(i). The finiteness of the supremum follows from inequality \eqref{pat-card} and 
Proposition~\ref{standardest-vol}.

(ii). Fix $y\in T$. It suffices to show $\nu^{P}(yQ; (D_{n})_{n\in\NN}) 
= \nu^{P}(Q; (D_{n})_{n\in\NN})$. Since 
$M_{D_{n}}(yQ) = M_{yD_{n}}(Q)$ and 
\begin{displaymath}
 \Big| \mathrm{card}\big(M_{yD_{n}}(Q)\big) - \mathrm{card}\big(M_{D_{n}}(Q) \big) \Big| 
 \le \mathrm{card}\big(M_{A_{n}}(Q)\big) ,
\end{displaymath}
where $A_{n} := \delta^{\{y\}}D_{n} \subseteq \partial^{\{y\}}D_{n}$, the claim follows from inequality \eqref{pat-card}  and Proposition~\ref{standardest-bdry}.
\end{proof}

\begin{proof}[Proof of Lemma~\ref{FLClocal}]
(i). If $X_\cP$ is FLC, then $\cP\subseteq X_\cP$ is FLC by definition. Conversely, assume that $\cP$
is FLC. Take $V\subseteq M$ compact and a corresponding finite set $\cF_{\cP}(V)\subseteq \cQ_{\cP}$ of patterns of $\cP$. Now let $Q$ be any $xV$-pattern of $X_\cP$. Then there is $P\in X_\cP$ such that $Q=P\cap x\mathring V$.
Since $P\in X_\cP$, there is a sequence $((x_n,P_n))_{n\in\mathbb N}\subseteq T\times\cP$ such that
$x_nP_n\to P$ as $n\to\infty$. Hence, for every $n\in\mathbb N$, the pattern 
$\wtilde Q_n:=P_n\cap x_n^{-1}x\mathring V$ is equivalent to some
pattern in $\cF_{\cP}(V)$, and $x_n\wtilde Q_n\to Q$ as $n\to\infty$. Since $\cF_{\cP}(V)$ is finite, 
there is $\wtilde Q\in\cF_{\cP}(V)$, a sequence $(y_k)_{k\in\mathbb N}$ in $T$ and a subsequence $(n_k)_{k\in\mathbb N}$ of $\mathbb N$ such that 
$\wtilde Q_{n_k}=y_k\wtilde Q$ for all $k\in\mathbb N$, implying that $x_{n_k}y_k\wtilde Q\to Q$
as $k\to\infty$. Local compactness of $M$ and properness of the group action imply that a subsequence
of $(x_{n_k}y_k)_{k\in\mathbb N}$ converges to some $z\in T$. Continuity of the group action then yields
$z\wtilde Q=Q$. Thus $z^{-1}Q\in\cF_{\cP}(V)$. 

(ii). 
For patterns $Q,\wtilde Q\in\mathcal P_r(M)$ define $\varepsilon(Q,\wtilde Q)$ by
\begin{displaymath}
\varepsilon(Q,\wtilde Q) := \inf \big\{\delta>0:\exists x\in T: Q\subseteq (x\wtilde Q)_\delta\text{~and~} x\wtilde Q\subseteq (Q)_\delta\big\}.
\end{displaymath}
If $\wtilde Q$ is not equivalent to $Q$, we have $\varepsilon(Q,\wtilde Q)>0$. Indeed,  write 
$Q=\{q_1,\ldots, q_k\}$ and $\wtilde Q=\{\wtilde q_1,\ldots,\wtilde q_k\}$ and assume that
$\varepsilon(Q,\wtilde Q)=0$. Invoking local compactness of $M$, we find a sequence $(x_n)_{n\in\mathbb N}
\subseteq T$ such that we have $x_n\wtilde q_i \to q_i$ for $i\in\{1,\ldots, k\}$ as $n\to\infty$ (possibly 
after some permutation of indices).
Local compactness of $M$ and properness of the group action imply that a subsequence of $(x_n)_{n\in\mathbb N}$ converges to some $x\in T$. By continuity of the group action, we thus get $x\wtilde Q=Q$, which is a contradiction.

Now take an arbitrary $Q\in\cQ_{\cP}$ and define the compact support $V:=\overline{(Q)_r}$ of $Q$ in $M$. Let $\cF_{\cP}(V)$ be a finite set of patterns corresponding to $V$ in the FLC condition and define $\varepsilon>0$ by
\begin{displaymath}
2\varepsilon:=\min\big\{\varepsilon(Q,\wtilde Q):\wtilde Q\in \cF_{\cP}(V) \text{~and~} 
\forall y\in T:Q\ne y\wtilde Q\big\}.
\end{displaymath} 
Now assume that there exist $x\in T$ and $\wtilde Q\in\cQ_{\cP}$ such that 
$x\wtilde Q\subseteq (Q)_\varepsilon$ and $Q\subseteq (x\wtilde Q)_\varepsilon$. 
Then $Q$ is equivalent to $\wtilde Q$, by definition of $\varepsilon$.

(iii). Assume that $\cP$ is not FLC. Then there exist a compact set $V_{0}\subseteq M$ and an infinite
collection $(Q_n)_{n\in\mathbb N}$ of mutually non-equivalent patterns of $\cP$ supported on $T$-shifted copies of $V_{0}$. Due to compactness of $\cP_r(M)$, a subsequence $\wtilde Q_k:=Q_{n_k}$, 
$k\in\mathbb N$, of $(Q_n)_{n\in\mathbb N}\subseteq \cQ_{\cP}$ converges to 
some $Q\in\cP_r(M)$. Let $V\subseteq M$ be a compact set satisfying $Q\subseteq \mathring V$.
Since $\cQ_{\cP}\wedge V$ is closed in $\cP_r(M)$ by assumption, we have $Q\in\cQ_{\cP}\wedge V$,
which implies $Q\in\cQ_{\cP}$. By construction,
we have $\varepsilon(Q,\wtilde Q_k)\to0$ as $k\to\infty$. If $Q$ is not equivalent to
any $\wtilde Q_k$, this contradicts local rigidity of $\cP$. Otherwise, $Q=x\wtilde Q_{\ell}$
for exactly one $\ell$ and some $x\in T$ and $\varepsilon(\wtilde Q_{\ell},\wtilde Q_k)=\varepsilon(Q,\wtilde Q_k)\to0$ as $k\to\infty$. This is also a contradiction to local rigidity of $\cP$, since
$\wtilde Q_{\ell}$ is not equivalent to $\wtilde Q_k$ for $k\ne\ell$.
\end{proof}

Before we can approach auxiliary results for the proof of  
Theorem~\ref{uniFLC} we introduce some more systematic notation for pattern collections of point sets.

\begin{definition}
	\label{pat-col}
 	Let $U \subseteq M$ and $D \subseteq T$. For $P \in \cP_{r}(M)$ we define 
 	\begin{displaymath}
 		\cQ_{P}(U;D) := \big\{\, Q \subseteq P: \exists x \in D  \text{~such that~}
		xQ \subseteq  U\big\}
	\end{displaymath}
	and, in slight abuse of notation for $\cP \subseteq \cP_{r}(M)$,
	\begin{displaymath}
		\cQ_{\cP}(U;D) := \bigcup_{P \in\cP} \cQ_{P}(U;D).
	\end{displaymath}
	It will also be convenient to fix in addition the number of points $k\in\NN$ of the patterns, 
	in symbols,
	\begin{displaymath}
		\cQ^{k}_{P}(U;D) :=\big\{\, Q \in \cQ_{P}(U;D) : \card(Q) =k \big\},
	\end{displaymath}
	and similarly $\cQ^{k}_{\cP}(U;D)$. In particular, we have $M_{D}(Q) = \cQ_{P}^{\card(Q)}(Q;D)$ 
	for every given point set $P$. 	We write $\cQ^{k}_{\cP}(U) := \cQ^{k}_{\cP}(U;T)$ 
	and $\cQ^{k}_{\cP}:= \cQ^{k}_{\cP}(M)$, which has already 
	been used before. 
\end{definition}

The next lemma and the subsequent proposition will be needed in the course of the proof of Theorem~\ref{uniFLC}. 
But they also enter the proof of Theorem~\ref{lem:hoflemma}, which is the main ingredient for the ergodic theorem of randomly coloured point sets. 

\begin{lemma}
 	\label{pattern-count}
	Given $P \in \cP_{r}(M)$ 
	and a relatively compact subset $U \subseteq M$, there exists a constant $\Gamma_{U} >0$ 
	which depends only on $U$ and the radius of relative discreteness $r$ -- but not on 
	$P \in \cP_{r}(M)$ -- such that for every 
	$k \in\NN$ and every compact subset $D \subseteq T$ the estimate 
	\begin{displaymath}
 		\card\big( \cQ^{k}_{P}(U;D)\big) \le \Gamma_{U}^{k-1} \card(P \cap D^{-1}U) 
	\end{displaymath}
	holds. 
\end{lemma}

\begin{proof}
 	We start by observing $\cQ^{k}_{P}(U;D) \subseteq \bigcup_{q \in P \cap D^{-1}U} A_{q}$,
	where
	\begin{align*}
 		A_{q} :=& \,\Big\{ Q \subseteq P  : \;\card(Q)=k, \,
		q \in Q \text{~~and~~} \exists \, x \in T 
		  \text{~~such that~~} xQ \subseteq U\Big\} \nonumber\\
		= & \, \Big\{ \{q,q_{2}, \ldots, q_{k}\} \subseteq P : \;\exists\, x \in S_{q,U} 
		\text{~~with~~}
		x q_{i} \in U \,\forall i \in \{2,\ldots, k\}\Big\} \nonumber\\
		\subseteq & \, \Big\{ \{q,q_{2}, \ldots, q_{k}\} \subseteq P :\;  q_{i} \in P \cap S_{q,U}^{-1}U \,\forall i \in \{2,\ldots, k\}\Big\}.
	\end{align*}
	This implies 
	\begin{align*}
 	  \card\big( \cQ^{k}_{P}(U;D)\big) &\le \sum_{q \in P \cap D^{-1}U}	
	  	\Big[\card\big( P \cap S_{q,U}^{-1} U\big) \Big]^{k-1} \nonumber\\
		& \le \card(P \cap D^{-1}U) \, \bigg[ \sup_{q \in P \cap TU}  
		\card\big( P \cap S_{q,U}^{-1} U \big) \bigg]^{k-1}.
	\end{align*}
	For every $q \in P \cap TU$ there exists $x_{q} \in T$ and $m_{q}\in U$ such that $q= x_{q}m_{q}$. 
	Thus, we conclude from the transformation property of transporters in Remark~\ref{S-transform} that
	\begin{align*}
 		\card\big( P \cap S_{q,U}^{-1} U \big) &= \card\big( P \cap x_{q}S_{m_{q},U}^{-1} U \big) 
		\le \card\big( x_{q}^{-1} P \cap S_{U,U}^{-1} U \big)  \nonumber\\
		&\le N \big( S_{U,U}^{-1} U\big),
	\end{align*}
	where the last inequality rests on \eqref{Nbound} and holds uniformly in $P \in\cP_{r}(M)$ and $x_{q} \in T$, and therefore uniformly in $q \in P \cap TU$. 
	Here, the application of \eqref{Nbound} is justified, because $S_{U,U}^{-1} U$ is relatively compact
	in $M$. This follows from Lemma~\ref{transporterU} and continuity of the group action. 
	So the claim holds with $\Gamma_{U}= N ( S_{U,U}^{-1} U)$.
\end{proof}

We write $L^{0}_{b,c}(M)$  to denote the set of all  real-valued, Borel-measurable and  bounded functions 
$\varphi$ on $M$, whose set-theoretic support $\{m\in M:\varphi(m)\ne0\}$ is relatively compact. 
For $\varphi\in L^{0}_{b,c}(M)$, we consider $f_\varphi:\cP_r(M)\to\mathbb R$ as in Definition~\ref{fphi}.

\begin{proposition}\label{comp}
Let $(D_n)_{n\in\mathbb N}$ be a F{\o}lner sequence in $T$. Fix $k\in\mathbb N$ and consider functions $\varphi_i\in L^{0}_{b,c}(M)$, $i\in\{1,\ldots,k\}$, whose set-theoretic supports $U_i$, $i\in\{1,\ldots,k\}$, are relatively compact and pairwise disjoint. Let $U := \bigcup_{i=1}^{k} U_{i}$. Then we have the equality
\begin{displaymath}
\int_{D_n}{\d}x \; \bigg(\prod_{i=1}^k
f_{\varphi_i}\bigg)(xP) = \sum_{Q \in
\cQ_{P}^{k}(U;D_{n})}  I(Q)  +
\mathocal{o}\big(\vol (D_n)\big),
\end{displaymath}
asymptotically as $n\to\infty$. Here the $\mathocal{o}(\vol (D_n))$-term can be chosen uniformly in $P\in \cP_{r}(M)$ and the leading term
\begin{equation}
  \label{I-def}
	I(Q) :=  \sum_{\pi\in\cS_{k}} \int_{T}\!\d x\, \prod_{i=1}^k \varphi_i(xq_{\pi(i)})
\end{equation}
involves a sum over all permutations from the symmetric group $\cS_{k}$ so that the fixed choice for enumerating the points of the pattern $Q= \{ q_{1}, \ldots, q_{k}\}$ is irrelevant.
\end{proposition}

\begin{proof}
We fix $P \in\cP_{r}(M)$ arbitrary. Note first that, for $p\in M$ and $i\in\{1,\ldots,k\}$ fixed, the function $x\mapsto \varphi_i(xp)$ is 
integrable, since $\varphi_i$ is measurable, bounded, has a relatively compact support, and since the group action is continuous and proper. Hence, $x \mapsto \prod_{i=1}^k \varphi_i(xq_{\pi(i)})$ is integrable, too. Moreover, since the supports
$U_{i}$, $i\in\{1,\ldots,k\}$, are pairwise disjoint, we have 
\begin{equation}
	\label{factorout}
 	\bigg(\prod_{i=1}^k f_{\varphi_i}\bigg)(xP) 
	= \prod_{i=1}^k \bigg(\sum_{p \in P}  \varphi_i(xp)\bigg)
	= \sum_{Q \in \cQ_{P}^{k}(U;D_{n})} \sum_{\pi\in\cS_{k}}
	\prod_{i=1}^k \varphi_i(xq_{\pi(i)})
\end{equation}
for every $x\in D_{n}$. By Lemma~\ref{pattern-count} the set $\cQ_{P}^{k}(U;D_{n})$ is finite, and  
integrating \eqref{factorout} gives
\begin{equation}
\label{eval2-start}
    \int_{D_n}\!{\d}x \; \bigg(\prod_{i=1}^k
    f_{\varphi_i}\bigg)(xP)
    =   \sum_{Q \in \cQ_{P}^{k}(U;D_{n})} \;\sum_{\pi\in\cS_{k}} \int_{D_{n}}\!\d x\, \prod_{i=1}^k \varphi_i(xq_{\pi(i)}).
\end{equation}
Now we wish to replace the sum over permutations on the right-hand side  of \eqref{eval2-start} by $I(Q)$ asymptotically as $n\to\infty$. This is achieved in analogy to the argument leading from \eqref{I-deff} to \eqref{I-rewrite}: we start with the observation 
\begin{equation}
\label{IQ-start}
 \int_{T}\!\d x\, \prod_{i=1}^k \varphi_i(xq_{\pi(i)})
  = \int_{S(Q)}\!\d x\, \prod_{i=1}^k \varphi_i(xq_{\pi(i)}),
\end{equation}
where we introduced
  \begin{displaymath}
    S(\what{Q}) := \{x\in T: x\what{Q}  \subseteq U\} = \bigcap_{\what{q} \in \what{Q}} S_{\what{q},U}\subseteq T
  \end{displaymath}
  for general $\what{Q} \subseteq M$. In this way we excluded parts of the domain of integration in \eqref{IQ-start} where the integrand vanishes anyway. In the special case $\what{Q} \subseteq U$ we have 
  $S(\what{Q}) \subseteq S_{\overline{U}, \, \overline{U}} =: L_{U}$, which is compact by Lemma~\ref{transporterK}. 
Next suppose that $Q\subseteq D_{n}^{-1}U$ (as is the case for $Q \in \cQ_{D_{n}}^{k}(P)$). Then there exists $y \equiv y(Q)\in D_{n}$ and $\what{Q}\subseteq U$ such that $Q=y^{-1}\what{Q}$. Hence we conclude from Remark~\ref{S-transform} that 
\begin{align}
	\label{S-reduce-1}
 	S(Q) = S(y^{-1}\what{Q}) = S(\what{Q})y &\subseteq L_{U}y \\
	\label{S-reduce}
	&\subseteq L_{U} D_{n}.
\end{align}
Therefore \eqref{IQ-start} and \eqref{S-reduce} yield the identity
  \begin{equation}
  	\label{IQ-rewrite}
   I(Q)
   = \sum_{\pi\in\cS_{k}}\int_{L_{U}D_{n}}\!\d x \;\prod_{i=1}^k \varphi_i(xq_{\pi(i)})
  \end{equation}
  for every $Q \in \cQ_{D_{n}}^{k}(P)$. From \eqref{IQ-rewrite} and \eqref{eval2-start} we deduce the estimate
  \begin{displaymath}
  \begin{split}
		\Bigg|\sum_{Q \in \cQ_{P}^{k}(U;D_{n})} \!\! &I(Q) \; - \int_{D_n}\!{\d}x \; \bigg(\prod_{i=1}^k
    f_{\varphi_i}\bigg)(xP)  \Bigg|   \\  
    & \le \sum_{\pi\in\cS_{k}}
    \int_{(L_{U}D_n)\setminus D_n}  \!{\d}x \;  \sum_{Q \in \cQ_{P}^{k}(U;\{x\})} \prod_{i=1}^k
    |\varphi_i(xq_{\pi(i)})| \\
    & \le  k! \,\vol \big((L_{U}D_n)\setminus D_n\big) \, 
    \sup_{x\in T} \Big[\card\big( \cQ_{P}^{k}(U;\{x\})\big) \Big] \prod_{i=1}^k
    \|\varphi_{i}\|_\infty.
    \end{split}
  \end{displaymath}
In order to get the first inequality above we have used the identity $\sum_{Q \in \cQ_{P}^{k}(U;D_{n})} \prod_{i=1}^k \varphi_i(xq_{\pi(i)}) = \sum_{Q \in \cQ_{P}^{k}(U;\{x\})} \prod_{i=1}^k
  \varphi_i(xq_{\pi(i)})$, which holds for every fixed $x\in T$.  
The assertion of the proposition now follows from the F\o{}lner property \eqref{def-folner}, the fact that $L_{U}$ is independent of $P$ and the estimate 
	\begin{displaymath}
	\card\big(\cQ_{P}^{k}(U;\{x\})\big) \le \Gamma_{U}^{k-1} \card(P \cap x^{-1}U) 
	\le  \Gamma_{U}^{k-1} N(U) < \infty,
\end{displaymath}
which is based on Lemma~\ref{pattern-count} and holds uniformly in $P\in \cP_{r}(M)$ and $x \in T$ by \eqref{Nbound}. 
\end{proof}

The following proposition refines the asymptotic evaluation of Proposition~\ref{comp}
in terms of pattern frequencies. For that reason, $T$ is assumed to be unimodular, and the FLC 
assumption is imposed.

\begin{proposition}\label{comp-freq}
Let $(D_n)_{n\in\mathbb N}$ be a van Hove sequence in the unimodular group $T$. Fix $k\in\mathbb N$ and consider 
functions $\varphi_i\in L^{0}_{b,c}(M)$, $i\in\{1,\ldots,k\}$, 
whose set-theoretic supports $U_i$, $i\in\{1,\ldots,k\}$, 
are relatively compact and pairwise disjoint. Set $U:= \bigcup_{i=1}^{k} U_{i}$. 
Let $\cP\subseteq \cP_r(M)$ be of finite local complexity and let $\cF_{X_{\cP}}^{k}(U)$ be
a maximal subset of mutually non-equivalent patterns in $\cQ^{k}_{X_{\cP}}(U)$.
Then we have for every $P\in X_\cP$ the asymptotic estimate
\begin{displaymath} %\label{finalres}
	\int_{D_n}\!\d x\;  \bigg(\prod_{i=1}^k 
	f_{\varphi_i}\bigg)(xP)=\sum_{Q\in\cF_{X_{\cP}}^{k}(U)} I(Q) \, 
	\card\big(M_{D_n}(Q)\big) + \mathocal{o}\big(\vol (D_n)\big)
\end{displaymath}
as $n \to\infty$. Here, the finite set $\cF_{X_{\cP}}^{k}(U)$ and the integral $I(Q)$ 
are independent of the particular choice of $P\in X_\cP$, and
the error term can be chosen uniformly in $P\in X_{\cP}$.
\end{proposition}

\begin{proof}
Fix $P\in X_{\cP}$. By Proposition~\ref{comp}, we have
  \begin{equation}\label{eq:countsum}
    \int_{D_n}{\d}x \; \bigg(\prod_{i=1}^k
    f_{\varphi_i}\bigg)(xP) = \sum_{\wtilde Q\in \cQ^k_{P}(U;D_{n})}  I(\wtilde Q)  +
     \mathocal{o}\big(\vol (D_n)\big),
  \end{equation}
asymptotically as $n\to\infty$, where the error term can be chosen uniformly in $P\in X_{\cP}$. 
In order to establish a connection to pattern frequencies, we partition the set 
$\cQ^k_{P}(U;D_{n})$ into subsets of equivalent patterns. Due to FLC of $X_\cP$,  
cf.~Lemma~\ref{FLClocal}, the set $\cF^{k}_{X_{\cP}}(U)$ is finite. 
Given an arbitrary pattern $Q\in \cF^{k}_{X_{\cP}}(U)$ we consider the collection 
$\cQ^{k}_{P\cap D_{n}^{-1}U}(Q) \subseteq  \cQ^{k}_{P\cap D_{n}^{-1}U} = \cQ_{P}^{k}(U;D_{n})$ of all its translates in $P\cap D_{n}^{-1}U$.
Then the sum in \eqref{eq:countsum} decomposes 
  \begin{equation}
  	\label{countsum-decomp}
    \int_{D_n}{\d}x \; \bigg(\prod_{i=1}^k
    f_{\varphi_i}\bigg)(xP) = \sum_{Q\in \cF^{k}_{X_{\cP}}(U)} \;\sum_{\wtilde Q\in \cQ^{k}_{P\cap D_{n}^{-1}U}(Q)}
       \!\!I(\wtilde Q)  +
     \mathocal{o}\big(\vol (D_n)\big). %\nonumber\\[-2ex]
  \end{equation}
But the integral $I( \wtilde Q)$ is independent of the particular choice of $\wtilde Q\in \cQ^{k}_{P\cap D_{n}^{-1}U}(Q)$, as we show now: by definition there exists $y=y(\wtilde Q)\in T$ and enumerations of the points in the two patterns $Q=\{q_1,\ldots,q_k\}$ and $\wtilde Q=\{\wtilde q_1,\ldots,\wtilde q_k\}$ such that $y \wtilde q_i=q_i$ for all $i\in\{1,\ldots, k\}$. Then we get
\begin{align}
	\label{int-inv}
	I(\wtilde Q) &= \sum_{\pi\in \mathcal{S}_k}\int_T\mathrm{d}x\, \prod_{i=1}^k \varphi_i 
		\left(x \wtilde q_{\pi(i)}\right)
	 = \sum_{\pi\in \mathcal{S}_k}\int_T\mathrm{d}x\, \prod_{i=1}^k\varphi_i
	 \left((yx)^{-1} q_{\pi(i)}\right)  \\
 &  = \sum_{\pi\in \mathcal{S}_k}\int_T \mathrm{d}x \,\prod_{i=1}^k \varphi_i\left(x q_{\pi(i)}\right) 
  = I(Q), \nonumber
\end{align}
where we used unimodularity of the group for the second and third equality
and left invariance of the Haar measure for the third equality.

In order to analyse the cardinality of $\cQ^{k}_{P\cap D_{n}^{-1}U}(Q)$ for $\card(Q) =k$, consider the set
\begin{displaymath}
S:=\{x\in T: x Q\subseteq U\} = \bigcap_{i=1}^k S_{q_{i}, U},
\end{displaymath}
which is relatively compact in $T$ by Lemma~\ref{transporterU}. 
Then we claim
\begin{displaymath}
\cQ^{k}_{P\cap D_{n}^{-1}U}(Q)=
\big\{\wtilde Q\subseteq P: 
\exists y\in S^{-1}D_n \text{~with~}\, y\wtilde Q= Q\big\}=M_{S^{-1}D_n}( Q).
\end{displaymath}
Indeed, to verify the inclusion $\cQ^{k}_{P\cap D_{n}^{-1}U}(Q) \subseteq M_{S^{-1}D_n}(Q)$, take 
$\wtilde Q\in \cQ^{k}_{P\cap D_{n}^{-1}U}(Q)$ and choose $x\in D_n$ and $y\in T$ such that 
$x\wtilde Q \subseteq U$ and $y\wtilde Q = Q$. Then we have  
$xy^{-1}\in S$. But this means that $y \in S^{-1}D_{n}$, whence $\wtilde Q\in M_{S^{-1}D_n}(Q)$. For the reverse inclusion, 
take $\wtilde Q\in M_{S^{-1}D_n}(Q)$ and choose $y\in S^{-1}D_n$ with $y\wtilde Q= Q$. Then 
$y=s^{-1}x$ for some $s\in S$ and some $x\in D_n$. Hence $x\wtilde Q=s Q\subseteq U$. 
This means that $Q\in \cQ^{k}_{P\cap D_{n}^{-1}U}(Q)$.

But the sets $M_{S^{-1}D_n}(Q)$ and $M_{D_n}(Q)$ are asymptotically 
of the same cardinality. This can be seen from
\begin{displaymath}
\begin{split}
M_{S^{-1}D_n}(Q)\, \sydi \, M_{D_n}(Q)&=\Big\{\wtilde Q\subseteq P: 
\exists x\in \big(\delta^{S^{-1}}D_n\big)^{-1}: x Q= \wtilde Q \Big\}\\
&\subseteq \Big\{\wtilde Q\subseteq P\cap  \big(\partial^{S^{-1}}D_n\big)^{-1} Q:\exists x\in T: x Q= \wtilde Q\Big\},
\end{split}
\end{displaymath}
where $\sydi$ denotes the symmetric difference: we argue as in the proof of Lemma~\ref{count-compare2}, compare Eqs.~\eqref{pat-aux} -- \eqref{pat-card}, to show
\begin{displaymath}
\card\left(M_{S^{-1}D_n}(Q)\,\sydi\, M_{D_n}(Q)\right)\le F(Q) \card \Big( P\cap  \big(\partial^{S^{-1}}D_n\big)^{-1}  Q\Big).
\end{displaymath}
A final appeal to Proposition~\ref{standardest-bdry} yields
\begin{displaymath}
\card\big( \cQ^{k}_{P\cap D_{n}^{-1}U}(Q) \big) = 
\card\big(M_{S^{-1}D_n}(Q)\big)=\card\big(M_{D_n}(Q)\big) + \mathocal{o}\big(\vol (D_n)\big)
\end{displaymath}
as $n\to\infty$, where the error term can be chosen uniformly in $P\in X_{\cP}$. This holds by the van Hove property of $(D_n)_{n\in\mathbb N}$, where we used unimodularity and properness of the group action. Thus, the claim follows together with \eqref{countsum-decomp} and \eqref{int-inv}. 
\end{proof}

\begin{proof}[Proof of Theorem~\ref{uniFLC}]
Let $\mu$ be a $T$-invariant Borel probability measure on $X_{\cP}$.
We first prove the asserted characterisation of ergodicity of $\mu$.  Our arguments
rely on Theorem~\ref{abstract-erg}, which requires a tempered F\o{}lner sequence. In addition, the van Hove property enters through Proposition~\ref{comp-freq}.

{\rm (i)} $\Rightarrow$ {\rm (ii)} W.l.o.g.~fix a non-empty pattern $Q=\{q_1,\ldots,q_k\}$, 
$k\in\mathbb N$, of $\cP$. By FLC of $X_\cP$, cf. Remarks~\ref{remFLC}, we may choose $\varepsilon\in ]0,r/2[$ such 
that all patterns of $X_\cP$ of cardinality $k$, which admit a support equivalent to the compact set $\overline{(Q)_\varepsilon}$, 
are equivalent to $Q$. For $i\in\{1,\ldots,k\}$ define the mutually disjoint sets $U_i:={B_\varepsilon(q_i)}$. Choose $\varphi_i\in C_c(M)$ of compact support $\overline{U_{i}}$ for $i\in\{1,\ldots,k\}$ and consider the function $f:=f_{\varphi_1}\cdot\ldots\cdot f_{\varphi_k}\in C(X_{\cP})$. 
Setting $U:= \bigcup_{i=1}^{k} U_{i}$, we can now apply Proposition~\ref{comp-freq} with $\cF^{k}_{X_{\cP}}(U) = \{Q\}$. This yields 
\begin{equation}
	\label{int-freq}
 	\int_{D_n}\!\d x \, f(xP) =  I(Q) \, \card\big(M_{D_{n}}(Q)\big) + \mathocal{o}\big(\vol(D_{n})\big),
\end{equation}
where $P$ enters only through $M_{D_{n}}(Q)$ on the right-hand side. 
Since $\mu$ is ergodic and $f\in L^1(X_{\cP},\mu)$, Theorem~\ref{abstract-erg}~(ii) guarantees the existence of a set 
$X\subseteq X_{\cP}$ of full $\mu$-measure such that for all $P\in X$ we have
\begin{displaymath}
\lim_{n\to\infty} \frac{1}{\vol (D_n)}\int_{D_n}\d x \, f(xP)=\mu(f),
\end{displaymath}
and this limit is independent of $P\in X$. Hence condition (ii) of the Theorem is satisfied.

{\rm (ii)} $\Rightarrow$ {\rm (i)} We will apply the characterisation of ergodicity in 
Theorem~\ref{abstract-erg}~(iii). First, we define a suitable $\|\cdot\|_{\infty}$-dense 
subset $\mathcal{D}$ of $C(X_\cP)$. It will be constructed from the set 
\begin{displaymath}
    \mathcal{D}_{0}:=\left\{f_\varphi:\varphi\in C_c(M), \mathrm{diam}\big(\mathrm{supp}
    (\varphi)\big)<r/2 \right\} \cup\left\{1\right\},
\end{displaymath}
where $1\in C(X_\cP)$ denotes the constant function equal to one, and with $f_\varphi$ 
as in Definition~\ref{fphi}. The set $\mathcal D_0$
separates points in $X_\cP$. Hence the Stone-Weierstra\ss{} theorem
\cite[Prob.~7R]{Kel75} assures that the algebra
$\mathcal{D}:=\mathrm{alg}(\mathcal{D}_{0})$ generated by 
$\mathcal{D}_{0}$ is dense in $C(X_\cP)$ with respect to the supremum norm.

W.l.o.g.~consider $f\in\mathcal D$ of the form
$f=f_{\varphi_1}\cdot\ldots\cdot f_{\varphi_k}$, where $k\in\mathbb N$ and 
$f_{\varphi_j}\in \mathcal{D}_{0}\setminus \{1\}$ for $j\in\{1,\ldots,k\}$. Write $V_i:=
\mathrm{supp}(\varphi_i)$ for the compact supports of the functions $\varphi_i$ 
for $i\in\{1,\ldots,k\}$. Note that $V_i\cap V_j\ne\varnothing$ implies that 
$f_{\varphi_i}\cdot f_{\varphi_j}=f_{\varphi_i\cdot\varphi_j}$. We thus assume w.l.o.g.
that $V_i\cap V_j=\varnothing$ for $i\ne j$. Write $U_i \subseteq V_i$ for the
set-theoretical support of the function $\varphi_i$, for $i\in\{1,\ldots,k\}$, and define
$U:=\bigcup_{i=1}^k U_i$.

Theorem~\ref{abstract-erg} guarantees the existence of $X'\subseteq X_{\cP}$ of full $\mu$-measure 
and of a $T$-invariant function
$f^\star\in L^1(X_{\cP},\mu)$ such that $\mu(f^\star)=\mu(f)$, and such that for all $P\in X'$ we have
\begin{equation}\label{ergshow}
\lim_{n\to\infty} \frac{1}{\vol (D_n)}\int_{D_n}\d x \,  f(xP)=f^\star(P).
\end{equation}
In order to show that the right-hand-side of \eqref{ergshow} % evaluates to '(5.26)' 
is constant in $P$, we will evaluate the left-hand-side of \eqref{ergshow} % evaluates to '(5.26)' 
using Proposition~\ref{comp-freq}. %evaluates to '5.6'  
To do so, we note that $\cF^{k}_{X_{\cP}}(U)$ is a finite set. Thus, there exists 
a set $X''\subseteq X_{\cP}$ of full $\mu$-measure such that hypothesis (ii) is satisfied 
for all $Q\in \cF^{k}_{X_{\cP}}(U)$ and for all $P\in X''$. 
Then the set $X:=X'\cap X''$ has full $\mu$-measure (and is in particular non-empty), 
and equation \eqref{ergshow}  % evaluates to '(5.26)'  
holds for all $P\in X$. Now, Proposition~\ref{comp-freq}. %evaluates to '5.6'
and hypothesis (ii) imply that the value $f^\star(P)$ is indeed independent of $P\in X$, 
which in turn yields $f^\star=\mu(f^\star)=\mu(f)$ on $X$. 
Therefore $\mu$ is ergodic. 
The asserted independence of the pattern frequency of the choice of the tempered van Hove sequence
follows from the corresponding independence in the Ergodic Theorem~\ref{abstract-erg}.

In the simpler case of unique ergodicity one can argue as above, now with an arbitrary
van Hove sequence. To prove {\rm (i)} $\Rightarrow$ {\rm (ii)},
one uses Theorem~\ref{abstract-unierg} (i). To prove {\rm (ii)} $\Rightarrow$ {\rm (i)}, one can apply
the characterisation of unique ergodicity in Theorem~\ref{abstract-unierg} (ii). Independence
of the choice of the van Hove sequence follows from Theorem~\ref{abstract-unierg}.

If $X_\cP$ is uniquely ergodic, then the convergence to $\nu(Q)$ in Definition~\ref{pattcount} is even 
uniform in $P\in X_\cP$, since, after dividing by $\vol(D_{n})$ in \eqref{int-freq}, the convergence on the left-hand side is uniform in $P\in X_\cP$ by 
Theorem~\ref{abstract-unierg}~(i), and since the error term can be chosen uniformly in $P\in X_\cP$.
\end{proof}

\begin{proof}[Proof of Proposition~\ref{freq-cond12}]
Let $(D_{n})_{n\in\NN}$ be a van Hove sequence in $T$.

(i) $\Rightarrow$ (ii). Note first that uniform convergence of $\nu(y,P)$ in $(y,P)\subseteq T\times\cP$,
with a limit independent of $(y,P)$, is equivalent to the existence of the limit 
\begin{equation}\label{freq-con-alt}
\lim_{n\to\infty}\nu_n^{y_n,P_n}(Q)=\lim_{n\to\infty}\frac{\mathrm{card}\big(\big\{\wtilde Q\subseteq P_n:
\exists x\in D_ny_n: x\wtilde Q=Q\big\}\big)}{\vol (D_n)}
\end{equation}
for every sequence $\left((y_n,P_n)\right)_{n\in\mathbb N}\subseteq T\times 
\cP$, with independence of the limit of $((y_n,P_n))_{n\in\mathbb N}$. Assume now that 
$X_\cP$ is uniquely ergodic and fix a pattern $Q \in\cQ_{\cP}$. Then condition \eqref{freq-con-alt} is satisfied, because for every sequence 
$((y_n,P_n))_{n\in\mathbb N}\subseteq T\times\cP$ we have
\begin{displaymath}
	\nu_n^{y_n,P_n}(Q)
	= \frac{\card\big(\big\{\wtilde Q\subseteq y_nP_n:\exists x\in D_n: x\wtilde Q=Q\big\}\big)}
		{\vol (D_n)} =\nu_n^{e,y_nP_n}(Q),
\end{displaymath}
and because the convergence in the limit underlying the definition of $\nu(Q)$ is uniform in 
$P\in X_\cP$ by unique ergodicity of $X_\cP$, see the second part of Theorem~\ref{uniFLC}.
Hence $\cP$ has uniform pattern frequencies.

(ii) $\Rightarrow$ (i). We use the characterisation (ii) in Theorem~\ref{abstract-unierg} with
the dense algebra of functions $\mathcal D$ from the proof of Theorem~\ref{uniFLC}, 
{\rm (ii)} $\Rightarrow$ {\rm (i)}. As explained there, it suffices to consider products $\prod_{i=1}^{k} f_{\varphi_{i}}$, $k\in\NN$, with $\varphi_{i} \in L^{0}_{b,c}(M)$ for $i\in\{1, \ldots, k\}$ having pairwise disjoint set-theoretical supports $U_{i}$. Given $P\in X_\cP$, we abbreviate
\begin{displaymath}
I_n(P):=\frac{1}{\vol (D_n)}\int_{D_n} \!\d x\,  \bigg(\prod_{i=1}^k 
f_{\varphi_i}\bigg)(xP)
\end{displaymath} 
and take a sequence $\left((y_m,P_m)\right)_{m\in\mathbb N}\subseteq T\times \cP$ 
such that $(y_mP_m)_{m\in\mathbb N}$ converges to $P$.
Then, for every $n\in\mathbb N$, the sequence $(I_n(y_mP_m))_{m\in\mathbb N}$ 
converges to $I_n(P)$ by dominated convergence. On the other hand,  uniform pattern frequencies, 
Remark~\ref{QP=QXP} and Proposition~\ref{comp-freq} imply that 
\begin{displaymath}
 	\lim_{n\to\infty} I_n(\wtilde P) = \sum_{Q\in \cF_{X_{\cP}}^{k}(U)} I(Q) \,\nu(Q) =: \mathcal{J},
\end{displaymath}
the convergence being uniform in $\wtilde P \in X_{\cP}$ and the limit $\mathcal{J}$ independent of $\wtilde P \in X_{\cP}$.  
In particular, uniformity allows the interchange of limits in 
\begin{displaymath}
\lim_{n\to\infty} I_n(P)=\lim_{n\to\infty}\lim_{m\to\infty} I_n(y_mP_m)
=\lim_{m\to\infty}\lim_{n\to\infty} I_n(y_mP_m)= \mathcal{J},
\end{displaymath} 
showing that $\lim_{n\to\infty} I_n(P)$ exists for every $P\in X_{\cP}$ and is independent of $P$.
\end{proof}

\begin{proof}[Proof of Proposition~\ref{cylinder1}]
With $1_A$ denoting the indicator function of a set $A$, the pointwise ergodic theorem 
Theorem~\ref{abstract-erg} {\rm (ii)} (together with a tempered subsequence $(D_n)_{n\in\mathbb N}$ 
of the given van Hove sequence in $T$) yields for $\mu$-a.a.~$P\in X_\cP$
\begin{displaymath}
\mu(C_\mathbold{U})=\mu(1_{C_\mathbold{U}})=\lim_{n\to\infty}\frac{1}{\vol (D_n)}\int_{D_n}{\rm d}x 
\,  
1_{C_\mathbold{U}}(xP).
\end{displaymath}
On the other hand, the indicator function of the cylinder set
$C_\mathbold{U}$ can be expressed as
\begin{displaymath}
1_{C_\mathbold{U}}=f_{1_{U_1}}\cdot\ldots\cdot f_{1_{U_k}},
\end{displaymath}
since $\mathrm{diam}(U_i)<r$. Apart from $Q$ itself, the region $(Q)_{\varepsilon}$ contains, 
up to equivalence, by hypothesis no other pattern of $\cQ_{X_{\cP}}$ of the same cardinality as $Q$. Therefore we can apply Proposition~\ref{comp-freq} 
with $\cF_{X_{\cP}}^{k}(U) = \{Q\}$, which yields for all $P\in X_{\cP}$
\begin{displaymath}
\int_{D_{n}}\!\d x\, 1_{C_\mathbold{U}}(xP) =
I(Q) \, \card\big(M_{D_n}(Q)\big) +
     \mathocal{o}\big(\vol (D_n)\big)
\end{displaymath}
as $n\to\infty$, where the integral $I(Q)$ is given by
\begin{displaymath}
I(Q)=\sum_{\pi\in \cS_k} \int_T \!\d x \, \prod_{i=1}^k 1_{U_i}(xq_{\pi(i)}).
\end{displaymath}
Therefore we conclude from Theorems~\ref{uniFLC} and~\ref{abstract-erg} that the pattern frequency $\nu(Q)$ exists for $\mu$-a.a.~$P\in X_\cP$, 
and for any such $P$ we have
\begin{displaymath}
\mu(C_\mathbold{U})=I(Q)\,\nu(Q).
\end{displaymath}

Next we show that $I(Q)=\mathrm{vol}(D_\varepsilon)$. To do so, we use the notation of 
Remark~\ref{S-transform} and introduce $T^\pi_{Q,U} := \bigcap_{i=1}^k S_{q_{\pi(i)}, U_i}$ for $\pi \in\cS_{k}$. Since each $U_{i}$ can accommodate at most one point of a pattern, we obtain 
\begin{align}
	\label{disjoint}
	D_\varepsilon &= \left\{x\in T: xQ\subseteq U\right\} \\
	&=\bigcup_{\pi\in\mathcal{S}_k} \big\{x\in T:xq_{\pi(i)}\in U_i\text{~~for all~~}i\in\{1,\ldots,k\}\big\}
	\nonumber \\
&= \bigcup_{\pi\in\mathcal{S}_k} T^\pi_{Q,U}
= \bigcup_{\pi\in\mathcal{S}_k(Q)} T^\pi_{Q,U}. \nonumber
\end{align}
 The restriction to $\mathcal{S}_k(Q)
\subseteq \mathcal S_k$ in the last equality of \eqref{disjoint} is justified, because if for some $\pi\in\cS_{k}$ there is $x\in T$ such that $xq_{\pi(i)}\in U_i$ 
for all $i$, then there must exist $x_{\pi}\in T$ such that $x_{\pi}q_{\pi(i)}=q_i$, due to our hypothesis on the smallness of 
$\varepsilon$ and the uniqueness of $Q$. Hence $\pi \in \mathcal S_k(Q)$. The representation \eqref{disjoint} also implies that $D_\varepsilon$ is open and relatively compact in $T$, compare Lemma~\ref{transporterU}.

On the other hand, since  $\pi \in\cS_{k} \setminus \cS_{k}(Q)$ does not contribute to $I(Q)$ either (by the same argument as above), we conclude
\begin{displaymath}
	I(Q)  = \sum_{\pi\in \cS_k(Q)}  \int_T \!\d x \, \prod_{i=1}^k 1_{U_i}(xq_{\pi(i)}) 
	= \sum_{\pi\in \cS_k(Q)}  \vol\big( T^\pi_{Q,U}\big).
\end{displaymath}
Thus, the desired equality $I(Q)=\mathrm{vol}(D_\varepsilon)$ follows 
if the rightmost union in \eqref{disjoint} is disjoint. To see this we take 
$\pi,\wtilde\pi\in\mathcal S_k(Q)$. By definition, there exist $x_\pi,x_{\wtilde\pi}\in T$ such that 
\begin{displaymath}
 	x_\pi q_{\pi(i)}=q_i \qquad \text{and} \qquad x_{\wtilde\pi}q_{\wtilde\pi(i)} =q_{i}
\end{displaymath}
for all $i \in \{1,\ldots,k\}$.  On account of Remark~\ref{S-transform}, this implies 
\begin{equation}
	\label{T-transform}
 	T^{\pi}_{Q,U} = T^{\mathrm{id}}_{Q,U}x_{\pi}  \qquad \text{and} \qquad 
	T^{\wtilde\pi}_{Q,U} = T^{\mathrm{id}}_{Q,U}x_{\wtilde\pi}.
\end{equation}
Assuming $T^\pi_{Q,U}\cap T^{\wtilde \pi}_{Q,U}\ne\varnothing$, we see that \eqref{T-transform} then 
ensures the existence of $y, \wtilde{y} \in T^{\mathrm{id}}_{Q,U}$ which
obey $y x_\pi= \wtilde{y} x_{\wtilde \pi}$. This implies in turn
\begin{displaymath}
y q_i=y x_\pi q_{\pi(i)}=\wtilde{y} x_{\wtilde \pi} q_{\pi(i)}= \wtilde{y} x_{\wtilde \pi} q_{\wtilde\pi\left((\wtilde\pi^{-1}\circ\pi)(i)\right)}= \wtilde{y} q_{(\wtilde\pi^{-1}\circ\pi)(i)}
\end{displaymath}
for all $i \in\{1,\ldots,k\}$. Hence, $y q_i\in U_i \cap U_{(\wtilde\pi^{-1}\circ\pi)(i)}$ for all 
$i \in\{1,\ldots,k\}$. Since $U_i\cap U_j=\varnothing$ for $i\ne j$, we infer that $\pi=\wtilde\pi$. Hence the 
rightmost union in \eqref{disjoint} is disjoint and $I(Q)=\mathrm{vol}(D_\varepsilon)$ holds. This completes the proof of the first statement of the proposition. 

To show the remaining statement of the proposition we 
assume that $T$ is Abelian and note that 
\begin{align*}
S_{ym, B_\varepsilon(ym)} &=\{x\in T:d(xym,ym)<\varepsilon\} 
	=\{x\in T: d(xm,m)<\varepsilon\}  \\
	&= S_{m, B_\varepsilon(m)}
\end{align*}
for all $y\in T$ and all $m\in M$ due to $T$-invariance of the metric. 
Hence, if the group acts also transitively on $M$, we infer 
$T^{\mathrm{id}}_{Q,U}= S_{m, B_{\varepsilon}(m)}$
for every $m\in M$. Together with \eqref{T-transform} and \eqref{disjoint} this implies
\begin{displaymath}
	D_\varepsilon=\bigcupdisjoint_{\pi\in\mathcal{S}_k(Q)} S_{m, B_{\varepsilon}(m)} x_\pi,
\end{displaymath}
and the statement follows from unimodularity of $T$ (which yields right 
invariance of the Haar measure).
\end{proof}

%%%%%%%%%%%%%%%%%%%%%%%%%%%%%%%%%%%%%%%%%%%%%%%%%%%%%%%%%%%%%%%%%%%%%%%%%%%%%
%%%%%%%%%%%%%%%%%%%%%%%%%%%%%%%%%%%%%%%%%%%%%%%%%%%%%%%%%%%%%%%%%%%%%%%%%%%%%
%
\section{Proofs of results in Section~\ref{sec:ran-col}}
\label{sec:proofs-ran-col}
%
%%%%%%%%%%%%%%%%%%%%%%%%%%%%%%%%%%%%%%%%%%%%%%%%%%%%%%%%%%%%%%%%%%%%%%%%%%%%%
%%%%%%%%%%%%%%%%%%%%%%%%%%%%%%%%%%%%%%%%%%%%%%%%%%%%%%%%%%%%%%%%%%%%%%%%%%%%%

\begin{proof}[Proof of Lemma~\ref{mhat-ass}]
	We only comment on properness of the induced action $\what\alpha$, since all other claims are evident.
 	Properness of $\alpha$ follows from Prop.~5(ii) in \cite[Chap.\ III.4.2]{Bour1}, where we choose $G=G'=T$, $\varphi = \mathrm{id}$, $X= \what{M}$, $X'=M$ and $\psi: \what{M} \rightarrow M$, $(m,a) \mapsto m$. 
\end{proof}

\begin{proof}[Proof of Proposition~\ref{cCr}]
  Compactness follows from closedness of $\cC_\mathcal{P}$ in the compact
  metrisable space $\cP_r(\what{M})$.  Let $(P_n^{(\omega_n)})_{n\in\mathbb
    N}\subseteq \mathcal C_\mathcal{P}$ be a sequence with
  $$
  	\lim_{n\to\infty}P_n^{(\omega_n)}=\what{P} \in\mathcal P_r(\what{M}).
	$$
  Let
  $P:=\pi(\what{P})\subseteq M$. We show that $P\in\mathcal P$ and that
  $\what{P}$ is a coloured point set, which implies $\what{P} \in
  \cC_\mathcal{P}$.

  Continuity of the projection $\pi$ yields $\lim_{n\to\infty} P_{n}
  =P$. Therefore we have $P\in\cP$ by closedness of $\cP$. Now assume that
  $\what{p}_{1} := (p,a_{1})\in \what{P}$ and $\what{p}_{2} := (p,a_{2}) \in
  \what{P}$, where $p\in P$ and $a_{1},a_{2} \in\AA$. 
  Thus, there exist two sequences $(\what{p}^{\,n}_{j})_{j\in\NN}$, $j=1,2$, such that 
  $\what{p}^{\,n}_{j} \in P_n^{(\omega_n)}$ for all $n\in\NN$ and 
  $\lim_{n\to\infty} \what{d}(\what{p}^{\,n}_{j}, \what{p}_{j}) = 0$ for both $j=1,2$. 
  Continuity of $\pi$ yields 
  $\lim_{n\to\infty} d(p^{n}_{1}, p) = 0 = \lim_{n\to\infty} d(p^{n}_{2}, p)$
  This implies   $d(p_{1}^{n}, p_{2}^{n}) < r$ for finally all $n\in\NN$. 
	Uniform discreteness of $P_{n}$ then
  yields $p_{1}^{n}= p_{2}^{n} =: p^{n}$, and we must have $a_{1}^{n} =
  \omega_{n}(p^{n}) = a_{2}^{n}$ for
  finally all $n\in\NN$. This shows $a_{1}=a_{2}$.
\end{proof}

\begin{proof}[Proof of Lemma~\ref{col-T-inv}]
  Let $Y:= \big\{ xP^{(\omega)}: x\in T, \, P^{(\omega)} \in
  \cC_{\cP}\big\}$. Since $\cC_{\cP} \subseteq \hXcP$ and $\hXcP$ is
  $T$-invariant and closed, we deduce $\overline{Y} \subseteq \hXcP$. 

  To prove the converse inclusion, let $P^{(\omega)} \in \hXcP$
  arbitrary. This means $P \in X_{\cP}$, so there exists a sequence
  $(P_{n})_{n\in\NN}$ in $\cP$ converging to $P$. By choosing appropriate
  $\omega_{n} \in \Omega_{P_{n}}$, we obtain a
  sequence $(P_{n}^{(\omega_{n})})_{n\in\NN}$ in $\cC_{\cP} \subseteq Y$
  which converges to $P^{\omega}$. Hence, $P^{(\omega)} \in \overline{Y}$.

  Continuity of the group action $\alpha_{\hXcP}$ follows from
  continuity of action $\what{\alpha}$ on $\what{M}$.
\end{proof}

\begin{proof}[Proof of Lemma~\ref{IAD}]
  We give a detailed proof for the first example only. 
  The proof for the second example follows along the same lines, where
  $T$-stationarity ensures $T$-covariance, the compactly supported strong mixing coefficient 
  ensures independence at a distance and where the continuous realisations 
  $\xi^{(\sigma)}(\cdot)$ ensure $C$-compatibility and hence $M$-compatibility.
  
  For the first example, $T$-covariance and independence at a distance are clear. 
  It remains to verify $C$-compatibility, from which $M$-compatibility follows. 
  First, we construct a $\|\cdot\|_{\infty}$-dense subset $\mathcal{D}$ of
  $C(\hXcP)$. For $\varphi\in
  C_c(M)$ and $\psi\in C_c(\AA)$ define $f_{\varphi,\psi}:\hXcP\to\mathbb R$ by
  \begin{equation}\label{fphipsi}
    f_{\varphi,\psi}(P^{(\omega)}):=\sum_{p\in P}
    \varphi(p)\cdot\psi\big(\omega(p)\big), \qquad P^{(\omega)} \in\hXcP.
  \end{equation}
  Continuity of $f_{\varphi,\psi}$ is obvious from the definition of the vague topology. 
  We also introduce the constant function $1\in C(\hXcP)$ equal to one and the set  
  \begin{equation}
    \label{D-null}
    \mathcal{D}_{0}:=\Big\{f_{\varphi,\psi}:\varphi\in C_c(M) \text{~with~}
    \diam\left(\supp(\varphi)\right)<r/2, \;
    \psi\in C_c(\AA)\Big\} \cup\{1\},
  \end{equation}
  which separates points in $\hXcP$. The Stone-Weierstra\ss{} theorem
  \cite[Prob.~7R]{Kel75} then assures that the algebra
  $\mathcal{D}:=\mathrm{alg}(\mathcal{D}_{0})$ generated by $\mathcal{D}_{0}$
  is dense in $C(\hXcP)$ with respect to the supremum norm.

  Since
  \begin{displaymath}
    E_{f}(P) -E_{f}(P') = \int_{\Omega_{P}}\!\d\PP_{P}(\omega)
    \int_{\Omega_{P'}}\!\d\PP_{P'}(\sigma) \, [f(P^{(\omega)}) - f(P'^{(\sigma)})]
  \end{displaymath}
  for all $P,P'\in X_{\cP}$ and since the algebra $\mathcal{D}$ is uniformly
  dense in $C(\hXcP)$, it suffices to prove continuity of $E_{f}$ for
  functions $f$ of the form $g_{k} := \prod_{i=1}^{k}
  f_{\varphi_{i},\psi_{i}}$, where $k\in\NN$ and $f_{\varphi_{i},\psi_{i}}
  \in\mathcal{D}_{0}$ for all $i\in\{1,\ldots, k\}$. Furthermore, since $X_\cP$
  is metrisable, it suffices to show sequential continuity of $E_{g_k}$.

  Fix $P\in X_\cP$ and take a sequence $(P_n)_{n\in\mathbb N}\subseteq X_\cP$
  which converges to $P$. Define the compact set $V:=\bigcup_{i=1}^k \supp(\varphi_i)$ and 
  the (finite) pattern 
  $Q:=P\cap \mathring V$. Then the pattern $Q$ and $
	(\mathring V)^c$ have a positive distance $\delta_0:=d(Q, (\mathring V)^c)>0$. 
For arbitrary fixed $\delta\in]0,\min(\delta_0,r)[$, we find by Lemma~\ref{con-thick} 
  an $N=N(\delta)$ such that we have for all $n\ge N$ the inclusions
  \begin{equation}\label{close}
  P_n\cap V \subseteq (P)_\delta, \qquad P\cap V\subseteq (P_n)_\delta.
  \end{equation}
  For $n\ge N$ we consider the finite patterns
    \begin{displaymath}
    Q_n :=  \Big\{ p \in P_n : \; 
    \exists\; q\in Q \text{~with~} d(p,q) <\delta \Big\}\subseteq \mathring V. 
  \end{displaymath}
  By \eqref{close}, there exists a bijection $h_n: Q \rightarrow Q_n$ with $d\big(q,h_n(q)\big) 
  <\delta$ for all $q \in Q$ and for all $n\ge N$. Thus, we get
  \begin{align}
    \label{PQ-close}
    E_{g_{k}}(P_n) &= \sum_{(q_{1},\ldots,q_{k})  \in Q^k} \bigg(
      \prod_{i=1}^{k} \varphi_{i}\big(h_n(q_{i})\big) \bigg) \int_{\Omega_{P_n}}
    \!\d\PP(\sigma)\, \prod_{i=1}^{k} \psi_{i}\Big(\sigma\big(h_n(q_{i})\big)\Big)
    \\ 
    &= \sum_{(q_{1},\ldots,q_{k}) \in Q^k} \bigg(
      \prod_{i=1}^{k} \varphi_{i}\big(h_n(q_{i})\big) \bigg) \int_{\Omega_{P}}
    \!\d\PP(\omega)\, \prod_{i=1}^{k} \psi_{i}\big(\omega(q_{i})\big), \nonumber
  \end{align}
  where the last equality follows from the fact that all random variables are
  independently and identically distributed. This implies for all $n\ge N$ the estimate
  \begin{align*}
    |E_{g_{k}}(P) -E_{g_{k}}(P_n)| \le \sum_{(q_{1},\ldots,q_{k}) \in Q^k} & \left| \prod_{i=1}^{k}
      \varphi_{i} (q_{i}) 
      - \prod_{i=1}^{k} \varphi_{i}\big(h_n(q_{i})\big) \right| \nonumber\ \\
   & \qquad \times \int_{\Omega_{P}} \!\d\PP(\omega)\, \prod_{i=1}^{k}
    \big| \psi_{i}\big(\omega(q_{i})\big) \big|  \nonumber\\
     \le\bigg(\prod_{i=1}^{k} \|\psi_{i}\|_\infty\bigg) \sum_{(q_{1},\ldots,q_{k}) \in Q^k} &\left| \prod_{i=1}^{k}
      \varphi_{i} (q_{i}) 
      - \prod_{i=1}^{k} \varphi_{i}\big(h_n(q_{i})\big) \right|.
  \end{align*}
  Since the functions $\varphi_{i}$ are continuous with compact support, we can
  make this difference as small as we want uniformly in $n\ge N$, by 
  choosing $\delta$ sufficiently close to zero. 
\end{proof}

\begin{proof}[Proof of Theorem~\ref{lem:hoflemma}]
  The map $Y_n:\Omega_P\to\mathbb R$ is continuous (hence measurable) for
  every $n\in \mathbb N$, as can be seen by applying Lebesgue's dominated
  convergence theorem.

  Below we prove \eqref{hofstatement} for random variables $Y_{n}$
  corresponding to functions $f$ in the $\|\cdot\|_{\infty}$-dense subalgebra
  $\mathcal{D} \subseteq C(\hXcP)$, which was introduced below Eq.\
  \eqref{D-null}. This and an $\varepsilon/3$-argument establish the lemma for
  all $f\in C(\hXcP)$ because, given an approximating sequence
  $(f_k)_{k\in\mathbb N}\subseteq \mathcal{D}$, we have
  $\big|Y_n^{(k)}(\omega)-Y_n(\omega)\big|\le \|f_k-f\|_\infty$ uniformly in
  $n$ and in $\omega$. Here, $Y_n^{(k)}$ denotes the random variable
  \eqref{Yn} corresponding to $f_{k}$.

  Thus, it suffices to prove \eqref{hofstatement} for random variables $Y_{n}$
  corresponding to functions $f$ of the form $f=f_{\varphi_1,\psi_1}
  \cdot\ldots\cdot f_{\varphi_k,\psi_k}$, where $k\in\mathbb N$ and
  $f_{\varphi_j,\psi_j}\in \mathcal{D}_{0}$ for $j\in\{1,\ldots,k\}$. To do so we
  fix $P\in X_{\cP}$ and $\omega\in\Omega_P$ arbitrary and set
  $U_j:=\mathrm{int}(\supp(\varphi_j))$, which is relatively compact for $j\in\{1,\ldots,k\}$.
  Then we can apply Proposition~\ref{comp} (with $\what{M}$ playing the r\^ole of $M$ there) and obtain
  \begin{equation}\label{explicit}
    \int_{D_n}{\!\d}x \; \bigg(\prod_{i=1}^k
    f_{\varphi_i,\psi_i}\bigg)(xP^{(\omega)}) =  \sum_{\pi \in\cS_{k}}  \sum_{Q \in\cQ_{P}^{k}(U;D_{n})} 	
    	\!\!I^{\pi}(Q) \, Z^{\pi}_{Q}(\omega) +
     \mathocal{o}\big(\vol (D_n)\big),
  \end{equation}
  asymptotically as $n\to\infty$, where the $\mathocal{o}(\vol (D_n))$-term can be chosen uniformly 
  in  $P\in X_{\cP}$ and $\omega\in\Omega_P$. 
  In \eqref{explicit} we have used the notation of Proposition~\ref{comp} and 
  Definition~\ref{pat-col}, except that we have singled out the sum over permutations $\pi$ from 
  the integral \eqref{I-def}, as well as the part involving the random variables
\begin{displaymath}
    \Omega_{P} \ni \omega \mapsto Z^{\pi}_{Q}(\omega) := \prod_{i=1}^k
    \psi_{i}\bigl( \omega(q_{\pi(i)})\bigr),
  \end{displaymath}
which amounts to setting 
  \begin{displaymath}
    I^{\pi}(Q) :=  \int_{T}\!\d x\, \prod_{i=1}^k \varphi_i(xq_{\pi(i)}).
  \end{displaymath}
   The lemma will now follow from \eqref{explicit} and the relation
  \begin{equation}
    \label{lln}
    \lim_{n\to\infty} \frac{1}{\vol(D_{n})}\, \sum_{Q \in\cQ_{P}^{k}(U;D_{n})} I^{\pi}(Q) \,
     \left[ Z^{\pi}_{Q}(\omega)
      - \int_{\Omega_{P}}\!\d\PP_{P}(\eta)\, Z^{\pi}_{Q}(\eta) \right] =0
  \end{equation}
  for $\PP_{P}$-a.a.\ $\omega\in\Omega_{P}$, every $P \in $ and every permutation $\pi \in\cS_{k}$.
  
  If the set $\cQ_{P}^{k}(U) :=\cQ_{P}^{k}(U;T)$ is finite, then \eqref{lln} follows from 
  $\mathrm{vol}(D_n)\to\infty$ as $n\to\infty$.
  Hence we assume in the remainder that the set $\cQ_{P}^{k}(U)$ is infinite. 
  The above relation \eqref{lln} will then follow
  from the strong law of large numbers, as we now show. We note, first, that the variances
  $\mathrm{Var}(Z^{\pi}_{Q}) \le \prod_{i=1}^{k}
  \|\psi_{i}\|_{\infty}^{2}$ are bounded uniformly in $Q$ (and $\pi$). 
  Second, the cardinality of the finite set $\cQ_{P}^{k}(U;D_{n})$  grows at
  most with $\vol(D_{n})$. This can be seen from the relative compactness of $U=\bigcup_{i=1}^{k} U_{i}$, 
  Lemma~\ref{pattern-count} and Proposition~\ref{standardest-vol},
  which require both unimodularity.
  Third, we show that the coefficients $I^{\pi}(Q)$ are uniformly bounded in
  $Q \in \cQ_{P}^{k}(U)$ (and $\pi \in\cS_{k}$). This will follow from \eqref{IQ-start}, which yields the estimate
  \begin{displaymath}
    \label{coeffbound}
    |I^{\pi}(Q)| 
    \le  \vol\big( S(Q)\big) \, \prod_{i=1}^{k} \| \varphi_{i}\|_{\infty}, 
  \end{displaymath}
	the inclusion \eqref{S-reduce-1}, compactness of $L_{U}:= S_{\overline{U},\, \overline{U}}$ by 
	Lemma~\ref{transporterK}
	and right invariance $\vol(L_{U} y) = \vol (L_{U}) < \infty$ of the Haar measure
  on the unimodular group $T$.  
  
  Having these three properties in mind, the desired relation \eqref{lln}
  follows from the strong law of large numbers and Kolmogorov's criterion
  \cite[Thm.~14.5]{Bau01}, provided we know that the family $\big(I^{\pi}(Q) Z^{\pi}_{Q}\big)_{Q \in \cQ^{k}_{P}(U)}$
  consists of pairwise independent random variables. 
 
    If pairwise independence happens not to be the case, then we argue below that the index set
    $\cQ^{k}_{P}(U)$ can be partitioned into a \emph{finite} number $J$
    of mutually disjoint subsets
  \begin{equation}
    \label{partition}
    \cQ^{k}_{P}(U) = \bigcupdisjoint_{j=1}^{J} F_{j}
  \end{equation}
  such that for each $j\in\{1,\ldots,J\}$ the subfamily $\big(I^{\pi}(Q) Z^{\pi}_{Q}\big)_{Q \in F_{j}}$ consists 
  of pairwise independent random variables. Assuming this decomposition for the
  time being, we rewrite \eqref{lln} as
  \begin{equation} \label{partition2}
    \lim_{n\to\infty}\sum_{j=1}^{J} \frac{\card
      (F_{j}^{(n)})}{\vol(D_{n})} \, \mathcal{Z}_{j}^{(n)}(\omega) =0
    \qquad\quad\text{for $\PP_{P}$-almost all $\omega\in\Omega_{P}$,}
  \end{equation}
  where $F_{j}^{(n)} := F_{j} \cap \cQ_{P}^{k}(U;D_{n})$ and 
  \begin{displaymath}
    \mathcal{Z}_{j}^{(n)}(\omega) := \frac{1}{\card(F_{j}^{(n)})}\,
    \sum_{Q \in F_{j}^{(n)}}  I^{\pi}(Q) \, \left[ Z^{\pi}_{Q}(\omega) 
      - \int_{\Omega_{P}}\!\d\PP_{P}(\eta)\, Z^{\pi}_{Q}(\eta) \right]
  \end{displaymath}
  for $j\in\{1,\ldots,J\}$. But \eqref{partition2} is indeed true. This follows from
  $\mathrm{vol}(D_n)\to\infty$ as $n\to\infty$ for those $j\in\{1,\ldots,J\}$ such 
  that $F_j$ is finite, and from the $\PP_{P}$-almost
  sure relation $\lim_{n\to\infty} \mathcal{Z}_{j}^{(n)} =0$ for those $j\in\{1,\ldots,J\}$ such 
  that $F_j$ is infinite, thanks to pairwise independence by the
  strong law of large numbers and Kolmogorov's criterion.

  It remains to verify the existence of the partition \eqref{partition}. This
  may be seen by a graph-colouring argument: construct a graph $\mathcal{T}$
  with infinite vertex set $\cQ^{k}_{P}(U)$. Two vertices $Q$ and
  $Q'$ of $\mathcal{T}$ are joined by an edge, if and only if $Z^{\pi}_{Q}$ and
  $Z^{\pi}_{Q'}$ are $\PP_P$-dependent. Clearly, a vertex colouring of $\mathcal{T}$ 
  (with finitely many colours and with adjacent vertices having different colours) provides 
  an example for the partition that we are seeking.   
  Due to uniform discreteness of $P$, independence at a
  distance of $\PP_{P}$ and because $U$ is
  contained in some compact set in $M$, we infer that the degree of any vertex
  in $\mathcal{T}$ is bounded by some number $d_{\mathcal{T}\!,\mathrm{max}} <
  \infty$. Thus, the vertex-colouring theorem \cite{Die05} ensures the
  existence of such a colouring with $J \le 1+
  d_{\mathcal{T}\!,\mathrm{max}}$ different colours.
\end{proof}

\begin{proof}[Proof of Theorem~\ref{theo:perg}]
  First, we prove the existence of a unique $T$-invariant Borel probability
  measure $\hatmu$ on $\hXcP$ which obeys \itemref{fubini}.
  
  Thanks to $M$-compatibility, Assumption~\ref{p-ass-wcomp}, the integral
  $$
  	I(f) := \int_{X_{\cP}}\!\d\mu(P)\, E_{f}(P)
	$$ 
	is well-defined and finite for
  every $f\in C(\hXcP)$. Moreover, the map $I: C(\hXcP) \rightarrow \RR, f
  \mapsto I(f)$ is a positive, bounded linear functional which is also
  normalised, $I(1)=1$, and $T$-invariant because of \eqref{col-ps-trans},
  $T$-covariance of $\PP_{P}$ and $T$-invariance of $\mu$. By the Riesz-Markov
  theorem there exists a unique Borel probability measure $\hatmu$ on $\hXcP$
  such that
  \begin{equation}
    \label{hatmu-rep}
    \hatmu(f) = \int_{X_{\cP}} \!\d\mu(P) \int_{\Omega_{P}}\!\d\PP_{P}(\omega) \;
    f(P^{(\omega)}) 
  \end{equation}
  for all $f \in C(\hXcP)$. Since $\hXcP$ is a compact metric space and
  $\hatmu$ is a Borel measure, the continuous functions $C(\hXcP)$ lie dense
  in $L^{1}(\hXcP,\hatmu)$ w.r.t.\ $\|\cdot\|_{1}$. Thus,
  given $f\in L^{1}(\hXcP,\hatmu)$ there exists a sequence
  $(f_{k})_{k\in\NN} \subseteq C(\hXcP)$ which converges pointwise and in
  $\|\cdot\|_{1}$-sense towards $f$. This and dominated convergence yield for
  all $f \in L^{\infty}(\hXcP,\hatmu)$ measurability of the map 
  $E_{f}: X_{\cP} \rightarrow \RR \cup \{\pm\infty\}$ and that
  \eqref{hatmu-rep} holds. Finally, these conclusions hold also for 
  $f \in L^{1}(\hXcP,\hatmu)$ by decomposing $f$
  into its positive and negative part and using monotone convergence for a
  sequence of $L^{\infty}$-approximants.

  In what remains we prove ergodicity of the $T$-invariant probability measure
  $\hatmu$. The additional statements about exceptional sets will be obtained
  along the way.  Fix $f \in C(\hXcP)$ arbitrary. On the one hand, the Ergodic
  Theorem~\ref{abstract-erg} for $\cQ=\hXcP$ provides the existence of
  $f^{\star} \in L^{1}(\hXcP,\hatmu)$ and of a $\hatmu$-null set
  $\what{N} \subseteq \hXcP$ such that
  \begin{equation}
    \label{total-erg}
    \lim_{n\to\infty} \frac{1}{\vol(D_{n})}  \int_{D_n}\!\d x\;
    f(xP^{(\omega)})  = f^{\star}(P^{(\omega)})
  \end{equation}
  for all $P^{(\omega)} \in \hXcP \setminus \what{N}$. On the other hand, we
  apply the Ergodic Theorem~\ref{abstract-erg} for $\cQ=X_{\cP}$ to the
  function $E_f\in L^{\infty}(X_{\cP},\mu)$ and combine it with
  Theorem~\ref{lem:hoflemma} (which requires unimodularity of $T$).  
  This yields the existence of a set
  $\wtilde{X}\subseteq X_{\cP}$ of full $\mu$-measure and, for every
  $P\in\wtilde{X}$, of a set $\wtilde{\Omega}_P\subseteq \Omega_P$ of full
  $\mathbb P_P$-measure such that the equality
  \begin{align}\label{fibred-erg}
    \lim_{n\to\infty} \frac{1}{\vol(D_{n})}  \int_{D_n}\!\d x\,
    f(xP^{(\omega)})  
    &=  \int_{X_{\cP}} \!\d\mu(Q) \int_{\Omega_{Q}}\!\d\PP_{Q}(\sigma) \,
    f(Q^{(\sigma)})  \\
    &= \hatmu(f)\nonumber
  \end{align}
  holds for all $P\in\wtilde{X}$ and for all $\omega\in\wtilde{\Omega}_P$.  In
  the uniquely ergodic case we rely on $C$-compatibility
  $E_{f} \in C(X_{\cP})$ and apply the Ergodic Theorem~\ref{abstract-unierg}
  instead. This gives \eqref{fibred-erg} with $\wtilde{X}=X_{\cP}$ (and
  without requiring temperedness for the F\o{}lner sequence). But
  \begin{align*}
    \int_{\hXcP} \!\d\hatmu & (P^{(\omega)}) \, |f^{\star}(P^{(\omega)}) -
    \hatmu(f)|  \nonumber\\ 
    & = \int_{\hXcP} \!\d\hatmu(P^{(\omega)}) \, 1_{\hXcP \setminus\what{N}}(P^{(\omega)}) \,
    |f^{\star}(P^{(\omega)}) -  \hatmu(f)|     \nonumber\\ 
    &=   \int_{X_{\cP}} \!\d\mu(P) \int_{\Omega_{P}}\!\d\PP_{P}(\omega) \,
     1_{\hXcP \setminus\what{N}}(P^{(\omega)}) \, |f^{\star}(P^{(\omega)}) -
     \hatmu(f)|  \nonumber \\ 
    &= 0
  \end{align*}
  on account of \eqref{fubini-eq}, \eqref{total-erg} and \eqref{fibred-erg},
  showing $\hatmu$-a.e. $f^{\star} = \hatmu(f)$ for all $f\in C(\hXcP)$. The
  implication ~(iii) $\Rightarrow$ (i)~ in the Ergodic
  Theorem~\ref{abstract-erg} for $\cQ=\hXcP$ now completes the proof.
\end{proof}

\begin{proof}[Proof of Proposition~\ref{col-cyl-prop}]
Let $(D_n)_{n \in\mathbb N}$ be a tempered subsequence of  a F{\o}lner sequence in $T$. 
By Theorem~\ref{theo:perg}  we have for $\mu$-a.a.\ $P \in X_\cP$ and for $\PP_{P}$-a.a.\ $\omega\in\Omega_{P}$ that 
 \begin{displaymath}
 \begin{split} 
     \hatmu(C_\mathbold{U}^\mathbold{A})&= \int_{\hXcP} \!\d\hatmu(Q^{(\sigma)}) \; 1_{C_\mathbold{U}^
\mathbold{A}}(Q^{(\sigma)}) 
     = \int_{X_{\cP}} \!\d\mu(Q) \int_{\Omega_{Q}}\!\d\PP_{Q}(\sigma) \;
      1_{C_\mathbold{U}^\mathbold{A}}(Q^{(\sigma)})\\
&=\lim_{n\to\infty} \frac{1}{\vol(D_{n})} \int_{D_n}\!\d x\; 
1_{C_\mathbold{U}^\mathbold{A}}(xP^{(\omega)}).
\end{split}
\end{displaymath}
Since the coloured cylinder set $C_\mathbold{U}^\mathbold{A}$ contains precisely all those 
coloured point sets which possess  
exactly one point in each of the $U_{i}$ (thanks to $\mathrm{diam}(U_i)<r$), and with corresponding colour value in $A_i$ for $i\in\{1,\ldots,k\}$, we can express its 
indicator function as
\begin{displaymath}
1_{C_\mathbold{U}^\mathbold{A}}=f_{1_{U_1}, 1_{A_{1}}}\cdot\ldots\cdot f_{1_{U_k}, 1_{A_{k}}}.
\end{displaymath}
Thus, we conclude from \eqref{explicit} -- which, according to the hypotheses of Proposition~\ref{comp}, is also valid for indicator functions $\varphi_{i}= 1_{U_{i}}$, $\psi_{i}=1_{A_{i}}$, $i\in\{1,\ldots, k\}$, of open, relatively compact sets -- 
that for every $P^{(\omega)}\in \hXcP$ the equality
 \begin{displaymath}
    \int_{D_n}{\d}x \,1_{C_\mathbold{U}^\mathbold{A}}(xP^{(\omega)}) = \sum_{Q\in \cQ^k_{P}(U;D_{n})} 
    \sum_{\pi\in\cS_k} I^{\pi}(Q) \, Z^{\pi}_{Q}(\omega) +\mathocal{o}\big(\vol (D_n)\big),
\end{displaymath}
holds asymptotically as $n\to\infty$. Here we used the notation introduced in \eqref{explicit}. 
Likewise, the law of large numbers \eqref{lln} continues to hold for 
$\varphi_{i}= 1_{U_{i}}$, $\psi_{i}=1_{A_{i}}$, $i\in\{1,\ldots, k\}$. This amounts to 
  \begin{displaymath}
    \lim_{n\to\infty} \frac{1}{\vol(D_{n})}\, \sum_{Q\in \cQ^k_{P}(U;D_{n})} 
    I^{\pi}(Q) \, \left[ Z_Q^{\pi}(\omega)
      - \PP(A_1)\cdot\ldots\cdot\PP(A_k) \right] =0
  \end{displaymath}
  for $\PP_{P}$-a.a.\ $\omega\in\Omega_{P}$, for every $P \in X_{\cP}$ and every permutation $\pi\in\cS_k$, because the expectation of $Z^{\pi}_{Q}$ factorises due to the product structure of $\PP_{P}$ and disjointness of the $U_{i}$. Now we benefit from $\PP_{P}$ being a product of identical factors and summarise the arguments so far as
  \begin{equation}
  	\label{cua-tmp}
 		\hatmu(C_\mathbold{U}^\mathbold{A}) = \PP(A_1)\cdot\ldots\cdot\PP(A_k) \;
		\lim_{n\to\infty} \frac{1}{\vol(D_{n})}\, \sum_{Q\in \cQ^k_{P}(U;D_{n})} I(Q), 
	\end{equation}
 where $I(Q) := \sum_{\pi\in\cS_k} I^{\pi}(Q) = \sum_{\pi\in\cS_k} \int_{T}\d x\, \prod_{i=1}^{k} 1_{U_{i}}(xq_{\pi(i)})$. Eq.\ \eqref{cua-tmp} holds for $\mu$-a.a.\ $P\in X_{\cP}$. Since $f_{1_{U_1}}\cdot\ldots\cdot f_{1_{U_k}} = 1_{C_\mathbold{U}}$, Proposition~\ref{comp} yields 
\begin{displaymath}
	\label{cua-tmp2}
 \lim_{n\to\infty} \frac{1}{\vol(D_{n})}\, \sum_{Q\in\cQ^{k}_{P}(U;D_{n})} I(Q) 
      =\lim_{n\to\infty} \frac{1}{\vol(D_{n})}\, \int_{D_n}{\d}x \,1_{C_\mathbold{U}}(xP) =  
      \mu(C_\mathbold{U}),
\end{displaymath}
where the last equality holds for $\mu$-a.a.\ $P \in X_\cP$ as a consequence of the Pointwise Ergodic 
Theorem~\ref{abstract-erg} applied to $\mu$ (Cor.~\ref{erg-cor}). The claim then follows together with \eqref{cua-tmp}.
\end{proof}

%%%%%%%%%%%%%%%%%%%%%%%%%%%%%%%%%%%%%%%%%%%%%%%%%%%%%%%%%%%%%%%%%%%%%%%%%%%%%
%%%%%%%%%%%%%%%%%%%%%%%%%%%%%%%%%%%%%%%%%%%%%%%%%%%%%%%%%%%%%%%%%%%%%%%%%%%%%
%
\section{Proofs of results in Section~\ref{secdyn}}
\label{sec:proofs-dyn}
%
%%%%%%%%%%%%%%%%%%%%%%%%%%%%%%%%%%%%%%%%%%%%%%%%%%%%%%%%%%%%%%%%%%%%%%%%%%%%%
%%%%%%%%%%%%%%%%%%%%%%%%%%%%%%%%%%%%%%%%%%%%%%%%%%%%%%%%%%%%%%%%%%%%%%%%%%%%%

\begin{proof}[Proof of Lemma~\ref{alpha-trans}] \quad \itemref{proper-trans}\quad 
Firstly, properness of the action $\alpha_{\VV}$ of $T$ on $\VV$ implies that the action $\alpha_{\VV \times\VV}$ of $T$ on $\VV \times\VV$, defined by $\alpha_{\VV \times\VV}(x, (v,w)) := x(v,w) := (xv,xw)$, is also proper. This follows from Prop.~5(ii) in \cite[Chap.\ III.4.2]{Bour1}, where we choose $G=G'=T$, $\varphi = \mathrm{id}$, $X= \VV \times\VV$, $X'=\VV$ and $\psi: \VV \times\VV \rightarrow \VV$, $(v,w) \mapsto v$.
Secondly, properness of the action $\alpha_{\VV \times\VV}$ of $T$ on $\VV \times\VV$  implies that the action $\alpha$ of $T$ on $M$ is proper. This follows from Prop.~5(i) in \cite[Chap.\ III.4.2]{Bour1}, where we choose $G=G'=T$, $\varphi = \mathrm{id}$, $X= \VV \times\VV$, $X'=M$ and $\psi: \VV \times\VV \rightarrow$, $(v,w) \mapsto \{v,w\}$. Here, the map $\psi$ is required to be continuous, onto and proper. While the first two properties are instantly clear, the third one follows from Prop.~2(c) in \cite[Chap.\ III.4.1]{Bour1}: setting there $X= \VV \times\VV$ and $K = \{e,\pi\}$, the permutation group of 2 objects which acts on $(v,w) \in \VV \times\VV$ according to $e(v,w) := (v,w)$ and $\pi (v,w) :=(w,v)$, we recognise $\psi$ as the canonical map $X \rightarrow X/K$.  
  
	\itemref{free-trans}\quad 
	Let $m_{v,w}\in M$ and $x\in T$ be given such that $xm_{v,w}=m_{v,w}$.
	This means that $\{(xv,xw),(xw,xv)\}=\{(v,w), (w,v)\}$. If $(xv,xw)=(v,w)$,
	we have $xv=v$, and freeness on $\VV$ implies $x=e$. Otherwise, we
	have $(xv,xw)=(w,v)$, implying $w=x(xw)=x^2w$. Freeness on $\VV$
	yields $x^2=e$, from which $x=e$ follows by assumption.
\end{proof}

\begin{proof}[Proof of Theorem~\ref{uniFLCgraph}]
Theorem~\ref{uniFLC} gives a characterisation of (unique) ergodicity in terms of uniform
\emph{pattern} frequencies. In order to prove Theorem~\ref{uniFLCgraph}, it suffices to show 
that the frequency of every pattern of $X_\mathcal{G}$ can be expressed in terms of frequencies 
of certain patches from $X_\mathcal{G}$.

Indeed, for every pattern $Q$ of $\cG$ which is not a patch there exists a uniquely
determined minimal patch $H$ of $\cG$ by ``adding the missing vertices on the diagonal''.
Then, the pattern $Q$ occurs if and only if the patch $H$ occurs.
\end{proof}

%%%%%%%%%%%%%%%%%%%%%%%%%%%%%%%%%%%%%%%%%%%%%%%%%%%%%%%%%%%%%%%%%%%%%%%%%
%%%%%%%%%%%%%%%%%%%%%%%%%%%%%%%%%%%%%%%%%%%%%%%%%%%%%%%%%%%%%%%%%%%%%%%%%

\section*{Acknowledgements}

We are grateful to Michael Baake and Daniel Lenz for stimulating discussions. 
We also thank Jean Bellissard and the anonymous referee for useful comments.
This work was partly supported by the German Research Council (DFG),
within CRC/Tr~12 (PM) and CRC~701 (CR).

%%%%%%%%%%%%%%%%%%%%%%%%%%%%%%%%%%%%%%%%%%%%%%%%%%%%%%%%%%%%%%%%%%%%%%%%%
%%%%%%%%%%%%%%%%%%%%%%%%%%%%%%%%%%%%%%%%%%%%%%%%%%%%%%%%%%%%%%%%%%%%%%%%%

\end{document}